\documentclass[fleqn]{article}
\usepackage[utf8]{inputenc}
\usepackage{latexsym}
\usepackage{amssymb}
\usepackage{amsmath}
\usepackage{amsfonts}
\usepackage{indentfirst}
\usepackage[a4paper, total={6in, 10in}]{geometry}
\newcommand{\newcontent}{\vspace{.08in}\noindent}

\usepackage{setspace}
\usepackage{multirow}
\usepackage{array}
\usepackage{booktabs}
\allowdisplaybreaks
\usepackage{graphicx}
\usepackage{caption}
\usepackage{subfigure}
\newtheorem{remark}{Remark}

\newtheorem{lemma}{Lemma}
\newtheorem{theorem}{Theorem}

\newtheorem{assumption}{Assumption}
\begin{document}

\title{\bf Doubly robust estimation of average treatment effect revisited
\footnote{Corresponding author: Lixing Zhu. E-mail addresses: lzhu@hkbu.edu.hk (L. Zhu). The research was supported by a grant from the University Grants Council of Hong Kong, Hong Kong, China. }}
\author{Keli Guo${\rm ^1}$, Chuyun Ye${\rm ^3}$, Jun Fan${\rm ^1}$ and Lixing Zhu${\rm ^{1,2}}$ \\~\\
{\small {\small {\it ${\rm ^1}$ Department of Mathematics, Hong Kong Baptist University, Hong Kong} }}\\
{\small {\small {\it ${\rm ^2}$ Center for Statistics and Data Science, Beijing Normal University, Zhuhai, China} }}\\
{\small {\small {\it ${\rm ^3}$ School of Statistics, Beijing Normal University, Beijing, China} }}
}
\date{}
\maketitle

\begin{abstract}
The research described herewith is to re-visit the classical doubly robust estimation of average treatment effect by conducting a systematic study on the comparisons, in the sense of asymptotic efficiency, among all possible combinations of the estimated propensity score and outcome regression. To this end, we consider all nine combinations under, respectively, parametric, nonparametric and semiparametric structures. The comparisons provide useful information on when and how to efficiently utilize the model structures in practice. Further, when there is model-misspecification, either propensity score or outcome regression, we also give the corresponding  comparisons. Three phenomena are observed. Firstly, when all models are correctly specified, any combination can achieve the same semiparametric efficiency bound, which coincides with the existing results of some combinations. Secondly, when the  propensity score is correctly modeled and estimated, but the outcome regression is misspecified parametrically or semiparametrically, the asymptotic variance is always larger than or equal to the semiparametric efficiency bound. Thirdly, in contrast, when the  propensity score is misspecified parametrically or semiparametrically, while the outcome regression is correctly modeled and estimated, the asymptotic variance is not necessarily larger than the semiparametric efficiency bound. In some cases, the ``super-efficiency" phenomenon occurs. We also conduct a small numerical study.

\newcontent
{\bf Keywords:} Average treatment effect; Doubly robust estimation; Misspecification; Semiparametric efficiency bound
\end{abstract}

\section{Introduction}
Estimating the average treatment effect (ATE) is an important issue in many fields including social science and health science. See \cite{pspara}.  There are two basic methodologies: inverse propensity score-based estimation (PS) and outcome regression-based estimation (OR). But the former requires correctly postulated propensity model and the latter needs correctly postulated regression model. To avoid these misspecifications, as a very promising method, the doubly robust estimation (DR) has been well studied to become an almost matured field. See, \cite{dr} and \cite{comment}. As well known, the most commonly used method is parametric modeling, see \cite{drpara+pspara} for example. As long as one of the models in DR is correctly specified, the estimation can be consistent. Alternatively, to avoid model mis-specification, nonparametric modeling is also applied,  see \cite{psnonpara}. Later, a compromise between parametric and nonparametric estimation rises in the context of missing data to give a semiparametic estimation, see a relevant reference \cite{semidr}.

In this paper,  we focus on investigating the estimation efficiencies of all possible combinations of PS and OR estimator obtained by respectively using parametric, semiparametric and nonparametric estimations of both $PS$ and $OR$ model. As such, the research described herewith does not provide much about methodology development, while gives  insightful observations for which combinations would be good choices for use when the models are correctly specified and when they are not. To this end, we will derive their asymptotic distributions and compare their asymptotic variances with the semiparametric efficiency bound in \cite{efficiencybound}.  Particularly, the messages about the estimations with misspecified models are  new and interesting. We consider both locally misspecified and globally misspecified scenarios. Here, the local misspecifition means that the misspecified model is only distinct from the correctly specified model at a rate converging to zero as the sample size $n \to \infty$, and the global misspecification means that the model cannot converge to a correctly specified model. We mainly discuss the local misspecifications for parametric models as they are often popularly used. The details can be found in Section~2. The main findings are listed as follows.

\begin{itemize}
\item When both $PS$ and $OR$ models are correctly specified, all combinations of PS and OR estimator share the same  asymptotic efficiency. This is expected and coincides with the existing studies for some of the combinations in the literature.
\item  When $OR$ model is globally misspecified parametrically and semiparametrically while $PS$ model is not, the consistency of any combination is unaffected, but  the asymptotic variance is in general enlarged except for the cases when nonparametric $PS$ model is applied. In other words, nonparametrically estimating $PS$ model helps improve the estimation efficiency. Under the local misspecification, the asymptotic efficiency can be achieved.
\item In contrast, when $PS$ model is globally misspecified parametrically and semiparametrically while $OR$ model is not, we cannot have a definitive result about whether the asymptotic efficiency is worsen comparing with the semi-parametric efficiency bound though the consistency is still guaranteed. In some cases, there is even a ``super-efficiency phenomenon" which the variance can even be smaller than the bound. We will give an example to show this phenomenon in Section~3. Again, when $OR$ model is estimated nonparametrically and a misspecified $PS$ model is used, the asymptotic efficiency still holds. As previously mentioned, nonparametrically estimating $OR$ model can improve the estimation efficiency. Again under the local misspecification, the asymptotic efficiency can be achieved.
\item From the above, we can see that nonparametric estimation does help on improving asymptotic efficiency. However, this does not mean that it is always recommendable, particularly in high-dimensional scenarios, because it makes the tuning parameter in nonparametric estimation very difficult to choose and clearly causes estimation inefficacy. To reduce the impact of misspecification, semiparametric models, particularly with dimension reduction structure, could be a good choice.
\end{itemize}
All findings are summarised in the following table in which the black cells mean without such  combinations.
\begin{figure}[htbp]
\centering
\includegraphics[width=14cm]{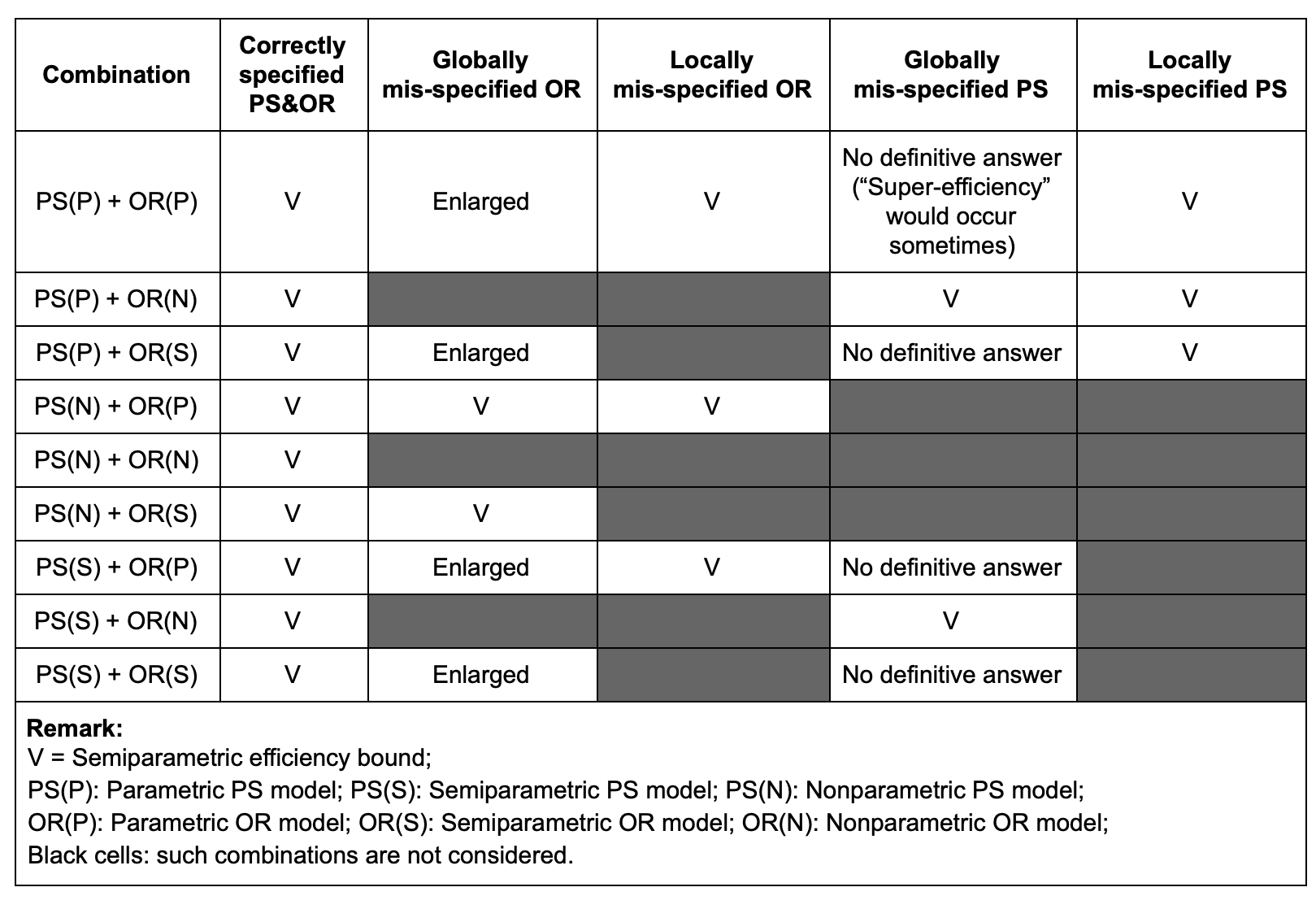}
\end{figure}

The rest of the paper is organized as follows: In Section 2, we introduce the counterfactual framework first to define average treatment effect and formalize notations, and then  discuss doubly robust estimators  and possible estimation methods for $PS$ and $OR$. In addition,  we also  introduce the concept of local mis-specification. In Section 3 we present the asymptotic properties of doubly robust estimators under various scenarios, and the comparisons between our conclusions and existing literature. Section 4 includes simulation studies and Section 5 summarizes the main conclusions in this article. Technical proofs are given in Appendix.

\section{Doubly Robust Estimation}

\subsection{Notation and setup}

Let $D$ be an indicator of observed treatment status ($D=1$ if treated, $D=0$ if untreated) and $X$ be a $p$-dimensional vector of covariates not affected by the treatment status with $p\geq2$. Let $\mathcal{X}$ be the support of $X$. We adopt the counterfactual outcome framework (see \cite{framework1} and \cite{framework2}) here to estimate the average treatment effect. Each individual is assumed to have potential outcomes: $Y(1)$, if the subject has received treatment, and $Y(0)$, if the subject hasn't received treatment. Let $Y$ be the observed outcome given by $(1-D)Y(0)+DY(1)$. In real situation, we can observe either $Y(1)$ or $Y(0)$ but not both of them for each individual in the sample, so it is impossible to observe the average treatment effect directly. The goal is to estimate the average treatment effect defined as
$$\Delta = \mathbb{E}[Y(1)]- \mathbb{E}[Y(0)]=\theta_1 - \theta_0.$$

\noindent We further make the following assumption throughout this paper.
\begin{assumption} \label{uncounfounded}
(Unconfoundedness) We assume that $D$ and $(Y(1), Y(0))$ are conditionally independent given $X$.
\end{assumption}

\noindent As mentioned previously, the prototypical doubly robust estimator proposed by \cite{dr} incorporates the information in both $PS$ model and $OR$ model so that it remains consistent even if one of the $PS$ or $OR$ model is misspecified. There are different choices of $PS$ model and $OR$ model including parametric model, nonparametric model and semiparametric model. In the next section, doubly robust estimators under different model combinations will be presented.

\subsection{Estimation procedures and further assumptions}

Define respectively the true $PS$ and $OR$ model as $P(D=1|X)=p(X)$, and $ \mathbb{E}[Y|X,D=1]= \mathbb{E}[Y(1)|X]=m_1(X)$ and $ \mathbb{E}[Y|X,D=0]= \mathbb{E}[Y(0)|X]=m_0(X)$. Then the average treatment effect can be identified by
$$\Delta=\theta_1-\theta_0= \mathbb{E}\left[\frac{DY}{p(X)}+\left(1-\frac{D}{p(X)}\right)m_1(X)-\frac{(1-D)Y}{1-p(X)}-\left(1-\frac{1-D}{1-p(X)}\right)m_0(X)\right].$$

Let $\big\{x_i,d_i,y_i\big\}^{n}_{i=1}$ be an independent random sample of size $n$ from the joint distribution of $(X,D,Y)$. Note that $x_i$ is a $p$-dimensional vector of covariates, $d_i$ is the binary indicator of treatment status and $y_i$ is the response of \textsl{i}-th individual. From \cite{dr}, the doubly robust estimator is defined as
\begin{equation*}
\begin{split}
\hat{\Delta}&=n^{-1}\sum_{i=1}^{n}\left[\frac{d_iy_i}{\hat{p}(x_i)}+\left(1-\frac{d_i}{\hat{p}(x_i)}\right)\hat{m}_1(x_i)\right]-n^{-1}\sum_{i=1}^{n}\left[\frac{(1-d_i)y_i}{1-\hat{p}(x_i)}+\left(1-\frac{1-d_i}{1-\hat{p}(x_i)}\right)\hat{m}_0(x_i)\right]\\
&=\hat{\theta}_1 - \hat{\theta}_0,
\end{split}
\end{equation*}

\noindent where $\hat{p}(x)$ is an estimated propensity score, $\hat{m}_1(x)$ and $\hat{m}_0(x)$ are estimated outcome regression models, which have different formulas under different model structures. As all combinations discussed in this paper are convergent to some quantities $\Delta^*$, we then write $\hat{\Delta} \to \Delta^*$ in probability as $n\to\infty$. Note that $\Delta=\Delta^*$ when either (but not necessarily both) $PS$ model or $OR$ model is correctly specified due to the double robustness property.

Firstly, when parametric models are considered, without loss of generality, we assume a logistic regression model $\Tilde{p}(x;\beta)=\frac{\exp(x^T\beta )}{1+\exp(x^T\beta)}$ with true parameter $\beta_0$ as the $PS$ model and linear regression models $\Tilde{m}_1(x;\gamma_1)=x^T\gamma_1$ and $\Tilde{m}_0(x;\gamma_0)=x^T\gamma_0$ with true parameter $\gamma_{1,0}$ and $\gamma_{0,0}$ as the $OR$ models. Maximum likelihood estimation (MLE) is used to estimate the unknown parameters. Denote the estimators respectively as $\hat{\beta}$, $\hat{\gamma}_1$ and $\hat{\gamma}_0$. We further make the following assumptions on these proposed models.

\begin{assumption} \label{psassumption}
Let $\Theta_{\beta} \subset \mathbb{R}^p$ be the parameter space for $\beta$ which is open and convex. We assume that the proposed propensity score model $\Tilde{p}(x;\beta):  \mathbb{R}^p \to \mathbb{R}$ is differentiable with respect to $\beta$. Further, we assume that $\Tilde{p}(x;\beta)$ is bounded away from 0 and 1 for any $\beta \in \Theta_{\beta}$.
\end{assumption}

\begin{assumption} \label{orassumption}
Let $\Theta_{\gamma_0}\subset \mathbb{R}^p$ and $\Theta_{\gamma_1}\subset \mathbb{R}^p$ be the parameter space for $\gamma_0$ and $\gamma_1$ respectively which are open and convex. We assume that the proposed outcome regression model $\Tilde{m}_j(x;\gamma_j): \mathbb{R}^p \to \mathbb{R}$ is differentiable with respect to $\gamma_j$, $j=0,1$.
\end{assumption}

\noindent According to \cite{mle1}, when models are correctly specified, we have $\sqrt{n}(\hat{\beta}-\beta_0)\xrightarrow[]{d}N(0,I^{-1}(\beta_0))$, $\sqrt{n}(\hat{\gamma}_1-\gamma_{1,0})\xrightarrow[]{d}N(0,I^{-1}(\gamma_{1,0}))$ and $\sqrt{n}(\hat{\gamma}_0-\gamma_{0,0})\xrightarrow[]{d}N(0,I^{-1}(\gamma_{0,0}))$, where $I(\beta_0)$, $I(\gamma_{1,0})$ and $I(\gamma_{0,0})$ are the Fisher information matrices. When models are misspecified, the convergence of MLE can also be obtained. See \cite{mle2}.  We have $\sqrt{n}(\hat{\beta}-\beta^*)\xrightarrow[]{d}N(0,V(\beta^*))$, where $V(\beta^*)$ is the information sandwich variance matrix. Note that $\beta^*$ is the value of $\beta$ which minimizes the Kullback–Leibler discrepancy with respect to $\beta$. Similarly, we have $\sqrt{n}(\hat{\gamma}_1-\gamma^*_1)\xrightarrow[]{d}N(0,V(\gamma^*_1))$ and $\sqrt{n}(\hat{\gamma}_0-\gamma^*_0)\xrightarrow[]{d}N(0,V(\gamma^*_0))$. Further, we introduce the concept of local misspecification for parametric models. Suppose the correctly specified models have the following forms:
\begin{equation} \label{localdef}
\begin{split}
p(x)&=\Tilde{p}(x;\beta_0)(1+\delta \times s(x)),\\
m_1(x)&=\Tilde{m}_1(x;\gamma_{1,0})+\delta_1\times s_1(x),\\ m_0(x)&=\Tilde{m}_0(x;\gamma_{0,0})+\delta_0 \times s_0(x).
\end{split}
\end{equation}
If $\delta$ is a nonzero fixed constant, we say that $p(x)$ is globally misspecified. If $\delta \to 0$, we say it is locally misspecified. Similarly, we can define the global and local misspecification for $m_1(x)$ and $m_0(x)$.

Secondly, when semiparametric models are considered, we propose the $PS$ model $g(\alpha^TX):=P(D=1|\alpha^TX)$ and the $OR$ models $r_1(\alpha^T_1X):=\mathbb{E}[Y(1)|\alpha^T_1X, D=1]$ and $r_0(\alpha^T_0X):=\mathbb{E}[Y(0)|\alpha^T_0X,\\ D=0]$ with dimension reduction structures $\alpha^TX$, $\alpha^T_1X$ and $\alpha^T_0X$ respectively. Similarly, we can define alternative $PS$ models $q_1(\alpha^T_1X):=P(D=1|\alpha^T_1X)$ and $q_0(\alpha^T_0X):=P(D=1|\alpha^T_0X)$. Note that $p(X)=P(D=1|X)=P(D=1|\alpha^TX)=g(\alpha^TX)$, $m_1(X)=E[Y(1)|X]=E[Y(1)|\alpha^T_1X]=E[Y(1)|\alpha^T_1X,D=1]=r_1(\alpha^T_1X)$ and $m_0(X)=E[Y(0)|X]=E[Y(0)|\alpha^T_0X]=E[Y(0)|\alpha^T_0X,D=0]=r_0(\alpha^T_0X)$ if and only if the dimension reduction structures are correctly specified. We assume that $\alpha$, $\alpha_1$ and $\alpha_0$ are vectors whose Euclidean norms equal 1. There are several available methods of obtaining root-$n$ consistent estimators for $\alpha$, $\alpha_1$ and $\alpha_0$ as mentioned in \cite{semidr}. Therefore, the impact of estimating $\alpha$, $\alpha_1$ and $\alpha_0$ is not considered in this paper.  The corresponding semiparametric estimators are $\hat{g}(\alpha^Tx)$,  $\hat{r}_1(\alpha^T_1x)$ and $\hat{r}_0(\alpha^T_0x)$ with
\begin{equation} \label{semidef}
\begin{split}
\hat{g}(\alpha^Tx)&=\frac{\sum_{j=1}^{n}d_jL_b(\alpha^Tx,\alpha^Tx_j)}{\sum_{j=1}^{n}L_b(\alpha^Tx,\alpha^Tx_j)},\\
\hat{r}_1(\alpha^T_1x)&=\frac{\sum_{j=1}^{n}d_jy_jK_{h_{m_1}}(\alpha^T_1x,\alpha^T_1x_j)}{\sum_{j=1}^{n}d_jK_{h_{m_1}}(\alpha^T_1x,\alpha^T_1x_j)},\\
\hat{r}_0(\alpha^T_0x)&=\frac{\sum_{j=1}^{n}(1-d_j)y_jK_{h_{m_0}}(\alpha^T_0x,\alpha^T_0x_j)}{\sum_{j=1}^{n}(1-d_j)K_{h_{m_0}}(\alpha^T_0x,\alpha^T_0x_j)},
\end{split}
\end{equation}
where $L_b(u,v)=\frac{1}{b}L\left(\frac{u-v}{b}\right)$, $K_{h_{m_1}}(u,v)=\frac{1}{h_{m_1}}K\left(\frac{u-v}{h_{m_1}}\right)$ and $K_{h_{m_0}}(u,v)=\frac{1}{h_{m_0}}K\left(\frac{u-v}{h_{m_0}}\right)$. Note that $K(\cdot): \mathbb{R} \to \mathbb{R}$, $L(\cdot): \mathbb{R} \to \mathbb{R}$ are kernel functions of order 2 and $b,h_{m_1},h_{m_0}$ are corresponding bandwidths. We further make the following assumption for the kernel functions and bandwidths.

\begin{assumption} \label{semiassumption}
Kernel functions $K(\cdot)$ and $L(\cdot)$ are symmetric around 0, compactly supported and at least twice continuously differentiable with $\int u^2K(u)du<\infty$ and $\int u^2L(u)du<\infty$.
The bandwidths $b,h_{m_1},h_{m_0}$ satisfy the following conditions as $n \to \infty$: (a) $b\to0$, $nb\to \infty$, $nb^3\to \infty$, $nb^4\to0$, $\log(n)/(nb^3)\to 0$; (b) $h_{m_1}, h_{m_0}\to0$, $nh_{m_1}, nh_{m_0}\to \infty$, $nh^3_{m_1}, nh^3_{m_0}\to \infty$, $nh^4_{m_1}, nh^4_{m_0}\to0$ and $\log(n)/(nh^3_{m_1})$, $\log(n)/(nh^3_{m_0}) \to 0$.
\end{assumption}

 Thirdly, when nonparametric models are considered, we assume the estimated $PS$ model $\hat{p}(x)$ and $OR$ models $\hat{m}_1(x)$, $\hat{m}_0(x)$ have the following form:
\begin{equation} \label{nondef}
\begin{split}
 \hat{p}(x)&=\frac{\sum_{j=1}^{n}d_j\Tilde{L}_{\Tilde{b}}(x,x_j)}{\sum_{j=1}^{n}\Tilde{L}_{\Tilde{b}}(x,x_j)}, \\
\hat{m}_1(x)&=\frac{\sum_{j=1}^{n}d_jy_j\Tilde{K}_{\Tilde{h}_{m_1}}(x,x_j)}{\sum_{j=1}^{n}d_j\Tilde{K}_{\Tilde{h}_{m_1}}(x,x_j)}, \quad
\hat{m}_0(x)=\frac{\sum_{j=1}^{n}(1-d_j)y_j\Tilde{K}_{\Tilde{h}_{m_0}}(x,x_j)}{\sum_{j=1}^{n}(1-d_j)\Tilde{K}_{\Tilde{h}_{m_0}}(x,x_j)},
\end{split}
\end{equation}
where $\Tilde{L}_{\Tilde{b}}(u,v)=\frac{1}{\Tilde{b}^p}\Tilde{L}\left(\frac{u-v}{\Tilde{b}}\right)$, $\Tilde{K}_{\Tilde{h}_{m_1}}(u,v)=\frac{1}{\Tilde{h}^p_{m_1}}\Tilde{K}\left(\frac{u-v}{\Tilde{h}_{m_1}}\right)$ and $\Tilde{K}_{\Tilde{h}_{m_0}}(u,v)=\frac{1}{\Tilde{h}^p_{m_0}}\Tilde{K}\left(\frac{u-v}{\Tilde{h}_{m_0}}\right)$. Note that $\Tilde{K}(\cdot):  \mathbb{R}^p \to \mathbb{R}$ and $\Tilde{L}(\cdot):  \mathbb{R}^p \to \mathbb{R}$ are kernel functions of order s, where $s>p$ is a positive integer. The corresponding bandwidths are $\Tilde{b},\Tilde{h}_{m_1},\Tilde{h}_{m_0}$. We further make the following assumption for the kernel functions and bandwidths.

\begin{assumption} \label{nonparaassumption}
Kernel functions $\Tilde{K}(\cdot)$ and $\Tilde{L}(\cdot)$ are symmetric around 0, compactly supported and at least s times continuously differentiable with $\int u^s\Tilde{K}(u)du<\infty$ and $\int u^s\Tilde{L}(u)du<\infty$. The bandwidths $\Tilde{b},\Tilde{h}_{m_1},\Tilde{h}_{m_0}$ satisfy the following conditions as $n \to \infty$: (a) $\Tilde{b}\to0$, $n\Tilde{b}^{p+2}\to \infty$, $n\Tilde{b}^{2s}\to0$, $\log(n)/(n\Tilde{b}^{p+s})\to 0$; (b) $\Tilde{h}_{m_1}, \Tilde{h}_{m_0}\to0$, $n\Tilde{h}^{p+2}_{m_1}, n\Tilde{h}^{p+2}_{m_0}\to \infty$, $n\Tilde{h}^{2s}_{m_1}, n\Tilde{h}^{2s}_{m_0}\to0$ and $\log(n)/(n\Tilde{h}^{p+s}_{m_1}), \log(n)/(n\Tilde{h}^{p+s}_{m_0}) \to 0$.
\end{assumption}

Furthermore, let $f(x): \mathbb{R}^p \to \mathbb{R}$ be the density function of $X$, $\Tilde{f}(\alpha^Tx): \mathbb{R} \to \mathbb{R}$ be the density function of $\alpha^TX$, $\Tilde{f}_1(\alpha^T_1x): \mathbb{R} \to \mathbb{R}$ be the density function of $\alpha^T_1X$ and $\Tilde{f}_0(\alpha^T_0x): \mathbb{R} \to \mathbb{R}$ be the density function of $\alpha^T_0X$. Recall that the true $PS$ model $p(x): \mathbb{R}^p \to \mathbb{R}$, the proposed semiparametric $PS$ model $g(\alpha^Tx): \mathbb{R} \to \mathbb{R}$ and the alternative $PS$ models $q_1(\alpha^T_1x): \mathbb{R} \to \mathbb{R}$ and $q_0(\alpha^T_0X): \mathbb{R} \to \mathbb{R}$ are defined as $p(X):=P(D=1|X)$, $g(\alpha^TX):=P(D=1|\alpha^TX)$, $q_1(\alpha^T_1X):=P(D=1|\alpha^T_1X)$ and $q_0(\alpha^T_0X):=P(D=1|\alpha^T_0X)$. These functions are useful in deriving the asymptotic distribution of $\hat{\Delta}$. We make the following assumption about $f(\cdot), \Tilde{f}(\cdot), \Tilde{f}_1(\cdot), \Tilde{f}_0(\cdot), p(\cdot), g(\cdot), q_1(\cdot),q_0(\cdot)$ throughout the paper.
\begin{assumption} \label{general}
Density functions $f(\cdot), \Tilde{f}(\cdot), \Tilde{f}_1(\cdot), \Tilde{f}_0(\cdot)$ and propensity score models $p(\cdot), g(\cdot), q_1(\cdot),q_0(\cdot)$ are bounded away from 0 and 1.
\end{assumption}

As a result, we can obtain the following nine estimators using different combinations of $PS$ and $OR$ estimator:
\begin{equation}\label{drdef}
\begin{split}
&\hat{\Delta}_1=\frac1n\sum_{i=1}^{n}\left\{\frac{d_iy_i}{\Tilde{p}(x_i;\hat{\beta})}+\left(1-\frac{d_i}{\Tilde{p}(x_i;\hat{\beta})}\right)\Tilde{m}_1(x_i;\hat{\gamma}_1)-\frac{(1-d_i)y_i}{1-\Tilde{p}(x_i;\hat{\beta})}+\left(1-\frac{1-d_i}{1-\Tilde{p}(x_i;\hat{\beta})}\right)\Tilde{m}_0(x_i;\hat{\gamma}_0)\right\}\\
&\hat{\Delta}_2=\frac1n\sum_{i=1}^{n}\left\{\frac{d_iy_i}{\Tilde{p}(x_i;\hat{\beta})}+\left(1-\frac{d_i}{\Tilde{p}(x_i;\hat{\beta})}\right)\hat{m}_1(x_i)-\frac{(1-d_i)y_i}{1-\Tilde{p}(x_i;\hat{\beta})}+\left(1-\frac{1-d_i}{1-\Tilde{p}(x_i;\hat{\beta})}\right)\hat{m}_0(x_i)\right\}\\
&\hat{\Delta}_3=\frac1n\sum_{i=1}^{n}\left\{\frac{d_iy_i}{\hat{p}(x_i)}+\left(1-\frac{d_i}{\hat{p}(x_i)}\right)\Tilde{m}_1(x_i;\hat{\gamma}_1)-\frac{(1-d_i)y_i}{1-\hat{p}(x_i)}+\left(1-\frac{1-d_i}{1-\hat{p}(x_i)}\right)\Tilde{m}_0(x_i;\hat{\gamma}_0)\right\}\\
&\hat{\Delta}_4=\frac1n\sum_{i=1}^{n}\left\{\frac{d_iy_i}{\hat{p}(x_i)}+\left(1-\frac{d_i}{\hat{p}(x_i)}\right)\hat{m}_1(x_i)-\frac{(1-d_i)y_i}{1-\hat{p}(x_i)}+\left(1-\frac{1-d_i}{1-\hat{p}(x_i)}\right)\hat{m}_0(x_i)\right\}\\
&\hat{\Delta}_5=\frac1n\sum_{i=1}^{n}\left\{\frac{d_iy_i}{\hat{g}(\alpha^Tx_i)}+\left(1-\frac{d_i}{\hat{g}(\alpha^Tx_i)}\right)\Tilde{m}_1(x_i;\hat{\gamma}_1)-\frac{(1-d_i)y_i}{1-\hat{g}(\alpha^Tx_i)}+\left(1-\frac{1-d_i}{1-\hat{g}(\alpha^Tx_i)}\right)\Tilde{m}_0(x_i;\hat{\gamma}_0)\right\}\\
&\hat{\Delta}_6=\frac1n\sum_{i=1}^{n}\left\{\frac{d_iy_i}{\Tilde{p}(x_i;\hat{\beta})}+\left(1-\frac{d_i}{\Tilde{p}(x_i;\hat{\beta})}\right)\hat{r}_1(\alpha^T_1x_i)-\frac{(1-d_i)y_i}{1-\Tilde{p}(x_i;\hat{\beta})}+\left(1-\frac{1-d_i}{1-\Tilde{p}(x_i;\hat{\beta})}\right)\hat{r}_0(\alpha^T_0x_i)\right\}\\
&\hat{\Delta}_7=\frac1n\sum_{i=1}^{n}\left\{\frac{d_iy_i}{\hat{g}(\alpha^Tx_i)}+\left(1-\frac{d_i}{\hat{g}(\alpha^Tx_i)}\right)\hat{m}_1(x_i)-\frac{(1-d_i)y_i}{1-\hat{g}(\alpha^Tx_i)}+\left(1-\frac{1-d_i}{1-\hat{g}(\alpha^Tx_i)}\right)\hat{m}_0(x_i)\right\}\\
&\hat{\Delta}_8=\frac1n\sum_{i=1}^{n}\left\{\frac{d_iy_i}{\hat{p}(x_i)}+\left(1-\frac{d_i}{\hat{p}(x_i)}\right)\hat{r}_1(\alpha^T_1x_i)-\frac{(1-d_i)y_i}{1-\hat{p}(x_i)}+\left(1-\frac{1-d_i}{1-\hat{p}(x_i)}\right)\hat{r}_0(\alpha^T_0x_i)\right\}\\
&\hat{\Delta}_9=\frac1n\sum_{i=1}^{n}\left\{\frac{d_iy_i}{\hat{g}(\alpha^Tx_i)}+\left(1-\frac{d_i}{\hat{g}(\alpha^Tx_i)}\right)\hat{r}_1(\alpha^T_1x_i)-\frac{(1-d_i)y_i}{1-\hat{g}(\alpha^Tx_i)}+\left(1-\frac{1-d_i}{1-\hat{g}(\alpha^Tx_i)}\right)\hat{r}_0(\alpha^T_0x_i)\right\}.
\end{split}
\end{equation}

We can show the consistencies of these estimators even if one of $PS$ or $OR$ model is misspecified, see Appendix 6.1. In the next section, we focus on studying their asymptotic distributions.

\section{Asymptotic distributions}
In this section, we derive the asymptotic distributions of the proposed estimators. The comparisons between their asymptotic variances and the semiparametric efficiency bound are also presented. Detailed proofs can be found in Appendix 6.2.

\begin{theorem}
	Suppose that the $PS$ and $OR$ models are correctly specified. Under Assumptions 1-6 in Section~2, for all nine combinations, we have $\sqrt{n}(\hat{\Delta}_k-\Delta)\xrightarrow{d} N(0,\Sigma_1)$ for $k=1,\cdots,9$ with
	$$\Sigma_1=\mathbb{E}\left\{\frac{\mathrm{Var}[Y(1)|X]}{p(X)}+\frac{\mathrm{Var}[Y(0)|X]}{1-p(X)}+\left[m_1(X)-\mathbb{E}[Y(1)]-m_0(X)+\mathbb{E}[Y(0)]\right]^2\right\},$$
	which is the same as the semiparametric efficiency bound shown by \cite{efficiencybound}.
\end{theorem}

\begin{remark}
	The results for $\hat{\Delta}_1$ (parametric+parametric) and $\hat{\Delta}_4$ (nonparametric+nonparametric) coincide with the results in the literature, see e.g. \cite{drpara+pspara} and \cite{denonpara}. For $\hat{\Delta}_9$ (semiparametric+semiparametric), the result is similar to that in \cite{semidr} in the context of missing data. The other results are newly derived in this paper.
\end{remark}

\begin{theorem}
	Suppose that the $PS$ model is correctly specified and the $OR$ model is globally misspecified with fixed nonzero $\delta_1$ and $\delta_0$.  We then have the estimators $\hat{\Delta}_k$ for $k=1, 3,5, 6, 8, 9$. Under Assumptions 1-6 in Section~2,
	\begin{equation}
	\begin{split}
	&\sqrt{n}(\hat{\Delta}_k-\Delta)\xrightarrow{d} N(0,\Sigma_1), \quad \mbox{for \, \, $k=3,8$}\\
	&\sqrt{n}(\hat{\Delta}_k-\Delta)\xrightarrow{d} N(0,\Sigma_2), \quad \mbox{for \, \, $k=1,5$}\\
	&\sqrt{n}(\hat{\Delta}_k-\Delta)\xrightarrow{d} N(0,\Sigma_3), \quad \mbox{for \, \, $k=6,9$}\\
	\end{split}
	\end{equation}
	where $\hat{\Delta}_3$ is nonparametric+misspecified parametric and $\hat{\Delta}_8$ is nonparametric+misspecified semiparametric, $\Sigma_1$ is defined in Theorem~1 and  $\Sigma_2$ and $\Sigma_3$  are as follows:
	\begin{equation*}
	\begin{split}
	&\Sigma_2 =\Sigma_1+\mathbb{E}\bigg\{\Big[\sqrt{\frac{1}{p(X)}-1}[\Tilde{m}_1(X;\gamma^*_1)-m_1(X)]+\sqrt{\frac{1}{1-p(X)}-1}[\Tilde{m}_0(X;\gamma^*_0)-m_0(X)]\\
	&\qquad+\sqrt{p(X)(1-p(X))}w(X)\Big]^2\bigg\} \geq \Sigma_1, \\
	&\Sigma_3=\Sigma_1+\mathbb{E}\bigg\{\Big[\sqrt{\frac{1}{p(X)}-1}[r_1(\alpha^T_1X)-m_1(X)]+\sqrt{\frac{1}{1-p(X)}-1}[r_0(\alpha^T_0X)-m_0(X)]\\
	&\qquad+\sqrt{p(X)(1-p(X))}w(X)\Big]^2\bigg\} \geq \Sigma_1.\\
	\end{split}
	\end{equation*}
	The equalities hold if $OR$ models are correctly specified. Note that $w(X)$ is different for different estimators. See Appendix 6.2 for details.
\end{theorem}

\begin{remark}
The result for $\hat{\Delta}_1$ (parametric+parametric) coincides with the results in \cite{misspecification-reference1} and \cite{misspecification-reference2}. The other results are newly derived in this paper. The results show some interesting phenomena. Firstly, due to the nonparametric estimation for the correctly specified $PS$ model, $\hat \Delta_k$ for $k=3,8$, we can achieve the asymptotic efficiency. Secondly, under locally misspecification of $OR$ models with $\delta_1\to 0$ and $\delta_0\to 0$, $\Sigma_2$ converges to $\Sigma_1$. That is, the asymptotic variances of $\hat{\Delta}_1$ and $\hat{\Delta}_5$ converge to $\Sigma_1$ as $\delta_1\to 0$ and $\delta_0\to 0$.
\end{remark}

\begin{theorem}
	Suppose that the $PS$ model is globally misspecified with fixed nonzero $\delta$, while the $OR$ model is correctly specified. Under Assumptions 1-6 in Section~2,
	\begin{equation}
	\begin{split}
	&\sqrt{n}(\hat{\Delta}_k-\Delta)\xrightarrow{d} N(0,\Sigma_1), \quad \mbox{for \, \, $k=2,7$},\\
	&\sqrt{n}(\hat{\Delta}_k-\Delta)\xrightarrow{d} N(0,\Sigma_4), \quad \mbox{for \, \, $k=1,6$},\\
	&\sqrt{n}(\hat{\Delta}_k-\Delta)\xrightarrow{d} N(0,\Sigma_5), \quad \mbox{for \, \, $k=5,9$},
	\end{split}
	\end{equation}
	
	\noindent where $\Sigma_4$ and $\Sigma_5$ are as follows:
	\begin{equation*}
	\begin{split}
	&\Sigma_4=\Sigma_1+\mathbb{E}\left\{\frac{1}{p(X)}\textrm{Var}[Y(1)|X]\left[\left(\frac{p(X)}{\Tilde{p}(X;\beta^*)}+w_1(X)p(X)\right)^2-1\right]\right\}\\
	&\qquad+\mathbb{E}\left\{\frac{1}{1-p(X)}\textrm{Var}[Y(0)|X]\left[\left(\frac{1-p(X)}{1-\Tilde{p}(X;\beta^*)}+w_0(X)(1-p(X))\right)^2-1\right]\right\},
	\end{split}
	\end{equation*}
	\begin{equation*}
	\begin{split}
	&\Sigma_5=\Sigma_1+\mathbb{E}\left\{\frac{1}{p(X)}\textrm{Var}[Y(1)|X]\left[\left(\frac{p(X)}{g(\alpha^TX)}+w_1(X)p(X)\right)^2-1\right]\right\}\\
	&\qquad+\mathbb{E}\left\{\frac{1}{1-p(X)}\textrm{Var}[Y(0)|X]\left[\left(\frac{1-p(X)}{1-g(\alpha^TX)}+w_0(X)(1-p(X))\right)^2-1\right]\right\}.
	\end{split}
	\end{equation*}
	Note that $w_1(X)$ and $w_0(X)$ are different for different estimators. See Appendix 6.2 for details.
\end{theorem}

\begin{remark}
	The results show some interesting phenomena. Firstly, again, due to the nonparametric estimation for the $OR$ model, the estimators $\hat \Delta_k$ for $k=2,7$ can achieve the asymptotic efficiency. Secondly, under locally misspecification of $PS$ model with $\delta\to 0$, $\Sigma_4$ converges to $\Sigma_1$. That is, the estimators $\hat{\Delta}_1$ and $\hat{\Delta}_5$ can also achieve the semiparametric efficiency bound. Thirdly, it is difficult to compare $\Sigma_4$ and $\Sigma_5$ with $\Sigma_1$. We still have difficulty to reach a general conclusion. This case is very different from the case with correctly specified $PS$ model and misspecified $OR$ model as stated in Theorem~2.
\end{remark}

Although a general comparison is very difficult to theoretically determine under what circumstances $\Sigma_4$ and $\Sigma_5$ would be smaller than $\Sigma_1$ and under what circumstances they are larger, we give a simple example to show that the asymptotic variance $\Sigma_4$ derived from $\hat \Delta_1$ could be smaller than $\Sigma_1$ in certain cases.
Suppose that the true propensity score function is simply a constant function $p(x)=p^*$ and the assumed propensity score model is also  a constant $g$, where $0< p^*,g < 1$. We further assume that $\mathbb{E}(X)=0, \mathrm{Var}[Y(1)|X]=\mathrm{Var}[Y(0)|X]$.
Then we have $w_1(x)=w_0(x)=0$, the formula of $\Sigma_4$ can be reduced to
$$\Sigma_4 - \Sigma_1=\mathbb{E}\left\{\mathrm{Var}[Y(1)|X]\right\}\left\{\frac{p^*}{g^2}+\frac{1-p^*}{(1-g)^2}-\frac{1}{p^*}-\frac{1}{1-p^*}\right\}.$$

\noindent For each fixed value $p^*$, we can determine whether the asymptotic variance is enlarged or not by looking at the function $f(g)=\frac{p^*}{g^2}+\frac{1-p^*}{(1-g)^2}-\frac{1}{p^*}-\frac{1}{1-p^*}$. That is, if $f(g)=0$, $\Sigma_4=\Sigma_1$; if $f(g)>0$, $\Sigma_4>\Sigma_1$; if $f(g)<0$, $\Sigma_4<\Sigma_1$. In Figure 1, we plot three curves of the $f(g)$ about $g$ with  $p^*=1/4, 1/2, 3/4$. First, we can see that $f(g)=0$ when $p^*=g=1/4, 1/2, 3/4$. This means that when the model is correctly specified, the variance can achieve the semiparametric efficiency bound. Second, when $p^*=1/2$, $f(g)\geq0$ for all values of $g$. In other words, misspecification always causes the variance enlargement.  In contrast, when $p^* \not =1/2$, the situation becomes different. From the curves with $p^*=1/4, 3/4$, we can see that the semiparametric efficiency bound can only be achieved at $g=1/4, 3/4$ accordingly, otherwise, there are ranges of $g$ such that the variances can even be smaller than the bound. This shows possible ``super-efficiency phenomenon" when the misspecification occurs.

\begin{figure}[htbp]
\centering
\includegraphics[width=13cm]{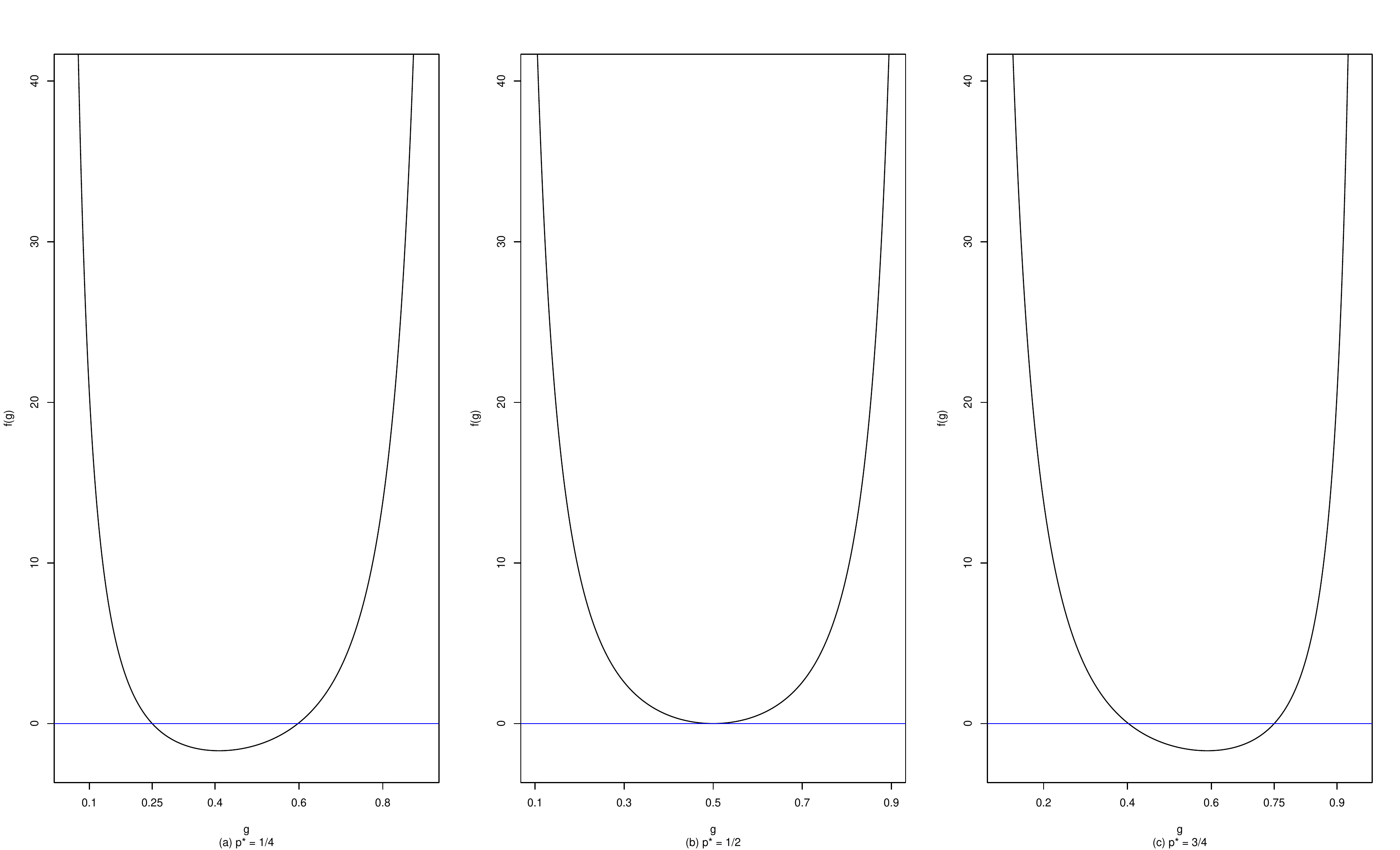}
\caption{The plots of $f(g)$ for each fixed $p^*$. The left and right panels illustrate that the asymptotic variance $\Sigma_4$ derived from $\hat \Delta_1$ could be smaller than the semiparametric efficiency bound $\Sigma_1$.}
\end{figure}

\section{Numerical investigation}
We conduct some Monte Carlo simulations to investigate the performances of these doubly robust estimators  in finite sample scenarios in terms of bias, standard deviation and mean squared error. The experiments are repeated 1000 times, and the sample size is taken to be 1000.

Suppose for each subject $i=1,2,...,n$, the 10-dimensional covariates $X_i=(x_{i1}, ..., x_{i10})^T$ is independently drawn from $N(0, I)$, where $I$ is the $10\times10$ identity matrix. The potential outcomes $Y(1)$ and $Y(0)$ follow $N(\mathbb{E}[Y(1)], 1)$ and $N(\mathbb{E}[Y(0)], 1)$ respectively, where
$$\mathbb{E}[Y(1)]=10+\beta^TX\;\; and\;\; \mathbb{E}[Y(0)]=5+\beta^TX.$$
Let $\beta$ be $(0.5,0.5,0.5,0.5,0,...,0)^T$. Following these two specific regression models, the true ATE is equal to 5. The true propensity scores are determined by a logistic regression model
$$P(D=1|X=x)=\frac{\exp(\alpha^Tx+s_0)}{1+\exp(\alpha^Tx+s_0)},$$
where $\alpha=\alpha^\prime/(1+s^2_1)^{1/2}$, $\alpha^\prime=\beta/\left\lVert \beta\right\rVert+(0,...,0,s_1)^T$. Similar to the setting in missing data, see \cite{simulation1}, let constant $s_0$ control the proportion of treated subjects and let constant $s_1$ control the closeness between $\alpha$ and $\beta$. When $s_1=0$, $\alpha$ and $\beta$ are the same. When $s_1=1$, the angle between $\alpha$ and $\beta$ is $45^\circ$. When $s_1$ is large enough, $\alpha$ and $\beta$ are vertical to each other. For each subject, the treatment indicator $d_i$ is generated from a Bernoulli distribution with parameter $P(D=1|X=x_i)$.

If the parametric method is used to model the propensity score, a logistic regression of $d_i$ on $x_{ij}\,'s$ is regarded as  a correctly specified $PS$ model. Similarly, a linear regression of $y_i$ on $x_{ij}\,'s$ is regarded as the correctly specified $OR$ model. Following the setting in \cite{simulation2}, we introduce covariates $Z_i=(z_{i1},...,z_{i10})^T$ as below:
\begin{equation*}
\begin{split}
&z_{i1}=\exp(x_{i1}/3),\:z_{i2}=\frac{x_{i2}}{1+\exp(x_{i1})}+10,\: z_{i3}=\left(\frac{x_{i1}x_{i3}}{25}+0.6\right)^3,\: z_{i4}=(x_{i2}+x_{i4}+20)^2\\
&z_{i5}=\exp(x_{i5}/3),\:z_{i6}=\frac{x_{i6}}{1+\exp(x_{i5})}+10,\: z_{i7}=\left(\frac{x_{i5}x_{i7}}{25}+0.6\right)^3,\: z_{i8}=(x_{i6}+x_{i8}+20)^2\\
&z_{i9}=\exp(x_{i9}/3),\: z_{i10}=\frac{x_{i10}}{1+\exp(x_{i9})}+10
\end{split}
\end{equation*}
Suppose $Z_i\;'s$ are used instead of $X_i\;'s$, a logistic regression of $d_i$ on $z_{ij}\,'s$ is a misspecified $PS$ model and a linear regression of $y_i$ on $z_{ij}\,'s$ means a misspecified $OR$ model. If the semiparametric method is applied, $\alpha^TX$ is a  correct dimension reduction structure for propensity while $\beta^TX$ is  a mis-specified dimension reduction structure. Similarly, for outcome regression models, $\beta^TX$ is a correct dimension reduction structure while $\alpha^TX$ leads to a mis-specified dimension reduction structure. The kernel functions $K(\cdot)$ and $L(\cdot)$ are taken to be Guassian kernels $K(t)=L(t)=(2\pi)^{(-1/2)}exp(-t^2/2)$. For the nonparametric method, the multiple Guassian kernel $K(t)=L(t)=(2\pi)^{(-p/2)}exp(-||t||^2/2)$ is adopted. Inspired by \cite{semidr}, we consider the effect of closeness between $\alpha$ and $\beta$ and the proportion of untreated subjects on the performances of these $DR$ estimators. We also investigate the impact of mis-specification.

\subsection{Effect of proportion of untreated}
In this section, we first investigate the impact of proportion of untreated subjects on the bias, standard deviation (std) and mean squared error (mse). The closeness of $\alpha$ and $\beta$ is fixed and set to be $45^\circ$. We consider three scenarios in which the proportion of untreated is chosen to be 25\%, 50\% and 75\%. The simulation results are summarised in Tables~1,  2 and  3. Compared to the case where the proportions of untreated are chosen to be around 25\% and 75\%, the stds and mses are the smallest when untreated subjects are about 50\%. In terms of bias, a balanced design with 50\% of untreated subjects gives the smallest biases  for  most of the estimators. For $\hat{\Delta}_1,\hat{\Delta}_3,\hat{\Delta}_5$ in Table~1 and $\hat{\Delta}_1$, $\hat{\Delta}_5$ in Table~3, the biases of these estimators in three scenarios are small enough so that we can ignore the impact from the proportion of untreated.

Secondly, we compare the performances of estimators when models are correctly specified, see Table~1. As we proved in theory, all estimators have the identical asymptotic variance. In the finite sample cases, we can observe that they perform similarly with regard to stds, and $\hat{\Delta}_8$ with nonparametric $PS$ and semiparametric $OR$ estimation works well in any scenario. In terms of biases, the biases of $\hat{\Delta}_2$, $\hat{\Delta}_4$ and $\hat{\Delta}_7$ are greater than that of other estimators in any scenario. This seems to mean that  the bias becomes larger if we use nonparametric  estimation for the $OR$ model. Among these three estimators, $\hat{\Delta}_4$ with both nonparametric $PS$ and $OR$ estimation has the largest bias, which shows the inefficiency in nonparametric estimation.

Finally, we explore the influence of misspecification, see Tables~2 and 3 and compare the results with Table~1. It is noteworthy that model misspecification has the least impact on std when there is 50\% of untreated subjects. Theoretically, $\hat{\Delta}_k$, $k=1,5,6,9$ are consistent but not efficient when $OR$ model is misspecified. In Table 2, we observe that the biases and stds of $\hat{\Delta}_1$ and $\hat{\Delta}_5$ are considerably enlarged. For $\hat{\Delta}_6$ and $\hat{\Delta}_9$, we only see slightly increases in the biases and stds. Recall that in theory  $\hat{\Delta}_3$ and $\hat{\Delta}_8$ can achieve the semiparametric efficiency bound when $OR$ is misspecified. However, in the finite sample case we conduct, we do not see this property for $\hat{\Delta}_3$ as its bias and std are significantly enlarged in Table 2. Theoretically, $\hat{\Delta}_k$, $k=1,5,6,9$ are consistent, but their efficiencies cannot be determined when $PS$ model is mis-specified. In Table~3, we observe that the stds of these estimators in the simulation are not necessarily enlarged, which coincides with our theoretical results. Recall that  $\hat{\Delta}_2$ and $\hat{\Delta}_7$ can achieve the semiparametric efficiency bound when $PS$ is mis-specified. Their performances in Table 3 are similar to their performances in Table 1, which supports our theory. Therefore, overall, the numerical results support the theoretical conclusions.

\subsection{Effect of closeness between $\alpha$ and $\beta$}
We now examine the influence of closeness between $\alpha$ and $\beta$. The proportion of treated subjects is  taken to be around $50\%$. We consider three scenarios in which the closeness between $\alpha$ and $\beta$ is chosen to be $0^\circ$, $45^\circ$ and $90^\circ$. The simulation results are presented in Tables~4, 5 and 6. We can observe that the larger the angle is, the smaller the bias is for most of the estimators, but there is some influence on $std$ and $mse$. For $\hat{\Delta}_1,\hat{\Delta}_3,\hat{\Delta}_5$ in Table 4 and $\hat{\Delta}_1$, $\hat{\Delta}_5$ in Table 6, the biases of these estimators in three scenarios are small so that we can ignore the impact from the  proportion of untreated.

Compare the performances of estimators when models are correctly specified (Table 4).  We can observe that these estimators perform similarly with regard to stds, but the most efficient estimator is again $\hat{\Delta}_8$ in any scenario. In terms of biases, the biases of $\hat{\Delta}_2$, $\hat{\Delta}_4$ and $\hat{\Delta}_7$ are significantly greater than those of other estimators when the angle is set to be $0^\circ$. However, as the angle increases, the gap becomes smaller.

Finally, we explore the influence of misspecification, see Tables~5 and  6 and compare the results with  those in Table 4. Note that the misspecification of semiparametric models no longer exists when the angle is set to be $0^\circ$.  From Table 5, the stds of $\hat{\Delta}_1$ and $\hat{\Delta}_5$ are considerably enlarged. The enlargements become smaller as the angle increases. The stds of $\hat{\Delta}_6$, $\hat{\Delta}_8$ and $\hat{\Delta}_9$ increase slightly. The enlargements become more seriously in scenario 3. Again, theoretically consistent and efficient estimator $\hat{\Delta}_3$ does not perform well in the limited numerical study. With increase of the angle, the enlargements on bias and std of $\hat{\Delta}_3$ reasonably reduce. Theoretically, $\hat{\Delta}_k$, $k=1,5,6,9$ are consistent, but we are unable to make a definitive comparison in terms of their asymptotic variances when the $PS$ model is misspecified. In accordance with the theory in this paper, the stds of these estimators in the simulation are not necessarily enlarged. Recall that the variance of $\hat{\Delta}_2$ and $\hat{\Delta}_7$ still achieves the semiparametric efficiency bound when $PS$ is misspecified. The comparison between the results in Table 6 and Table 4, we can see the coincidence with the theory.

In summary, the proportion of untreated has an impact on biases, stds and mses while the closeness between $\alpha$ and $\beta$ has an impact on biases only. The misspecification of $PS$ model seems to have less impact for bias, std and mse than the misspecification of $OR$ model. This effect is much more serious when the parametric estimation of $OR$ model is used.

\begin{table}[!htbp]
\centering
\caption{Correctly specified $PS$ and $OR$ models, $45^{\circ}$ between $\alpha$ and $\beta$}
\begin{tabular}{l*{10}{c}}
\toprule
Estimator &  \multicolumn{3}{c}{25\% of untreated subjects} & \multicolumn{3}{c}{50\% of untreated subjects} & \multicolumn{3}{c}{75\%  of untreated subjects} \\
\midrule
{}   & bias & std & mse & bias & std & mse & bias & std & mse \\
$\hat{\Delta}_1$  & -0.0024 & 0.0970 & 0.0094 & -0.0041 & 0.0759 & 0.0058 & -0.0027 & 0.0974 & 0.0095\\
$\hat{\Delta}_2$  & 0.1868 & 0.0925 & 0.0435 & 0.1448 & 0.0749 & 0.0266 & 0.1852 & 0.0940 & 0.0431\\
$\hat{\Delta}_3$  & -0.0024 & 0.0897 & 0.0081 & -0.0030 & 0.0722 & 0.0052 & -0.0013 & 0.0898 & 0.0081\\
$\hat{\Delta}_4$  & 0.3146 & 0.0930 & 0.1076 & 0.2937 & 0.0768 & 0.0922 & 0.3130 & 0.0923 & 0.1065\\
$\hat{\Delta}_5$  & -0.0023 & 0.0922 & 0.0085 & -0.0033 & 0.0735 & 0.0054 & -0.0022 & 0.0923 & 0.0085\\
$\hat{\Delta}_6$  & 0.0052 & 0.0954 & 0.0091 & 0.0005 & 0.0761 & 0.0058 & 0.0045 & 0.0956 & 0.0091\\
$\hat{\Delta}_7$  & 0.2144 & 0.0962 & 0.0552 & 0.1762 & 0.0797 & 0.0374 & 0.2126 & 0.0967 & 0.0545\\
$\hat{\Delta}_8$  & 0.0426 & 0.0760 & 0.0076 & 0.0322 & 0.0673 & 0.0056 & 0.0428 & 0.0763 & 0.0076\\
$\hat{\Delta}_9$  & 0.0140 & 0.0871 & 0.0078 & 0.0084 & 0.0724 & 0.0053 & 0.0138 & 0.0870 & 0.0078\\
\bottomrule
\end{tabular}
\end{table}

\begin{table}[!htbp]
\centering
\caption{Misspecified $OR$ model, $45^{\circ}$ between $\alpha$ and $\beta$}
\begin{tabular}{l*{10}{c}}
\toprule
Estimator &  \multicolumn{3}{c}{25\% of untreated subjects} & \multicolumn{3}{c}{50\% of untreated subjects} & \multicolumn{3}{c}{75\%  of untreated subjects} \\
\midrule
{}   & bias & std & mse & bias & std & mse & bias & std & mse \\
$\hat{\Delta}_1$  & 0.0240 & 0.4105 & 0.1689 & -0.0099 & 0.2204 & 0.0486 & 0.0129 & 0.4012 & 0.1610\\
$\hat{\Delta}_3$  & 1.4122 & 0.5267 & 2.2714 & 0.3671 & 0.3397 & 0.2500 & -0.9007 & 0.7190 & 1.3276\\
$\hat{\Delta}_5$  & 0.3649 & 0.3400 & 0.2486 & 0.0841 & 0.1790 & 0.0391 & -0.2117 & 0.3252 & 0.1505\\
$\hat{\Delta}_6$  & 0.0061 & 0.1019 & 0.0104 & 0.0026 & 0.0776 & 0.0060 & 0.0058 & 0.1041 & 0.0109\\
$\hat{\Delta}_8$  & 0.0621 & 0.0893 & 0.0118 & 0.0490 & 0.0788 & 0.0086 & 0.0633 & 0.0901 & 0.0121\\
$\hat{\Delta}_9$  & 0.0167 & 0.1076 & 0.0119 & 0.0132 & 0.0913 & 0.0085 & 0.0168 & 0.1101 & 0.0124\\
\bottomrule
\end{tabular}
\end{table}

\begin{table}[!htbp]
\centering
\caption{Misspecified $PS$ model, $45^{\circ}$ between $\alpha$ and $\beta$}
\begin{tabular}{l*{10}{c}}
\toprule
Estimator &  \multicolumn{3}{c}{25\% of untreated subjects} & \multicolumn{3}{c}{50\% of untreated subjects} & \multicolumn{3}{c}{75\%  of untreated subjects} \\
\midrule
{}   & bias & std & mse & bias & std & mse & bias & std & mse \\
$\hat{\Delta}_1$  & -0.0017 & 0.1110 & 0.0123 & -0.0037 & 0.0766 & 0.0059 & -0.0019 & 0.0971 & 0.0094\\
$\hat{\Delta}_2$  & 0.2345 &  0.0925 & 0.0635 & 0.2071 & 0.0751 & 0.0485 & 0.2387 & 0.0948 & 0.0660\\
$\hat{\Delta}_5$  & -0.0029 & 0.0901 & 0.0081 & -0.0030 & 0.0723 & 0.0052 & -0.0013 & 0.0900 & 0.0081\\
$\hat{\Delta}_6$  & 0.0174 & 0.1081 & 0.0120 & 0.0128 & 0.0757 & 0.0059 & 0.0205 & 0.0932 & 0.0091\\
$\hat{\Delta}_7$  & 0.2168 & 0.0912 & 0.0553 & 0.1801 & 0.0751 & 0.0381 & 0.2158 & 0.0914 & 0.0549\\
$\hat{\Delta}_9$  & 0.0135 & 0.0802 & 0.0066 & 0.0082 & 0.0694 & 0.0049 & 0.0141 & 0.0806 & 0.0067\\
\bottomrule
\end{tabular}
\end{table}

\begin{table}[!htbp]
\centering
\caption{Correctly specified $PS$ and $OR$ models, 50\% untreated subjects}
\begin{tabular}{l*{10}{c}}
\toprule
Estimator &  \multicolumn{3}{c}{$0^{\circ}$ between $\alpha$ and $\beta$} & \multicolumn{3}{c}{$45^{\circ}$ between $\alpha$ and $\beta$} & \multicolumn{3}{c}{$90^{\circ}$ between $\alpha$ and $\beta$} \\
\midrule
{}   & bias & std & mse & bias & std & mse & bias & std & mse \\
$\hat{\Delta}_1$  & 0.0002 & 0.0713 & 0.0051 & -0.0041 & 0.0759 & 0.0058 & -0.0031 & 0.0707 & 0.0050\\
$\hat{\Delta}_2$  & 0.2102 & 0.0707 & 0.0492 & 0.1448 & 0.0749 & 0.0266 & -0.0034 & 0.0694 & 0.0048\\
$\hat{\Delta}_3$  & 0.0011 & 0.0685 & 0.0047 & -0.0030 & 0.0722 & 0.0052 & -0.0028 & 0.0671 & 0.0045\\
$\hat{\Delta}_4$  & 0.4193 & 0.0729 & 0.1811 & 0.2937 & 0.0768 & 0.0922 & -0.0041 & 0.0728 & 0.0053\\
$\hat{\Delta}_5$  & 0.0004 & 0.0696 & 0.0048 & -0.0033 & 0.0735 & 0.0054 & -0.0032 & 0.0688 & 0.0047\\
$\hat{\Delta}_6$  & 0.0088 & 0.0714 & 0.0052 & 0.0005 & 0.0761 & 0.0058 & -0.0033 & 0.0706 & 0.0050\\
$\hat{\Delta}_7$  & 0.2526 & 0.0737 & 0.0693 & 0.1762 & 0.0797 & 0.0374 & -0.0046 & 0.0774 & 0.0060\\
$\hat{\Delta}_8$  & 0.0538 & 0.0652 & 0.0071 & 0.0322 & 0.0673 & 0.0056 & -0.0032 & 0.0630 & 0.0040\\
$\hat{\Delta}_9$  & 0.0175 & 0.0697 & 0.0052 & 0.0084 & 0.0724 & 0.0053 & -0.0038 & 0.0673 & 0.0045\\
\bottomrule
\end{tabular}
\end{table}

\begin{table}[!htbp]
\centering
\caption{Misspecified $OR$ model, 50\% untreated subjects}
\begin{tabular}{l*{10}{c}}
\toprule
Estimator &  \multicolumn{3}{c}{$0^{\circ}$ between $\alpha$ and $\beta$} & \multicolumn{3}{c}{$45^{\circ}$ between $\alpha$ and $\beta$} & \multicolumn{3}{c}{$90^{\circ}$ between $\alpha$ and $\beta$} \\
\midrule
{}   & bias & std & mse & bias & std & mse & bias & std & mse \\
$\hat{\Delta}_1$  & 0.0125 & 0.2632 & 0.0693 & -0.0099 & 0.2204 & 0.0486 & -0.0021 & 0.1715 & 0.0294\\
$\hat{\Delta}_3$  & 0.4617 & 0.4240 & 0.3928 & 0.3671 & 0.3397 & 0.2500 & 0.0649 & 0.2960 & 0.0917\\
$\hat{\Delta}_5$  & 0.1104 & 0.2016 & 0.0528 & 0.0841 & 0.1790 & 0.0391 & 0.0081 & 0.1681 & 0.0283\\
$\hat{\Delta}_6$  & 0.0088 & 0.0714 & 0.0052 & 0.0026 & 0.0776 & 0.0060 & -0.0024 & 0.0742 & 0.0055\\
$\hat{\Delta}_8$  & 0.0538 & 0.0652 & 0.0071 & 0.0490 & 0.0788 & 0.0086 & -0.0050 & 0.0825 & 0.0068\\
$\hat{\Delta}_9$  & 0.0175 & 0.0697 & 0.0052 & 0.0132 & 0.0913 & 0.0085 & -0.0067 & 0.1002 & 0.0101\\
\bottomrule
\end{tabular}
\end{table}

\begin{table}[!htbp]
\centering
\caption{Misspecified $PS$ model, 50\% untreated subjects}
\begin{tabular}{l*{10}{c}}
\toprule
Estimator &  \multicolumn{3}{c}{$0^{\circ}$ between $\alpha$ and $\beta$} & \multicolumn{3}{c}{$45^{\circ}$ between $\alpha$ and $\beta$} & \multicolumn{3}{c}{$90^{\circ}$ between $\alpha$ and $\beta$} \\
\midrule
{}   & bias & std & mse & bias & std & mse & bias & std & mse \\
$\hat{\Delta}_1$  & 0.0013 & 0.0716 & 0.0051 & -0.0037 & 0.0766 & 0.0059 & -0.0036 & 0.0773 & 0.0060\\
$\hat{\Delta}_2$  & 0.2957 & 0.0714 & 0.0925 & 0.2071 & 0.0751 & 0.0485 & -0.0040 & 0.0727 & 0.0053\\
$\hat{\Delta}_5$  & 0.0004 & 0.0696 & 0.0048 & -0.0030 & 0.0723 & 0.0052 & -0.0028 & 0.0673 & 0.0045\\
$\hat{\Delta}_6$  & 0.0267 & 0.0714 & 0.0058 & 0.0128 & 0.0757 & 0.0059 & -0.0039 & 0.0772 & 0.0060\\
$\hat{\Delta}_7$  & 0.2526 & 0.0737 & 0.0693 & 0.1801 & 0.0751 & 0.0381 & -0.0039 & 0.0696 & 0.0049\\
$\hat{\Delta}_9$  & 0.0175 & 0.0697 & 0.0052 & 0.0082 & 0.0694 & 0.0049 & -0.0032 & 0.0633 & 0.0040\\
\bottomrule
\end{tabular}
\end{table}

\newpage

\section{Discussion}
In this paper, the classical doubly robust estimation for ATE is revisited. We consider nine combinations of the estimated $PS$ model and $OR$ model under parametric, semiparametric and nonparametric model structures. When the models are correctly specified,  these nine estimators reach the same semiparametric efficiency bound. In other words, these nine estimators are all asymptotically efficient. Under the locally misspecified parametric $PS$ or $OR$ model which converges to the underlying parametric model, the estimators can still achieve the semiparametric efficiency bound. Further, when the $OR$ model is globally misspecified and  the $PS$ model is correctly specified, the asymptotic variance is always greater than or equal to the semiparametric efficiency bound. Yet, when the $PS$ model is globally misspecified and the $OR$ model is correctly specified, the situation becomes complicated. The asymptotic variance may not be always enlarged and in some cases, could be even smaller than the semiparametric efficiency bound. This phenomenon is interesting and worth a further study.

\section{Appendix}

The theorems in Section 3 and the following proofs are based on Assumptions \ref{uncounfounded}-\ref{general} in Section 2.

\subsection{Double robustness}
The double robustness of DR estimators has been proved in previous literature, see \cite{drproof} for example. When both of the PS model and OR model are correctly specified, the estimator is robust for sure. For completeness, we still provide brief calculations when one of $PS$ model or $OR$ model is misspecified to demonstrate the double robustness of these nine estimators.
\begin{lemma}
The estimators $\hat{\Delta}_k,\: k=1,\cdots,9$ are doubly robust when either (but not necessarily both) $PS$ model or $OR$ model is correctly specified.
\end{lemma}

\noindent \textbf{Proof} Recall that the true average treatment effect is defined as $\Delta=\theta_1-\theta_0$. For a doubly robust estimator of ATE, we also have $\hat{\Delta}=\hat{\theta}_1-\hat{\theta}_0$. We suppose $\hat{\Delta}$ converges to some quantities $\Delta^*=\theta^*_1-\theta^*_0$. If the double robustness of $\hat{\Delta}$ holds, we should have $\Delta^*=\Delta$. For simplicity, we only show the double robustness of $\hat{\theta}_1$. Note that $\hat{\theta}_1=\frac1n\sum\limits_{i=1}^{n}\left[\frac{d_iy_i}{\hat{p}(x_i)}+\left(1-\frac{d_i}{\hat{p}(x_i)}\right)\hat{m}_1(x_i)\right]$, where $\hat{p}(x)$ is an estimator of $p(x)$ and $\hat{m}_1(x)$ is an estimator of $m_1(x)$.

If $PS$ model is correctly specified and $OR$ model is misspecified, $\hat{\theta}_1$ estimates the population mean $\theta^*_1=\mathbb{E}\left\{\frac{DY}{p(X)}+\left(1-\frac{D}{p(X)}\right)m^*_1(X)\right\}$, where $p(x)$ is the true $PS$ model and $m^*_1(x)$ is a misspecified $OR$ model. We have
\begin{equation*}
\begin{split}
\theta^*_1&=\mathbb{E}\left\{\frac{DY}{p(X)}+\left(1-\frac{D}{p(X)}\right)m^*_1(X)\right\} \\
&=\mathbb{E}\left\{\mathbb{E}\left[\frac{DY}{p(X)}+\left(1-\frac{D}{p(X)}\right)m^*_1(X)\Big|X\right]\right\}\\
&=\mathbb{E}\left\{\frac{\mathbb{E}[D(DY(1)+(1-D)Y(0))|X]}{p(X)}+\left(1-\frac{\mathbb{E}(D|X)}{p(X)}\right)m^*_1(X)\right\}\\
&=\mathbb{E}\left\{\frac{p(X)\mathbb{E}[Y(1)|X]}{p(X)}+\left(1-\frac{p(X)}{p(X)}\right)m^*_1(X)\right\}\\
&=\mathbb{E}\left\{\mathbb{E}[Y(1)|X]\right\}
=\mathbb{E}[Y(1)]
=\theta_1,
\end{split}
\end{equation*}
which follows the conclusion that the double robustness of $\hat{\theta}_1$ holds.

If $OR$ model is correctly specified and $PS$ model is misspecified, $\hat{\theta}_1$ estimates the population mean $\theta^*_1=\mathbb{E}\left\{\frac{DY}{p^*(X)}+\left(1-\frac{D}{p^*(X)}\right)m_1(X)\right\}$, where $p^*(x)$ is a misspecified $PS$ model and $m_1(x)$ is the true $OR$ model. We have
\begin{equation*}
\begin{split}
\theta^*_1&=\mathbb{E}\left\{\frac{DY}{p^*(X)}+\left(1-\frac{D}{p^*(X)}\right)m_1(X)\right\} \\
&=\mathbb{E}\left\{\mathbb{E}\left[\frac{DY}{p^*(X)}+\left(1-\frac{D}{p^*(X)}\right)m_1(X)\Big|X\right]\right\}\\
&=\mathbb{E}\left\{\frac{\mathbb{E}[D(DY(1)+(1-D)Y(0))|X]}{p^*(X)}+\left(1-\frac{\mathbb{E}(D|X)}{p^*(X)}\right)m_1(X)\right\}\\
&=\mathbb{E}\left\{\frac{p(X)m_1(x)}{p^*(X)}+\left(1-\frac{p(X)}{p^*(X)}\right)m_1(X)\right\}\\
&=\mathbb{E}\left\{m_1(X)\right\}
=\mathbb{E}[Y(1)]
=\theta_1,
\end{split}
\end{equation*}
which follows the conclusion that the double robustness of $\hat{\theta}_1$ holds. Double robustness of $\hat{\theta}_0$ can be derived similarly. Therefore, the double robustness of aforementioned estimators holds.

\subsection{Proof of Main Results}
In this section, we derive the asymptotic distributions of the nine doubly robust estimators presented in (\ref{drdef}) and provide proofs for Theorems~1, 2 and 3. Note that we only focus on the cases when the doubly robustness holds with $\Delta^*=\theta^*_1-\theta^*_0=\Delta=\theta_1-\theta_0$ from Section 6.1. For simplicity, here we only present the details in deriving the form of $\sqrt{n}(\hat{\theta}_1-\theta_1)$. Similar method can be applied to derive $\sqrt{n}(\hat{\theta}_0-\theta_0)$. Then the form of $\sqrt{n}(\hat{\Delta}-\Delta)$ can be derived. Consequently, the asymptotic distribution of $\sqrt{n}(\hat{\Delta}-\Delta)$ can be obtained.

\subsubsection{Parametric $PS$ model and $OR$ model}
When $PS$ model and $OR$ models are both parametric, we have
\begin{eqnarray}\label{6.2.1}
\sqrt{n}(\hat{\theta}_1-\theta_1) &=& \sqrt{n}\left\{n^{-1}\sum_{i=1}^{n}\left[\frac{d_iy_i}{\Tilde{p}(x_i;\hat{\beta})}+\left(1-\frac{d_i}{\Tilde{p}(x_i;\hat{\beta})}\right)\Tilde{m}_1(x_i;\hat{\gamma}_1)\right]-\theta_1\right\} \nonumber\\
&=& \frac{1}{\sqrt{n}}\sum_{i=1}^{n}\left\{\frac{d_iy_i}{\Tilde{p}(x_i;\beta^*)}+\left(1-\frac{d_i}{\Tilde{p}(x_i;\beta^*)}\right)\Tilde{m}_1(x_i;\gamma^*_1)-\theta_1\right\} \nonumber\\
&&+\frac{1}{\sqrt{n}}\sum_{i=1}^{n}d_i\left(\frac{1}{\Tilde{p}(x_i;\hat{\beta})}-\frac{1}{\Tilde{p}(x_i;\beta^*)}\right)[y_i-m_1(x_i)] \nonumber\\
&&+\frac{1}{\sqrt{n}}\sum_{i=1}^{n}d_i\left(\frac{1}{\Tilde{p}(x_i;\hat{\beta})}-\frac{1}{\Tilde{p}(x_i;\beta^*)}\right)[m_1(x_i)-\Tilde{m}_1(x_i;\gamma^*_1)] \nonumber\\
&&+\frac{1}{\sqrt{n}}\sum_{i=1}^{n}\left(1-\frac{d_i}{p(x_i)}\right)[\Tilde{m}_1(x_i;\hat{\gamma}_1)-\Tilde{m}_1(x_i;\gamma^*_1)] \nonumber\\
&&+\frac{1}{\sqrt{n}}\sum_{i=1}^{n}d_i\left(\frac{1}{p(x_i)}-\frac{1}{\Tilde{p}(x_i;\beta^*)}\right)[\Tilde{m}_1(x_i;\hat{\gamma}_1)-\Tilde{m}_1(x_i;\gamma^*_1)] \nonumber\\
&&+\frac{1}{\sqrt{n}}\sum_{i=1}^{n}d_i\left(\frac{1}{\Tilde{p}(x_i;\beta^*)}-\frac{1}{\Tilde{p}(x_i;\hat{\beta})}\right)[\Tilde{m}_1(x_i;\hat{\gamma}_1)-\Tilde{m}_1(x_i;\gamma^*_1)] \nonumber\\
&:=& \sqrt{n}A_n + \sqrt{n}B_n + \sqrt{n}C_n + \sqrt{n}D_n + \sqrt{n}E_n + \sqrt{n}F_n.
\end{eqnarray}

\noindent Under Assumptions \ref{psassumption} and \ref{orassumption}, we have $\hat{\beta},\beta^* \in \Theta_{\beta}$, $\hat{\gamma}_0,\gamma^*_0 \in \Theta_{\gamma_0}$ and $\hat{\gamma}_1,\gamma^*_1 \in \Theta_{\gamma_1}$. Define $L(\hat{\beta};\beta^*)$ to be the line segment in $ \Theta_{\beta}$ between $\hat{\beta}$ and $\beta^*$, so $L(\hat{\beta};\beta^*)$ consists of vectors with the form $\hat{\beta} + t(\hat{\beta}-\beta^*),0\leq t \leq 1$. Similarly, we can define $L(\hat{\gamma}_0;\gamma^*_0)$ and $L(\hat{\gamma}_1;\gamma^*_1)$. According to the Mean Value Theorem, see \cite{mvt}, $\exists \: \bar{\beta}\in L(\hat{\beta};\beta^*)$ such that $\frac{1}{\Tilde{p}(x;\hat{\beta})}-\frac{1}{\Tilde{p}(x;\beta^*)}=(\hat{\beta}-\beta^*)^T\frac{\partial\frac{1}{\Tilde{p}(x;\beta)}}{\partial\beta}\big|_{\beta=\bar{\beta}}$. Similarly, we have $\exists \: \bar{\gamma}_0\in L(\hat{\gamma}_0;\gamma^*_0)$ such that $\Tilde{m}_0(x;\hat{\gamma}_0)-\Tilde{m}_0(x;\gamma^*_0)=(\hat{\gamma}_0-\gamma^*_0)^T\frac{\partial \Tilde{m}_0(x;\gamma_0)}{\partial\gamma_0}\big|_{\gamma_0=\bar{\gamma}_0}$ and $\exists \: \bar{\gamma}_1\in L(\hat{\gamma}_1;\gamma^*_1)$ such that $\Tilde{m}_1(x;\hat{\gamma}_1)-\Tilde{m}_1(x;\gamma^*_1)=(\hat{\gamma}_1-\gamma^*_1)^T\frac{\partial \Tilde{m}_1(x;\gamma_1)}{\partial\gamma_1}\big|_{\gamma_1=\bar{\gamma}_1}$. Then we have
\begin{equation*}
\begin{split}
\sqrt{n}B_n &=\frac{1}{\sqrt{n}}\sum_{i=1}^{n}d_i\left(\frac{1}{\Tilde{p}(x_i;\hat{\beta})}-\frac{1}{\Tilde{p}(x_i;\beta^*)}\right)[y_i-m_1(x_i)] \\
&=\frac{1}{\sqrt{n}}\sum_{i=1}^{n}(\hat{\beta}-\beta^*)^T \frac{\partial\frac{1}{\Tilde{p}(x_i;\beta)}}{\partial\beta}\bigg|_{\beta=\bar{\beta}} d_i[y_i-m_1(x_i)]\\
&=(\hat{\beta}-\beta^*)^T \left\{\sqrt{n}\mathbb{E}\left[\frac{\partial\frac{1}{\Tilde{p}(X;\beta)}}{\partial\beta}\bigg|_{\beta=\bar{\beta}}D(Y-m_1(X))\right]+O_p(1)\right\}\\
&=(\hat{\beta}-\beta^*)^T \left\{\sqrt{n}\mathbb{E}\left[\mathbb{E}\left(\frac{\partial\frac{1}{\Tilde{p}(X;\beta)}}{\partial\beta}\bigg|_{\beta=\bar{\beta}}D(Y-m_1(X))\bigg|X\right)\right]+O_p(1)\right\}\\
&=(\hat{\beta}-\beta^*)^T \left\{\sqrt{n}\mathbb{E}\left[\frac{\partial\frac{1}{\Tilde{p}(X;\beta)}}{\partial\beta}\Big|_{\beta=\bar{\beta}}\left(\mathbb{E}(DY|X)-\mathbb{E}(D|X)m_1(X)\right)\right]+O_p(1)\right\}\\
&=(\hat{\beta}-\beta^*)^T \left\{\sqrt{n}\mathbb{E}\left[\frac{\partial\frac{1}{\Tilde{p}(X;\beta)}}{\partial\beta}\bigg|_{\beta=\bar{\beta}}(p(X)m_1(X)-p(X)m_1(X))\right]+O_p(1)\right\}\\
&:=(\hat{\beta}-\beta^*)^T O_p(1),
\end{split}
\end{equation*}
where the third equality is obtained by applying the Central Limit Theorem. Due to the consistency of maximum likelihood estimation (see \cite{mle1} and \cite{mle2}), $\hat{\beta}$ converges to $\beta^*$ at rate $O(n^{-1/2})$, so $\sqrt{n}B_n=o_p(1)$. Similarly, we can obtain
\begin{equation*}
\begin{split}
&\sqrt{n}C_n=\sqrt{n}(\hat{\beta}-\beta^*)^T \mathbb{E}\left\{\frac{\partial\frac{1}{\Tilde{p}(X;\beta)}}{\partial\beta}\bigg|_{\beta=\bar{\beta}}p(X)[m_1(X)-\Tilde{m}_1(X;\gamma^*_1)]\right\}+o_p(1),\\
&\sqrt{n}E_n=\sqrt{n}(\hat{\gamma}_1-\gamma^*_1)^T \mathbb{E}\left\{ \frac{\partial \Tilde{m}_1(x;\gamma_1)}{\partial\gamma_1}\bigg|_{\gamma_1=\bar{\gamma}_1}\left(1-\frac{p(X)}{\Tilde{p}(X;\beta^*)}\right)\right\}+o_p(1),\\
&\sqrt{n}D_n=o_p(1), \sqrt{n}F_n=o_p(1).
\end{split}
\end{equation*}
\noindent Now we can consider different cases as follows.

\newcontent (a) \textbf{Correctly specified $PS$ model and $OR$ model}

In this case, we have $p(x)=\Tilde{p}(x;\beta_0)=\Tilde{p}(x;\beta^*)$,
$m_1(x)=\Tilde{m}_1(x;\gamma_{1,0})=\Tilde{m}_1(x;\gamma^*_1)$ and $m_0(x)=\Tilde{m}_0(x;\gamma_{0,0})=\Tilde{m}_0(x;\gamma^*_0)$. Then $\sqrt{n}C_n=o_p(1)$ and $\sqrt{n}E_n=o_p(1)$. Combining the terms in (\ref{6.2.1}), we have
\begin{equation*}
\begin{split}
\sqrt{n}(\hat{\theta}_1-\theta_1)=\frac{1}{\sqrt{n}}\sum_{i=1}^{n}\left\{\frac{d_iy_i}{p(x_i)}+\left(1-\frac{d_i}{p(x_i)}\right)m_1(x_i)-\theta_1\right\}+o_p(1).
\end{split}
\end{equation*}
Similarly, we can derive the form of $\sqrt{n}(\hat{\theta}_0-\theta_0)$. As a result,
\begin{equation*}
\begin{split}
&\sqrt{n}(\hat{\Delta}_1-\Delta) \\
&=\sqrt{n}[(\hat{\theta}_1-\hat{\theta}_0)-(\theta_1-\theta_0)]\\
&=\frac{1}{\sqrt{n}}\sum_{i=1}^{n}\left\{\frac{d_iy_i}{p(x_i)}+\left(1-\frac{d_i}{p(x_i)}\right)m_1(x_i)-\frac{(1-d_i)y_i}{1-p(x_i)}-\left(1-\frac{1-d_i}{1-p(x_i)}\right)m_0(x_i)-(\theta_1-\theta_0)\right\}+o_p(1)\\
&=\frac{1}{\sqrt{n}}\sum_{i=1}^{n}\Phi(x_i,y_i,d_i)+o_p(1).
\end{split}
\end{equation*}

\noindent Note that $\mathbb{E}\left\{\Phi(X,Y,D)\right\}=0$. We further assume $\mathbb{E}\left\{\Phi(X,Y,D)^2\right\} < \infty$. It follows from the Central Limit Theorem and Slutsky's Theorem that $\sqrt{n}(\hat{\Delta}_1-\Delta)$ converges in distribution to $N(0,\Sigma_1)$ with
\begin{equation*}
\Sigma_1 =\mathbb{E}\left\{\Phi(X,Y,D)^2\right\}=\mathbb{E}\left\{\frac{\mathrm{Var}[Y(1)|X]}{p(X)}+\frac{\mathrm{Var}[Y(0)|X]}{1-p(X)}+[m_1(X)-\theta_1-m_0(X)+\theta_0]^2\right\},
\end{equation*}
\noindent which is the same as the semiparametric efficiency bound discussed by \cite{efficiencybound}.

\newcontent (b) \textbf{Correctly specified $PS$ model and misspecified $OR$ model}

In this case, we have $p(x)=\Tilde{p}(x;\beta_0)=\Tilde{p}(x;\beta^*)$, $m_1(x)\neq \Tilde{m}_1(x;\gamma_{1,0})\neq\Tilde{m}_1(x;\gamma^*_1)$ and $m_0(x)\neq \Tilde{m}_0(x;\gamma_{0,0})\neq\Tilde{m}_0(x;\gamma^*_0)$. Thus, $\sqrt{n}E_n=o_p(1)$. Combining the terms in (\ref{6.2.1}), we have
\begin{equation*}
\begin{split}
&\sqrt{n}(\hat{\theta}_1-\theta_1)\\
&=\frac{1}{\sqrt{n}}\sum_{i=1}^{n}\left\{\frac{d_iy_i}{p(x_i)}+\left(1-\frac{d_i}{p(x_i)}\right)\Tilde{m}_1(x_i;\gamma^*_1)-\theta_1\right\}\\
&\quad+\sqrt{n}(\hat{\beta}-\beta_0)^T \mathbb{E}\left\{\frac{\partial\frac{1}{\Tilde{p}(X;\beta)}}{\partial\beta}\bigg|_{\beta=\bar{\beta}}p(X)[m_1(X)-\Tilde{m}_1(X;\gamma^*_1)]\right\}+o_p(1).
\end{split}
\end{equation*}
Similarly, we can derive the form of $\sqrt{n}(\hat{\theta}_0-\theta_0)$. As a result,
\begin{equation}  \label{localor}
\begin{split}
&\sqrt{n}(\hat{\Delta}_1-\Delta) \\
&=\sqrt{n}[(\hat{\theta}_1-\hat{\theta}_0)-(\theta_1-\theta_0)]\\
&=\frac{1}{\sqrt{n}}\sum_{i=1}^{n}\left\{\frac{d_iy_i}{p(x_i)}+\left(1-\frac{d_i}{p(x_i)}\right)\Tilde{m}_1(x_i;\gamma^*_1)-\frac{(1-d_i)y_i}{1-p(x_i)} -\left(1-\frac{1-d_i}{1-p(x_i)}\right)\Tilde{m}_0(x_i;\gamma^*_0)-(\theta_1-\theta_0)\right\} \\
&\quad+\sqrt{n}(\hat{\beta}-\beta_0)^T\mathbb{E}\left\{\frac{\partial\frac{1}{\Tilde{p}(X;\beta)}}{\partial\beta}\bigg|_{\beta=\bar{\beta}}p(X)[m_1(X)-\Tilde{m}_1(X;\gamma^*_1)]-\frac{\partial\frac{1}{1-\Tilde{p}(X;\beta)}}{\partial\beta}\bigg|_{\beta=\bar{\beta}}(1-p(X))[m_0(X)-\Tilde{m}_0(X;\gamma^*_0)]\right\} \\
&\quad+o_p(1).
\end{split}
\end{equation}

Now we consider the case of local misspecification of $OR$ models. According to the definition of locally misspecified $OR$ models in (\ref{localdef}), we have $m_1(x)=\Tilde{m}_1(x;\gamma_{1,0})+\delta_1\times s_1(x)$ and $m_0(x)=\Tilde{m}_0(x;\gamma_{0,0})+\delta_0 \times s_0(x)$. Recall that the underlying $OR$ models are defined as $E\left[Y(1)|X\right]=m_1(X)$ and $E\left[Y(0)|X\right]=m_1(X)$. That is, $y_i(1)=m_1(x_i)+\epsilon_{i(1)}$ with $i.i.d.$ random error $\epsilon_{i(1)},i=1,...,n_1$  from $N(0,\sigma^2_{(1)})$ and $y_j(0)=m_0(x_j)+\epsilon_{j(0)}$ with $i.i.d.$ random error $\epsilon_{j(0)},j=1,...,n_0$ from $N(0,\sigma^2_{(0)})$. We further assume $\sigma_{(1)}$ and $\sigma_{(0)}$ are nonzero constants. Note that $n_1+n_0=n$. The proposed $OR$ models are $\Tilde{m}_1(x;\gamma_{1})=x^T\gamma_{1}$ and $\Tilde{m}_0(x;\gamma_{0})=x^T\gamma_{0}$. Then we can obtain the loglikelihood function
$$l(\gamma_{1})=-n \log(\sqrt{2\pi}\sigma_{(1)})-\frac{1}{2\sigma^2_{(1)}}\sum^{n_1}_{i=1}[y_i(1)-x^T_i\gamma_{1}]^2,$$
and the score function
$$\frac{\partial l(\gamma_{1})}{\partial\gamma_{1}}=\frac{1}{\sigma^2_{(1)}}\sum^{n_1}_{i=1}x_i[y_i(1)-x^T_i\gamma_{1}].$$
\noindent Then we can solve $\mathbb{E}\left[\frac{\partial l(\gamma_{1})}{\partial\gamma_{1}}\right]=0$ for $\gamma_1$ and the resulting $\gamma_1$ is $\gamma^*_1$. That is,
\begin{equation*}
\begin{split}
0&=\mathbb{E}\left\{X\left[Y(1)-X^T\gamma^*_1\right]\right\}=\mathbb{E}\left\{X\left[\mathbb{E}\left(Y(1)|X\right)-X^T\gamma^*_1\right]\right\}\\
&=\mathbb{E}\left\{X\left[m_1(X)-X^T\gamma^*_1\right]\right\}=\mathbb{E}\left\{X\left[\Tilde{m}_1(X;\gamma_{1,0})+\delta_1\times s_1(X)-X^T\gamma^*_1\right]\right\}\\
&=\mathbb{E}\left[XX^T\right](\gamma_{1,0}-\gamma^*_1)+\delta_1\mathbb{E}\left[Xs_1(X)\right],
\end{split}
\end{equation*}
\noindent which leads to $\gamma^*_1=\gamma_{1,0}+O(\delta_1)$ given that $\mathbb{E}\left[XX^T\right]$, $\mathbb{E}\left[Xs_1(X)\right]$ are bounded away from zero and infinity, and $\mathbb{E}\left[XX^T\right]$ is invertible. Similarly, we can obtain $\gamma^*_0-\gamma_{0,0}=O(\delta_0)$. Under Assumption \ref{orassumption}, we have $\gamma_{0,0},\gamma^*_0 \in \Theta_{\gamma_0}$ and $\gamma_{1,0},\gamma^*_1 \in \Theta_{\gamma_1}$. Define $L(\gamma_{0,0};\gamma^*_0)$ to be the line segment in $\Theta_{\gamma_0}$ between $\gamma_{0,0}$ and $\gamma^*_0$, so $L(\gamma_{0,0};\gamma^*_0)$ consists of vectors with the form $\gamma_{0,0}+ t(\gamma_{0,0}-\gamma^*_0),0\leq t \leq 1$. Similarly, we can define $L(\gamma_{1,0};\gamma^*_1)$. According to the Mean Value Theorem, see \cite{mvt}, $\exists \: \Tilde{\gamma}_0\in L(\gamma_{0,0};\gamma^*_0)$ such that $\Tilde{m}_0(x;\gamma_{0,0})-\Tilde{m}_0(x;\gamma^*_0)=(\gamma_{0,0}-\gamma^*_0)^T\frac{\partial \Tilde{m}_0(x;\gamma_0)}{\partial\gamma_0}\big|_{\gamma_0=\Tilde{\gamma}_0}$ and $\exists \: \Tilde{\gamma}_1\in L(\gamma_{1,0};\gamma^*_1)$ such that $\Tilde{m}_1(x;\gamma_{1,0})-\Tilde{m}_1(x;\gamma^*_1)=(\gamma_{1,0}-\gamma^*_1)^T\frac{\partial \Tilde{m}_1(x;\gamma_1)}{\partial\gamma_1}\big|_{\gamma_1=\Tilde{\gamma}_1}$. Then we have \begin{eqnarray*}m_1(x)-\Tilde{m}_1(x;\gamma^*_1)&=&
\Tilde{m}_1(x;\gamma_{1,0})+\delta_1\times s_1(x)-\Tilde{m}_1(x;\gamma^*_1)\\
&=&\delta_1\times s_1(x)+ (\gamma_{1,0}-\gamma^*_1)^T\frac{\partial \Tilde{m}_1(x;\gamma_1)}{\partial\gamma_1}\Big|_{\gamma_1=\Tilde{\gamma}_1}=O_p(\delta_1),\\ m_0(x)-\Tilde{m}_0(x;\gamma^*_0)&=&\Tilde{m}_0(x;\gamma_{0,0})+\delta_0\times s_0(x)-\Tilde{m}_0(x;\gamma^*_0)\\
&=&\delta_0\times s_0(x)+(\gamma_{0,0}-\gamma^*_0)^T\frac{\partial \Tilde{m}_0(x;\gamma_0)}{\partial\gamma_0}\Big|_{\gamma_0
=\Tilde{\gamma}_0}=O_p(\delta_0).
\end{eqnarray*}
\noindent So the second term in (\ref{localor}) can be written as
$\sqrt{n}(\hat{\beta}-\beta_0)[O_p(\delta_1)+O_p(\delta_0)]$
which converges to 0 in probability as $\delta_1\to0$ and $\delta_0\to0$ by Slutsky's theorem due to the consistency of MLE. Let $\Phi(x_i,y_i,d_i):=\frac{d_iy_i}{p(x_i)}+\left(1-\frac{d_i}{p(x_i)}\right)\Tilde{m}_1(x_i;\gamma^*_1)-\frac{(1-d_i)y_i}{1-p(x_i)} -\left(1-\frac{1-d_i}{1-p(x_i)}\right)\Tilde{m}_0(x_i;\gamma^*_0)-(\theta_1-\theta_0)$. Note that $\mathbb{E}\left\{\Phi(X,Y,D)\right\}=0$. We further assume $\mathbb{E}\left\{\Phi(X,Y,D)^2\right\} < \infty$. It follows from the Central Limit Theorem and Slutsky's Theorem that the first term converges in distribution to $N(0, \Sigma)$ with $$\Sigma=\mathrm{Var}\left[\frac{DY}{p(X)}+\left(1-\frac{D}{p(X)}\right)\Tilde{m}_1(X;\gamma^*_1)-\frac{(1-D)Y}{1-p(X)} -\left(1-\frac{1-D}{1-p(X)}\right)\Tilde{m}_0(X;\gamma^*_0)-(\theta_1-\theta_0)\right].$$ Note that $\Sigma$ converges to $\Sigma_1$ as $\delta_1\to0$ and $\delta_0\to0$. Consequently, the asymptotic variance of $\hat{\Delta}_1$ converges to $\Sigma_1$.

In the following, we consider the case when $OR$ models are globally misspecified. Since $PS$ model is correctly specified, we have $PS$ model $\Tilde{p}(x;\beta_0)=p(x)=\frac{\exp(x^T\beta_0)}{1+\exp(x^T\beta_0)}$. Then we can obtain
$$log\;likelihood\;function: l(\beta_0)=\sum^n_{i=1}\left[d_ix^T_i\beta_0-\log(1+e^{x^T_i\beta_0})\right];$$
$$score\; vector: S(\beta_0)=\frac{\partial l(\beta_0)}{\partial\beta_0}=\sum^n_{i=1}x_i[d_i-p(x_i)]=\sum^n_{i=1}S(\beta_0;x_i,d_i);$$
$$observed\;information\;matrix: -\frac{\partial^2 l(\beta_0)}{\partial\beta_0\partial\beta^T_0}.$$

\noindent By the weak law of large numbers, $-\frac{1}{n}\frac{\partial^2 l(\beta_0)}{\partial\beta_0\partial\beta^T_0}$ converges in probability to the Fisher information matrix $I(\beta_0)$. Applying Taylor series expansion, we have
\begin{equation*}
\begin{split}
0=S(\hat{\beta})&=S(\beta_0)+\frac{\partial^2 l(\beta_0)}{\partial\beta_0\partial\beta^T_0}(\hat{\beta}-\beta_0)+o_p(n^{-1/2})\\
&=n^{-1}S(\beta_0)-\left[-\frac{1}{n}\frac{\partial^2 l(\beta_0)}{\partial\beta_0\partial\beta^T_0}\right](\hat{\beta}-\beta_0)+o_p(n^{-1/2})\\
&=n^{-1}S(\beta_0)-I(\beta_0)(\hat{\beta}-\beta_0)+o_p(1).
\end{split}
\end{equation*}

\noindent Hence,
$$\sqrt{n}(\hat{\beta}-\beta_0)=I^{-1}(\beta_0)\frac{1}{\sqrt{n}}S(\beta_0)+o_p(1).$$

\noindent Therefore, we can further write
$$\sqrt{n}(\hat{\beta}-\beta_0)=\frac{1}{\sqrt{n}}\sum^n_{i=1}I^{-1}(\beta_0)S(\beta_0;x_i)+o_p(1)=\frac{1}{\sqrt{n}}\sum^n_{i=1}I^{-1}(\beta_0) x_i\left[d_i-p(x_i)\right]+o_p(1).$$

\noindent As a result,
\begin{equation*}
\begin{split}
&\sqrt{n}(\hat{\Delta}_1-\Delta) \\
&=\frac{1}{\sqrt{n}}\sum_{i=1}^{n}\bigg\{\frac{d_iy_i}{p(x_i)}+\left(1-\frac{d_i}{p(x_i)}\right)\Tilde{m}_1(x_i;\gamma^*_1)-\frac{(1-d_i)y_i}{1-p(x_i)} -\left(1-\frac{1-d_i}{1-p(x_i)}\right)\Tilde{m}_0(x_i;\gamma^*_0)-(\theta_1-\theta_0) \\
&\qquad+w(x_i)[p(x_i)-d_i] \bigg\}+o_p(1)\\
&=\frac{1}{\sqrt{n}}\sum_{i=1}^{n}\Phi(x_i,y_i,d_i)+o_p(1),
\end{split}
\end{equation*}
where $w(x_i)=\left(I^{-1}(\beta_0)x_i\right)^T \mathbb{E}\Big\{p(X)\frac{\partial\frac{1}{\Tilde{p}(X;\beta)}}{\partial\beta}\Big|_{\beta=\bar{\beta}}[\Tilde{m}_1(X;\gamma^*_1)-m_1(X)]-[1-p(X)]\frac{\partial\frac{1}{1-\Tilde{p}(X;\beta)}}{\partial\beta}\Big|_{\beta=\bar{\beta}}[\Tilde{m}_0(X;\gamma^*_0)-m_0(X)]\Big\}$.

\noindent Note that $\mathbb{E}\left\{\Phi(X,Y,D)\right\}=0$. We further assume $\mathbb{E}\left\{\Phi(X,Y,D)^2\right\} < \infty$. It follows from the Central Limit Theorem and Slutsky's Theorem that $\sqrt{n}(\hat{\Delta}_1-\Delta) $ converges in distribution to $N(0,\Sigma)$ with
\begin{equation*}
\begin{split}
\Sigma &=\mathbb{E}\left\{\Phi(X,Y,D)^2\right\}\\
&=\Sigma_1+\mathbb{E}\bigg\{\Big[\sqrt{\frac{1}{p(X)}-1}[\Tilde{m}_1(X;\gamma^*_1)-m_1(X)]+\sqrt{\frac{1}{1-p(X)}-1}[\Tilde{m}_0(X;\gamma^*_0)-m_0(X)]\\
&\quad+\sqrt{p(X)(1-p(X))}w(X)\Big]^2\bigg\} \geq \Sigma_1.
\end{split}
\end{equation*}

\noindent Therefore, the asymptotic variance of $\hat{\Delta}_1$ is enlarged compared to the semiparametric efficiency bound $\Sigma_1$ when $OR$ models are globally misspecified.

\newcontent (c) \textbf{Misspecified $PS$ model and correctly specified $OR$ model}

In this case, we have $p(x)\neq\Tilde{p}(x;\beta_0)\neq\Tilde{p}(x;\beta^*)$, $m_1(x)=\Tilde{m}_1(x;\gamma_{1,0})=\Tilde{m}_1(x;\gamma^*_1)$ and $m_0(x)=\Tilde{m}_0(x;\gamma_{0,0})=\Tilde{m}_0(x;\gamma^*_0)$. So $\sqrt{n}C_n=o_p(1)$. Combining the terms in (\ref{6.2.1}), we have
\begin{equation*}
\begin{split}
&\sqrt{n}(\hat{\theta}_1-\theta_1)\\
&=\frac{1}{\sqrt{n}}\sum_{i=1}^{n}\left\{\frac{d_iy_i}{\Tilde{p}(x_i;\beta^*)}+\left(1-\frac{d_i}{\Tilde{p}(x_i;\beta^*)}\right)m_1(x_i)-\theta_1\right\}\\
&\quad+\sqrt{n}(\hat{\gamma}_1-\gamma_{1,0})^T \mathbb{E}\left\{\frac{\partial \Tilde{m}_1(X;\gamma_1)}{\partial \gamma_1}\bigg|_{\gamma_1=\bar{\gamma}_1}\left(1-\frac{p(X)}{\Tilde{p}(X;\beta^*)}\right)\right\}\\
&\quad+o_p(1).
\end{split}
\end{equation*}
Similarly, we can derive the form of $\sqrt{n}(\hat{\theta}_0-\theta_0)$. As a result,
\begin{equation} \label{localps}
\begin{split}
&\sqrt{n}(\hat{\Delta}_1-\Delta) \\
&=\sqrt{n}[(\hat{\theta}_1-\hat{\theta}_0)-(\theta_1-\theta_0)]\\
&=\frac{1}{\sqrt{n}}\sum_{i=1}^{n}\bigg\{\frac{d_iy_i}{\Tilde{p}(x_i;\beta^*)}+\left(1-\frac{d_i}{\Tilde{p}(x_i;\beta^*)}\right)m_1(x_i)-\frac{(1-d_i)y_i}{1-\Tilde{p}(x_i;\beta^*)}-\left(1-\frac{1-d_i}{1-\Tilde{p}(x_i;\beta^*)}\right)m_0(x_i)\\
&\qquad -(\theta_1-\theta_0)\bigg\} \\
&\quad+\sqrt{n}(\hat{\gamma}_1-\gamma_{1,0})^T \mathbb{E}\left\{\frac{\partial \Tilde{m}_1(X;\gamma_1)}{\partial \gamma_1}\bigg|_{\gamma_1=\bar{\gamma}_1}\left(1-\frac{p(X)}{\Tilde{p}(X;\beta^*)}\right)\right\}\\
&\quad+\sqrt{n}(\hat{\gamma}_0-\gamma_{0,0})^T \mathbb{E}\left\{\frac{\partial \Tilde{m}_0(X;\gamma_0)}{\partial \gamma_0}\bigg|_{\gamma_0=\bar{\gamma}_0}\left(\frac{1-p(X)}{1-\Tilde{p}(X;\beta^*)}-1\right)\right\} \\
&\quad+o_p(1).
\end{split}
\end{equation}

Now we consider the case of local misspecification of $PS$ model. Recall that the underlying $PS$ model is defined as $P(D=1|X)=p(X)$. The locally misspecified $PS$ model is defined as $p(x)=\Tilde{p}(x;\beta_0)(1+\delta \times s(x))$ in (\ref{localdef}). The proposed $PS$ model is $\Tilde{p}(x;\beta)=\frac{\exp(x^T\beta)}{1+\exp(x^T\beta)}$. Then we can obtain the loglikelihood function
$$l(\beta)=\sum^n_{i=1}\left[d_ix^T_i\beta-\log(1+e^{x^T_i\beta})\right],$$
and the score function
$$\frac{\partial l(\beta)}{\partial\beta}=\sum^n_{i=1}x_i[d_i-\Tilde{p}(x_i;\beta)].$$
\noindent Then we can solve $\mathbb{E}\left[\frac{\partial l(\beta)}{\partial\beta}\right]=0$ for $\beta$ and the resulting $\beta$ is $\beta^*$. That is,
\begin{equation*}
\begin{split}
0&=\mathbb{E}\left\{X\left[D-\Tilde{p}(X;\beta^*)\right]\right\}\\
&=\mathbb{E}\left\{X\left[\mathbb{E}(D|X)-\Tilde{p}(X;\beta^*)\right]\right\}\\
&=\mathbb{E}\left\{X\left[p(X)-\Tilde{p}(X;\beta^*)\right]\right\}\\
&=\mathbb{E}\left\{X\left[\Tilde{p}(X;\beta_0)(1+\delta \times s(X))-\Tilde{p}(X;\beta^*)\right]\right\}\\
&=\mathbb{E}\left\{X\left[\frac{\partial \Tilde{p}(x;\beta)}{\partial\beta}\Big|_{\beta=\Tilde{\beta}}\right]^T\right\}(\beta_0-\beta^*)+\delta\mathbb{E}\left\{X\Tilde{p}(X;\beta_0)s(X)\right\},
\end{split}
\end{equation*}
\noindent where the last equality is obtained by Mean Value Theorem. Under Assumption \ref{psassumption}, we have $\beta_0,\beta^*\in \Theta_{\beta}$. Define $L(\beta_0;\beta^*)$ to be the line segment in $\Theta_{\beta}$ between $\beta_0$ and $\beta^*$, so $L(\beta_0;\beta^*)$ consists of vectors with the form $\beta_0+ t(\beta_0-\beta^*),0\leq t \leq 1$. According to the Mean Value Theorem, see \cite{mvt}, $\exists \: \Tilde{\beta}\in L(\beta_0;\beta^*)$ such that $\Tilde{p}(x;\beta_0)-\Tilde{p}(x;\beta^*)=(\beta_0-\beta^*)^T\frac{\partial \Tilde{p}(x;\beta)}{\partial\beta}\big|_{\beta=\Tilde{\beta}}$. Therefore, we can obtain $\beta^*-\beta_0=O(\delta)$ given that $\mathbb{E}\left\{X\left[\frac{\partial \Tilde{p}(x;\beta)}{\partial\beta}\big|_{\beta=\Tilde{\beta}}\right]^T\right\}$, $\mathbb{E}\left\{X\Tilde{p}(X;\beta_0)s(X)\right\}$ are bounded away from zero and infinity, and $\mathbb{E}\left\{X\left[\frac{\partial \Tilde{p}(x;\beta)}{\partial\beta}\big|_{\beta=\Tilde{\beta}}\right]^T\right\}$ is invertible. By Taylor series expansion, we have $$\Tilde{p}(x;\beta^*)-p(x)=\Tilde{p}(x;\beta^*)-\Tilde{p}(x;\beta_0)(1+\delta \times s(x))=(\beta^*-\beta_0)^T\frac{\partial \Tilde{p}(x;\beta)}{\partial\beta}\Big|_{\beta=\Tilde{\beta}}-\delta \times s(x)\Tilde{p}(x;\beta_0)=O_p(\delta).$$
\noindent So the second term and the third term in (\ref{localps}) can be written as $\sqrt{n}(\hat{\gamma}_1-\gamma_{1,0})O_p(\delta)$ and $\sqrt{n}(\hat{\gamma}_0-\gamma_{0,0})O_p(\delta)$ respectively. These two terms converge to 0 in probability as $\delta\to0$ again by Slutsky's theorem due to the consistency of MLE. Let $\Phi(x_i,y_i,d_i):=\frac{d_iy_i}{\Tilde{p}(x_i;\beta^*)}+\left(1-\frac{d_i}{\Tilde{p}(x_i;\beta^*)}\right)m_1(x_i)-\frac{(1-d_i)y_i}{1-\Tilde{p}(x_i;\beta^*)}-\left(1-\frac{1-d_i}{1-\Tilde{p}(x_i;\beta^*)}\right)m_0(x_i)-(\theta_1-\theta_0)$. Note that $\mathbb{E}\left\{\Phi(X,Y,D)\right\}=0$. We further assume $\mathbb{E}\left\{\Phi(X,Y,D)^2\right\} < \infty$. It follows from the Central Limit Theorem and Slutsky's Theorem that the first term converges in distribution to $N(0, \Sigma)$ with $$\Sigma=\mathrm{Var}\left[\frac{DY}{\Tilde{p}(X;\beta^*)}+\left(1-\frac{D}{\Tilde{p}(X;\beta^*)}\right)m_1(X)-\frac{(1-D)Y}{1-\Tilde{p}(X;\beta^*)} -\left(1-\frac{1-D}{1-\Tilde{p}(X;\beta^*)}\right)m_0(X)-(\theta_1-\theta_0)\right].$$ Note that $\Sigma$ converges to $\Sigma_1$ as $\delta\to0$. Consequently, the asymptotic variance of $\hat{\Delta}_1$ converges to $\Sigma_1$.

In the following, we consider the case when $PS$ model is globally misspecified. Since $OR$ models are correctly specified, we have $m_1(x)=\Tilde{m}_1(x;\gamma_{1,0})=x^T\gamma_{1,0}$ and $m_0(x)=\Tilde{m}_0(x;\gamma_{0,0})=x^T\gamma_{0,0}$. That is, $y_i(1)=x^T_i\gamma_{1,0}+\epsilon_{i(1)}$ with $i.i.d.$ random error $\epsilon_{i(1)},i=1,...,n_1$  from $N(0,\sigma^2_{(1)})$ and $y_i(0)=x^T_i\gamma_{0,0}+\epsilon_{i(0)}$ with $i.i.d.$ random error $\epsilon_{i(0)},i=1,...,n_0$ from $N(0,\sigma^2_{(0)})$. We further assume $\sigma_{(1)}$ and $\sigma_{(0)}$ are nonzero constants. Note that $n_1+n_0=n$. Then we can obtain

\newcontent $log\;likelihood\;functions:$
$$l(\gamma_{1,0})=-n \log(\sqrt{2\pi}\sigma_{(1)})-\frac{1}{2\sigma^2_{(1)}}\sum^{n_1}_{i=1}[y_i(1)-x^T_i\gamma_{1,0}]^2,\;l(\gamma_{0,0})=-n \log(\sqrt{2\pi}\sigma_{(0)})-\frac{1}{2\sigma^2_{(0)}}\sum^{n_0}_{i=1}[y_i(0)-x^T_i\gamma_{0,0}]^2;$$
$score\; vectors$:
$$S(\gamma_{1,0})=\frac{\partial l(\gamma_{1,0})}{\partial\gamma_{1,0}}=\frac{1}{\sigma^2_{(1)}}\sum^{n_1}_{i=1}x_i[y_i(1)-x^T_i\gamma_{1,0}]=\frac{1}{\sigma^2_{(1)}}\sum^{n}_{i=1}x_id_i[y_i-m_1(x_i)]=\sum^{n}_{i=1}S(\gamma_{1,0};x_i,y_i,d_i),$$
$$S(\gamma_{0,0})=\frac{\partial l(\gamma_{0,0})}{\partial\gamma_{0,0}}=\frac{1}{\sigma^2_{(0)}}\sum^{n_0}_{i=1}x_i[y_i(0)-x^T_i\gamma_{0,0}]=\frac{1}{\sigma^2_{(0)}}\sum^{n}_{i=1}x_i(1-d_i)[y_i-m_0(x_i)]=\sum^{n}_{i=1}S(\gamma_{0,0};x_i,y_i,d_i);$$
$observed\;information\; matrices$:
$$-\frac{\partial^2 l(\gamma_{1,0})}{\partial\gamma_{1,0}\partial\gamma^T_{1,0}} \; and \; -\frac{\partial^2 l(\gamma_{0,0})}{\partial\gamma_{0,0}\partial\gamma^T_{0,0}}.$$

\noindent By the weak law of large numbers, $-\frac{1}{n}\frac{\partial^2 l(\gamma_{1,0})}{\partial\gamma_{1,0}\partial\gamma^T_{1,0}}$ converges in probability to the Fisher information matrix $I(\gamma_{1,0})$ and $-\frac{1}{n}\frac{\partial^2 l(\gamma_{0,0})}{\partial\gamma_{0,0}\partial\gamma^T_{0,0}}$ converges in probability to the Fisher information matrix $I(\gamma_{0,0})$. Applying Taylor series expansion, we have
\begin{equation*}
\begin{split}
0=S(\hat{\gamma}_1)&=S(\gamma_{1,0})+\frac{\partial^2 l(\gamma_{1,0})}{\partial\gamma_{1,0}\partial\gamma^T_{1,0}}(\hat{\gamma}_1-\gamma_{1,0})+o_p(n^{-1/2})\\
&=n^{-1}S(\gamma_{1,0})-\left[-\frac{1}{n}\frac{\partial^2 l(\gamma_{1,0})}{\partial\gamma_{1,0}\partial\gamma^T_{1,0}}\right](\hat{\gamma}_1-\gamma_{1,0})+o_p(n^{-1/2})\\
&=n^{-1}S(\gamma_{1,0})-I(\gamma_{1,0})(\hat{\gamma}_1-\gamma_{1,0})+o_p(1),
\end{split}
\end{equation*}
\begin{equation*}
\begin{split}
0=S(\hat{\gamma}_0)&=S(\gamma_{0,0})+\frac{\partial^2 l(\gamma_{0,0})}{\partial\gamma_{0,0}\partial\gamma^T_{0,0}}(\hat{\gamma}_0-\gamma_{0,0})+o_p(n^{-1/2})\\
&=n^{-1}S(\gamma_{0,0})-\left[-\frac{1}{n}\frac{\partial^2 l(\gamma_{0,0})}{\partial\gamma_{0,0}\partial\gamma^T_{0,0}}\right](\hat{\gamma}_0-\gamma_{0,0})+o_p(n^{-1/2})\\
&=n^{-1}S(\gamma_{0,0})-I(\gamma_{0,0})(\hat{\gamma}_0-\gamma_{0,0})+o_p(1).
\end{split}
\end{equation*}
Hence,
$$\sqrt{n}(\hat{\gamma}_1-\gamma_{1,0})=I^{-1}(\gamma_{1,0})\frac{1}{\sqrt{n}}S(\gamma_{1,0})+o_p(1),$$
$$\sqrt{n}(\hat{\gamma}_0-\gamma_{0,0})=I^{-1}(\gamma_{0,0})\frac{1}{\sqrt{n}}S(\gamma_{0,0})+o_p(1).$$

\noindent Therefore, we can further write
\begin{equation*}
\begin{split}
\sqrt{n}(\hat{\gamma}_1-\gamma_{1,0})&=\frac{1}{\sqrt{n}}\sum^{n}_{i=1}I^{-1}(\gamma_{1,0})S(\gamma_{1,0};x_i,y_i,d_i)+o_p(1)\\
&=\frac{1}{\sqrt{n}}\sum^{n}_{i=1}\frac{1}{\sigma^2_{(1)}}I^{-1}(\gamma_{1,0})x_id_i[y_i-m_1(x_i)]+o_p(1),\\
\sqrt{n}(\hat{\gamma}_0-\gamma_{0,0})&=\frac{1}{\sqrt{n}}\sum^{n}_{i=1}I^{-1}(\gamma_{0,0})S(\gamma_{0,0};x_i,y_i,d_i)+o_p(1)\\
&=\frac{1}{\sqrt{n}}\sum^{n}_{i=1}\frac{1}{\sigma^2_{(0)}}I^{-1}(\gamma_{0,0})x_i(1-d_i)[y_i-m_0(x_i)]+o_p(1).
\end{split}
\end{equation*}

\noindent As a result,
\begin{equation*}
\begin{split}
&\sqrt{n}(\hat{\Delta}_1-\Delta) \\
&=\frac{1}{\sqrt{n}}\sum_{i=1}^{n}\bigg\{\frac{d_iy_i}{\Tilde{p}(x_i;\beta^*)}+\left(1-\frac{d_i}{\Tilde{p}(x_i;\beta^*)}\right)m_1(x_i)-\frac{(1-d_i)y_i}{1-\Tilde{p}(x_i;\beta^*)}-\left(1-\frac{1-d_i}{1-\Tilde{p}(x_i;\beta^*)}\right)m_0(x_i)-(\theta_1-\theta_0) \\
&\qquad+w_1(x_i)d_i[y_i-m_1(x_i)]+w_0(x_i)(1-d_i)[y_i-m_0(x_i)]\bigg\}+o_p(1)\\
&:=\frac{1}{\sqrt{n}}\sum_{i=1}^{n}\Phi(x_i,y_i,d_i)+o_p(1),
\end{split}
\end{equation*}
where
\begin{equation*}
\begin{split}
&w_1(x_i)= \frac{1}{\sigma^2_{(1)}}\left(I^{-1}(\gamma_{1,0})x_i\right)^T \mathbb{E}\left\{\frac{\partial \Tilde{m}_1(X;\gamma_1)}{\partial \gamma_1}\Big|_{\gamma_1=\bar{\gamma}_1}\left(1-\frac{p(X)}{\Tilde{p}(X;\beta^*)}\right)\right\},\\
&w_0(x_i)=\frac{1}{\sigma^2_{(0)})}\left(I^{-1}(\gamma_{0,0})x_i\right)^T \mathbb{E} \left\{ \frac{\partial \Tilde{m}_0(X;\gamma_0)}{\partial \gamma_0}\Big|_{\gamma_0=\bar{\gamma}_0}\left(\frac{1-p(X)}{1-\Tilde{p}(X;\beta^*)}-1\right) \right\}.
\end{split}
\end{equation*}

\noindent Note that $\mathbb{E}\left\{\Phi(X,Y,D)\right\}=0$. We further assume $\mathbb{E}\left\{\Phi(X,Y,D)^2\right\} < \infty$. It follows from the Central Limit Theorem and Slutsky's Theorem that $\sqrt{n}(\hat{\Delta}_1-\Delta) $ converges in distribution to $N(0,\Sigma)$ with
\begin{equation*}
\begin{split}
\Sigma &=\mathbb{E}\left\{\Phi(X,Y,D)^2\right\}\\
&=\Sigma_1 +\mathbb{E}\left\{\frac{1}{p(X)}\mathrm{Var}[Y(1)|X]\left[\left(\frac{p(X)}{\Tilde{p}(X;\beta^*)}+w_1(X)p(X)\right)^2-1\right]\right\}\\
&\quad+\mathbb{E}\left\{\frac{1}{1-p(X)}\mathrm{Var}[Y(0)|X]\left[\left(\frac{1-p(X)}{1-\Tilde{p}(X;\beta^*)}+w_0(X)(1-p(X))\right)^2-1\right]\right\}.
\end{split}
\end{equation*}
\noindent Therefore, the asymptotic variance of $\hat{\Delta}_1$ is is not necessarily enlarged compared to the semiparametric efficiency bound $\Sigma_1$ when $PS$ model is globally misspecified.

\subsubsection{Parametric $PS$ model and nonparametric $OR$ model}

When $PS$ model is parametric and $OR$ model is nonparametric, we have
\begin{eqnarray}\label{6.2.2}
\sqrt{n}(\hat{\theta}_1-\theta_1) &=& \sqrt{n}\left\{n^{-1}\sum_{i=1}^{n}\left[\frac{d_iy_i}{\Tilde{p}(x_i;\hat{\beta})}+\left(1-\frac{d_i}{\Tilde{p}(x_i;\hat{\beta})}\right)\hat{m}_1(x_i)\right]-\theta_1\right\} \nonumber\\
&=&\frac{1}{\sqrt{n}}\sum_{i=1}^{n}\left\{\frac{d_iy_i}{\Tilde{p}(x_i;\beta^*)}+\left(1-\frac{d_i}{\Tilde{p}(x_i;\beta^*)}\right)m_1(x_i)-\theta_1\right\} \nonumber\\
&&+\frac{1}{\sqrt{n}}\sum_{i=1}^{n}d_i\left(\frac{1}{\Tilde{p}(x_i;\hat{\beta})}-\frac{1}{\Tilde{p}(x_i;\beta^*)}\right)[y_i-m_1(x_i)] \nonumber\\
&&+\frac{1}{\sqrt{n}}\sum_{i=1}^{n}\left(1-\frac{d_i}{\Tilde{p}(x_i;\beta^*)}\right)[\hat{m}_1(x_i)-m_1(x_i)] \nonumber\\
&&+\frac{1}{\sqrt{n}}\sum_{i=1}^{n}d_i\left(\frac{1}{\Tilde{p}(x_i;\beta^*)}-\frac{1}{\Tilde{p}(x_i;\hat{\beta})}\right)[\hat{m}_1(x_i)-m_1(x_i)] \nonumber\\
&:=&\sqrt{n}A_n+\sqrt{n}B_n+\sqrt{n}C_n+\sqrt{n}D_n.
\end{eqnarray}

\noindent In Section 6.2.1, we already showed $\sqrt{n}B_n=o_p(1)$.

\begin{equation*}
\begin{split}
\sqrt{n}C_n &=\frac{1}{\sqrt{n}}\sum_{i=1}^{n}\left(1-\frac{d_i}{\Tilde{p}(x_i;\beta^*)}\right)\frac{\sum_{j=1}^{n}d_j(y_j-m_1(x_i))\Tilde{K}_{\Tilde{h}_{m_1}}(x_i,x_j)}{\sum_{j=1}^{n}d_j\Tilde{K}_{\Tilde{h}_{m_1}}(x_i,x_j)}\\
&=\frac{1}{\sqrt{n}}\sum_{i=1}^{n}\left(1-\frac{d_i}{\Tilde{p}(x_i;\beta^*)}\right)\frac{\sum_{j=1}^{n}d_j(y_j-m_1(x_i))\Tilde{K}_{\Tilde{h}_{m_1}}(x_i,x_j)}{n\frac{\sum_{j=1}^{n}d_j\Tilde{K}_{\Tilde{h}_{m_1}}(x_i,x_j)}{\sum_{j=1}^{n}\Tilde{K}_{\Tilde{h}_{m_1}}(x_i,x_j)}\frac{1}{n}\sum_{j=1}^{n}\Tilde{K}_{\Tilde{h}_{m_1}}(x_i,x_j)}\\
&=\frac{1}{\sqrt{n}}\sum_{i=1}^{n}\left(1-\frac{d_i}{\Tilde{p}(x_i;\beta^*)}\right)\frac{\sum_{j=1}^{n}d_j(y_j-m_1(x_i))\Tilde{K}_{\Tilde{h}_{m_1}}(x_i,x_j)}{n\hat{p}(x_i)\hat{f}(x_i)},
\end{split}
\end{equation*}
where $\hat{p}(x)$ is a nonparametric estimation of $PS$ model and $\hat{f}(x)$ is the kernel density estimator of $f(x)$. By standard arguments in nonparametric estimation (see \cite{nonparametric-reference1}, \cite{nonparametric-reference2} and \cite{nonparametric-reference3}), we have $\sup_{x \in \mathcal{X}}|\hat{p}(x)-p(x)|=O_p(\Tilde{h}^s_{m_1}+\sqrt{\frac{\log(n)}{n\Tilde{h}^p_{m_1}}})$ and $\sup_{x \in \mathcal{X}}|\hat{f}(x)-f(x)|=O_p(\Tilde{h}^s_{m_1}+\sqrt{\frac{\log(n)}{n\Tilde{h}^p_{m_1}}})$ under Assumption \ref{nonparaassumption}. Then we can have

\begin{equation*}
\begin{split}
\sqrt{n}C_n &=\frac{1}{\sqrt{n}}\sum_{i=1}^{n}\left(1-\frac{d_i}{\Tilde{p}(x_i;\beta^*)}\right)\frac{\sum_{j=1}^{n}d_j(y_j-m_1(x_i))\Tilde{K}_{\Tilde{h}_{m_1}}(x_i,x_j)}{np(x_i)f(x_i)}+o_p(1)\\
&=\frac{1}{n\sqrt{n}}\sum_{i=1}^{n}\sum_{j=1}^{n}\frac{1-\frac{d_i}{\Tilde{p}(x_i;\beta^*)}}{p(x_i)f(x_i)}d_j[y_j-m_1(x_i)]\Tilde{K}_{\Tilde{h}_{m_1}}(x_i,x_j)+o_p(1).
\end{split}
\end{equation*}

\noindent We further write $\sqrt{n}C_n$ in the form of U-statistics. Let
$H_{ij}=\frac{1-\frac{d_i}{\Tilde{p}(x_i;\beta^*)}}{p(x_i)f(x_i)}d_j[y_j-m_1(x_i)]$, we have
\begin{equation*}
\sqrt{n}C_n=\frac{n-1}{\sqrt{n}}\frac{1}{n(n-1)}\sum_{i=1}^{n}\sum_{j\neq  i}^{n}\left[\frac{H_{ij}+H_{ji}}{2}\right]\Tilde{K}_{\Tilde{h}_{m_1}}(x_i,x_j)+o_p(1)=\frac{n-1}{\sqrt{n}}U_n+o_p(1),
\end{equation*}
where $U_n:=\frac{1}{n(n-1)}\sum_{i=1}^{n}\sum_{j\neq  i}^{n}\left[\frac{H_{ij}+H_{ji}}{2}\right]\Tilde{K}_{\Tilde{h}_{m_1}}(x_i,x_j)$. Next, we compute the conditional expectation of $\left[\frac{H_{ij}+H_{ji}}{2}\right]\Tilde{K}_{\Tilde{h}_{m_1}}(x_i,x_j)$. We first compute
\begin{align*}
&\mathbb{E}[H_{ij}\Tilde{K}_{\Tilde{h}_{m_1}}(x_i,x_j)|x_j,y_j,d_j]\\
&=\mathbb{E}\left\{\mathbb{E}\left[\frac{1-\frac{d_i}{\Tilde{p}(x_i;\beta^*)}}{p(x_i)f(x_i)}d_j[y_j-m_1(x_i)]\Tilde{K}_{\Tilde{h}_{m_1}}(x_i,x_j)\bigg|x_i,x_j,y_j,d_j\right]\bigg|x_j,y_j,d_j\right\}\\
&=\mathbb{E}\left\{\frac{1-\frac{p(x_i)}{\Tilde{p}(x_i;\beta^*)}}{p(x_i)f(x_i)}d_j[y_j-m_1(x_i)]\Tilde{K}_{\Tilde{h}_{m_1}}(x_i,x_j)\bigg|x_j,y_j,d_j\right\}\\
&=\int w(x)d_j[y_j-m_1(x)]\Tilde{K}_{\Tilde{h}_{m_1}}(x,x_j)f(x)dx.
\end{align*}
Let $t=\frac{x-x_j}{\Tilde{h}_{m_1}}$ and further apply Taylor series expansion. Based on Assumption \ref{nonparaassumption}, we have
\begin{align*}
&\mathbb{E}[H_{ij}\Tilde{K}_{\Tilde{h}_{m_1}}(x_i,x_j)|x_j,y_j,d_j]\\
&=\int w(x_j+t\Tilde{h}_{m_1})d_j[y_j-m_1(x_j+t\Tilde{h}_{m_1})]\Tilde{K}(t)f(x_j+t\Tilde{h}_{m_1})dt\\
&=w(x_j)d_j[y_j-m_1(x_j)]f(x_j)+O(\Tilde{h}^s_{m_1})\\
&=\left(\frac{1}{p(x_j)}-\frac{1}{\Tilde{p}(x_j;\beta^*)}\right)d_j[y_j-m_1(x_j)] + O(\Tilde{h}^s_{m_1}).
\end{align*}

\noindent Similarly, we can obtain
\begin{equation*}
\begin{split}
&\mathbb{E}[H_{ji}\Tilde{K}_{\Tilde{h}_{m_1}}(x_i,x_j)|x_j,y_j,d_j]\\
&=\mathbb{E}\left\{\frac{1-\frac{d_j}{\Tilde{p}(x_j;\beta^*)}}{p(x_j)f(x_j)}d_i[y_i-m_1(x_j)]\Tilde{K}_{\Tilde{h}_{m_1}}(x_i,x_j)\bigg|x_j,y_j,d_j\right\}\\
&=\frac{1-\frac{d_j}{\Tilde{p}(x_j;\beta^*)}}{p(x_j)f(x_j)}\mathbb{E}\left\{\mathbb{E}\left[d_i(y_i-m_1(x_j))\Tilde{K}_{\Tilde{h}_{m_1}}(x_i,x_j)|x_i,x_j,y_j,d_j\right]|x_j,y_j,d_j\right\}\\
&=\frac{1-\frac{d_j}{\Tilde{p}(x_j;\beta^*)}}{p(x_j)f(x_j)}\mathbb{E}\left\{p(x_i)[m_1(x_i)-m_1(x_j)]\Tilde{K}_{\Tilde{h}_{m_1}}(x_i,x_j)|x_j,y_j,d_j\right\}\\
&=\frac{1-\frac{d_j}{\Tilde{p}(x_j;\beta^*)}}{p(x_j)f(x_j)}\int p(x)[m_1(x)-m_1(x_j)]\Tilde{K}_{\Tilde{h}_{m_1}}(x,x_j)f(x) dx\\
&=\frac{1-\frac{d_j}{\Tilde{p}(x_j;\beta^*)}}{p(x_j)f(x_j)}\int p(x_j+t\Tilde{h}_{m_1})[m_1(x_j+t\Tilde{h}_{m_1})-m_1(x_j)]\Tilde{K}(t)f(x_j+t\Tilde{h}_{m_1})dt\\
&=O(\Tilde{h}^s_{m_1}).
\end{split}
\end{equation*}

\noindent The conditional expectation of $\left[\frac{H_{ij}+H_{ji}}{2}\right]\Tilde{K}_{\Tilde{h}_{m_1}}(x_i,x_j)$ is $\frac{1}{2}\left(\frac{1}{p(x_j)}-\frac{1}{p(x_j;\beta^*)}\right)d_j[y_j-m_1(x_j)]+O(\Tilde{h}^s_{m_1})$. It follows that $\mathbb{E}\left[\frac{H_{ij}+H_{ji}}{2}\Tilde{K}_{\Tilde{h}_{m_1}}(x_i,x_j)\right]=O(\Tilde{h}^s_{m_1})$. Then we can calculate the projection of $U_n$. Based on Assumption \ref{nonparaassumption}, we have $E\big[||\frac{H_{ij}+H_{ji}}{2}\Tilde{K}_{\Tilde{h}_{m_1}}(x_i,x_j)||^2\big]=o(n)$. Applying Lemma 3.1 of \cite{Ustatistics} under Assumption \ref{nonparaassumption}, we obtain $\sqrt{n}C_n=\frac{1}{\sqrt{n}}\sum_{j=1}^n\Big\{\big(\frac{1}{p(x_j)}-\frac{1}{\Tilde{p}(x_j;\beta^*)}\big)d_j[y_j-m_1(x_j)]\Big\}+o_p(1)$. Following the derivation of $\sqrt{n}C_n$, we have
\begin{align*}
\sqrt{n}D_n&=\frac{1}{\sqrt{n}}\sum_{i=1}^{n}d_i\left(\frac{1}{\Tilde{p}(x_i;\beta^*)}-\frac{1}{\Tilde{p}(x_i;\hat{\beta})}\right)[\hat{m}_1(x_i)-m_1(x_i)]\\
&=\sqrt{n}(\beta^*-\hat{\beta})^T \frac{1}{n}\sum_{i=1}^{n}d_i\frac{\partial\frac{1}{\Tilde{p}(x_i;\beta)}}{\partial\beta}\bigg|_{\beta=\bar{\beta}}\frac{\sum_{j=1}^{n}d_j[y_j-m_1(x_i)]\Tilde{K}_{\Tilde{h}_{m_1}}(x_i,x_j)(x_i,x_j)}{\sum_{j=1}^{n}d_j\Tilde{K}_{\Tilde{h}_{m_1}}(x_i,x_j)(x_i,x_j)}\\
&=\sqrt{n}(\beta^*-\hat{\beta})^T \frac{1}{n}\sum_{i=1}^{n}d_i\frac{\partial\frac{1}{\Tilde{p}(x_i;\beta)}}{\partial\beta}\bigg|_{\beta=\bar{\beta}}\frac{\sum_{j=1}^{n}d_j[y_j-m_1(x_i)]\Tilde{K}_{\Tilde{h}_{m_1}}(x_i,x_j)}{n\frac{\sum_{j=1}^{n}d_jK_{h_{m_1}}(x_i,x_j)}{\sum_{j=1}^{n}\Tilde{K}_{\Tilde{h}_{m_1}}(x_i,x_j)}\frac{1}{n}\sum_{j=1}^{n}\Tilde{K}_{\Tilde{h}_{m_1}}(x_i,x_j)}\\
&=\sqrt{n}(\beta^*-\hat{\beta})^T \frac{1}{n}\sum_{i=1}^{n}d_i\frac{\partial\frac{1}{\Tilde{p}(x_i;\beta)}}{\partial \beta}\bigg|_{\beta=\bar{\beta}}\frac{\sum_{j=1}^{n}d_j[y_j-m_1(x_i)]\Tilde{K}_{\Tilde{h}_{m_1}}(x_i,x_j)}{n\hat{p}(x_i)\hat{f}(x_i)}\\
&=\sqrt{n}(\beta^*-\hat{\beta})^T \frac{1}{n}\sum_{i=1}^{n}d_i\frac{\partial\frac{1}{\Tilde{p}(x_i;\beta)}}{\partial \beta}\bigg|_{\beta=\bar{\beta}}\frac{\sum_{j=1}^{n}d_j[y_j-m_1(x_i)]\Tilde{K}_{\Tilde{h}_{m_1}}(x_i,x_j)}{np(x_i)f(x_i)}+o_p(1)\\
&=\sqrt{n}(\beta^*-\hat{\beta})^T \frac{n-1}{n\sqrt{n}}\left\{\frac{1}{\sqrt{n}}\sum_{j=1}^n\frac{\partial\frac{1}{\Tilde{p}(x_j;\beta)}}{\partial \beta}\bigg|_{\beta=\bar{\beta}}d_j[y_j-m_1(x_j)]+O(\sqrt{n}\Tilde{h}^s_{m_1})+o_p(1)\right\}+o_p(1).
\end{align*}
\noindent Recall that $\bar{\beta}$ is defined in Section 6.2.1. It follows from the Central Limit Theorem and Slutsky's Theorem that $\sqrt{n}D_n=o_p(1)$ under Assumptions \ref{psassumption} and \ref{nonparaassumption}. Consequently, combing the terms in (\ref{6.2.2}), we have
\begin{equation*}
\begin{split}
&\sqrt{n}(\hat{\theta}_1-\theta_1)\\
&=\frac{1}{\sqrt{n}}\sum_{i=1}^{n}\left\{\frac{d_iy_i}{\Tilde{p}(x_i;\beta^*)}+\left(1-\frac{d_i}{\Tilde{p}(x_i;\beta^*)}\right)m_1(x_i)-\theta_1\right\}\\
&\quad+\frac{1}{\sqrt{n}}\sum_{i=1}^n\left\{\left(\frac{1}{p(x_i)}-\frac{1}{\Tilde{p}(x_i;\beta^*)}\right)d_i[y_i-m_1(x_i)]\right\}\\
&\quad+o_p(1)\\
&=\frac{1}{\sqrt{n}}\sum_{i=1}^{n}\left\{\frac{d_iy_i}{p(x_i)}+\left(1-\frac{d_i}{p(x_i)}\right)m_1(x_i)-\theta_1\right\}+o_p(1).
\end{split}
\end{equation*}
Similarly, we can derive the form of $\sqrt{n}(\hat{\theta}_0-\theta_0)$. Therefore,
\begin{equation*}
\begin{split}
&\sqrt{n}(\hat{\Delta}_2-\Delta)\\
&=\sqrt{n}[(\hat{\theta}_1-\hat{\theta}_0)-(\theta_1-\theta_0)]\\
&=\frac{1}{\sqrt{n}}\sum_{i=1}^{n}\left\{\frac{d_iy_i}{p(x_i)}+\left(1-\frac{d_i}{p(x_i)}\right)m_1(x_i)-\frac{(1-d_i)y_i}{1-p(x_i)}-\left(1-\frac{1-d_i}{1-p(x_i)}\right)m_0(x_i)-(\theta_1-\theta_0)\right\}+o_p(1)\\
&=\frac{1}{\sqrt{n}}\sum_{i=1}^{n}\Phi(x_i,y_i,d_i)+o_p(1).
\end{split}
\end{equation*}

\noindent Note that $\mathbb{E}\left\{\Phi(X,Y,D)\right\}=0$. We further assume $\mathbb{E}\left\{\Phi(X,Y,D)^2\right\} < \infty$. It follows from the Central Limit Theorem and Slutsky's Theorem that $\sqrt{n}(\hat{\Delta}_2-\Delta)$ converges in distribution to $N(0,\Sigma_1)$. In other words, the asymptotic variance of $\hat{\Delta}_2$ achieves the semiparametric efficiency bound no matter whether $PS$ model is globally/locally misspecified or not.

\subsubsection{Nonparametric $PS$ model and parametric $OR$ model}

When $PS$ model is nonparametric and $OR$ model is parametric, we have
\begin{eqnarray}\label{6.2.3}
\sqrt{n}(\hat{\theta}_1-\theta_1)&=&\sqrt{n}\left\{n^{-1}\sum_{i=1}^{n}\left[\frac{d_iy_i}{\hat{p}(x_i)}+\left(1-\frac{d_i}{\hat{p}(x_i)}\right)\Tilde{m}_1(x_i;\hat{\gamma}_1)\right]-\theta_1\right\} \nonumber\\
&=&\frac{1}{\sqrt{n}}\sum_{i=1}^{n}\left\{\frac{d_iy_i}{p(x_i)}+\left(1-\frac{d_i}{p(x_i)}\right)\Tilde{m}_1(x_i;\gamma^*_1)-\theta_1\right\}\nonumber\\
&&+\frac{1}{\sqrt{n}}\sum_{i=1}^{n}d_i\left(\frac{1}{\hat{p}(x_i)}-\frac{1}{p(x_i)}\right)[y_i-\Tilde{m}_1(x_i;\gamma^*_1)] \nonumber\\
&&+\frac{1}{\sqrt{n}}\sum_{i=1}^{n}\left(1-\frac{d_i}{\hat{p}(x_i)}\right)[\Tilde{m}_1(x_i;\hat{\gamma}_1)-\Tilde{m}_1(x_i;\gamma^*_1)] \nonumber\\
&:=& \sqrt{n}A_n + \sqrt{n}B_n + \sqrt{n}C_n.
\end{eqnarray}

\newcontent Under Assumption \ref{nonparaassumption}, due to the consistency of $\hat{p}(x)$ to $p(x)$ mentioned in Section 6.2.2, we have $B_n=B^*_n+o_p(n^{-1/2})$ and $C_n=C^*_n+o_p(n^{-1/2})$ with
\begin{equation*}
\begin{split}
&B^*_n=n^{-1}\sum_{i=1}^{n}d_i\frac{1}{p^2(x_i)}[p(x_i)-\hat{p}(x_i)][y_i-\Tilde{m}_1(x_i;\gamma^*_1)],\\
&C^*_n=n^{-1}\sum_{i=1}^{n}\left(1-\frac{d_i}{p(x_i)}\right)[\Tilde{m}_1(x_i;\hat{\gamma}_1)-\Tilde{m}_1(x_i;\gamma^*_1)].
\end{split}
\end{equation*}

\noindent Consequently,
\begin{equation*}
\begin{split}
\sqrt{n}B_n&=\sqrt{n}B^*_n+o_p(1) \\
&=\frac{1}{\sqrt{n}}\sum_{i=1}^{n}d_i\frac{1}{p^2(x_i)}\left[p(x_i)-\frac{\sum_{j=1}^{n}d_j\Tilde{L}_{\Tilde{b}}(x_i,x_j)}{\sum_{j=1}^{n}\Tilde{L}_{\Tilde{b}}(x_i,x_j)}\right][y_i-\Tilde{m}_1(x_i;\gamma^*_1)]+o_p(1)\\
&=\frac{1}{\sqrt{n}}\sum_{i=1}^{n}d_i\frac{1}{p^2(x_i)}\frac{\sum_{j=1}^{n}[p(x_i)-d_j]\Tilde{L}_{\Tilde{b}}(x_i,x_j)}{n\frac{1}{n}\sum_{j=1}^{n}\Tilde{L}_{\Tilde{b}}(x_i,x_j)}[y_i-\Tilde{m}_1(x_i;\gamma^*_1)]+o_p(1)\\
&=\frac{1}{\sqrt{n}}\sum_{i=1}^{n}d_i\frac{1}{p^2(x_i)}\frac{\sum_{j=1}^{n}[p(x_i)-d_j]\Tilde{L}_{\Tilde{b}}(x_i,x_j)}{n\hat{f}(x_i)}[y_i-\Tilde{m}_1(x_i;\gamma^*_1)]+o_p(1),
\end{split}
\end{equation*}
where $\hat{f}(x)$ is the kernel density estimator of $f(x)$. Under Assumption \ref{nonparaassumption}, due to the consistency of $\hat{f}(x)$ to $f(x)$ as previously mentioned in Section 6.2.2, we have
\begin{equation*}
\begin{split}
\sqrt{n}B_n &=\frac{1}{\sqrt{n}}\sum_{i=1}^{n}d_i\frac{1}{p^2(x_i)}\frac{\sum_{j=1}^{n}[p(x_i)-d_j]\Tilde{L}_{\Tilde{b}}(x_i,x_j)}{nf(x_i)}[y_i-\Tilde{m}_1(x_i;\gamma^*_1)]+o_p(1)\\
&=\frac{1}{n\sqrt{n}}\sum_{i=1}^{n}\sum_{j=1}^{n}\frac{d_i[y_i-\Tilde{m}_1(x_i;\gamma^*_1)]}{p^2(x_i)f(x_i)}[p(x_i)-d_j]\Tilde{L}_{\Tilde{b}}(x_i,x_j)+o_p(1).
\end{split}
\end{equation*}

\noindent We further write $\sqrt{n}B_n$ in the form of U-statistics. Let
$H_{ij}=\frac{d_i(y_i-\Tilde{m}_1(x_i;\gamma^*_1))}{p^2(x_i)f(x_i)}[p(x_i)-d_j]$, we have:
\begin{equation*}
\sqrt{n}B_n=\frac{n-1}{\sqrt{n}}\frac{1}{n(n-1)}\sum_{i=1}^{n}\sum_{j\neq i}^{n}\left[\frac{H_{ij}+H_{ji}}{2}\right]\Tilde{L}_{\Tilde{b}}(x_i,x_j)+o_p(1)=\frac{n-1}{\sqrt{n}}U_n +o_p(1),
\end{equation*}
where $U_n:=\frac{1}{n(n-1)}\sum_{i=1}^{n}\sum_{j\neq i}^{n}\left[\frac{H_{ij}+H_{ji}}{2}\right]\Tilde{L}_{\Tilde{b}}(x_i,x_j)$. Next, we compute the conditional expectation of $\left[\frac{H_{ij}+H_{ji}}{2}\right]\Tilde{L}_{\Tilde{b}}(x_i,x_j)$. We first compute
\begin{equation*}
\begin{split}
&\mathbb{E}[H_{ij}\Tilde{L}_{\Tilde{b}}(x_i,x_j)|x_j,y_j,d_j]\\
&=\mathbb{E}\left\{\frac{d_i(y_i-\Tilde{m}_1(x_i;\gamma^*_1))}{p^2(x_i)f(x_i)}[p(x_i)-d_j]\Tilde{L}_{\Tilde{b}}(x_i,x_j)\bigg|x_j,y_j,d_j\right\}\\
&=\mathbb{E}\left\{\mathbb{E}\left[\frac{d_i(y_i-\Tilde{m}_1(x_i;\gamma^*_1))}{p^2(x_i)f(x_i)}[p(x_i)-d_j]\Tilde{L}_{\Tilde{b}}(x_i,x_j)\bigg|x_i,x_j,y_j,d_j\right]\bigg|x_j,y_j,d_j\right\}\\
&=\mathbb{E}\left\{\frac{m_1(x_i)-\Tilde{m}_1(x_i;\gamma^*_1)}{p(x_i)f(x_i)}[p(x_i)-d_j]\Tilde{L}_{\Tilde{b}}(x_i,x_j)\bigg|x_j,y_j,d_j\right\}\\
&=\int w(x)[p(x)-d_j]\Tilde{L}_{\Tilde{b}}(x,x_j)f(x)dx.
\end{split}
\end{equation*}
Let $t=\frac{x-x_j}{\Tilde{b}}$ and further apply Taylor series expansion. Based on Assumption \ref{nonparaassumption}, we have
\begin{equation*}
\begin{split}
&\mathbb{E}[H_{ij}\Tilde{L}_{\Tilde{b}}(x_i,x_j)|x_j,y_j,d_j]\\
&=\int w(x_j+t\Tilde{b})[p(x_j+t\Tilde{b})-d_j]\Tilde{L}(t)f(x_j+t\Tilde{b})dt\\
&=w(x_j)[p(x_j)-d_j]f(x_j)+O(\Tilde{b}^s)\\
&=\frac{m_1(x_j)-\Tilde{m}_1(x_j;\gamma^*_1)}{p(x_j)f(x_j)}[p(x_j)-d_j]f(x_j)+O(\Tilde{b}^s)\\
&=\left(1-\frac{d_j}{p(x_j)}\right)[m_1(x_j)-\Tilde{m}_1(x_j;\gamma^*_1)]+O(\Tilde{b}^s).
\end{split}
\end{equation*}

\noindent Similarly, we can obtain
\begin{align*}
&\mathbb{E}[H_{ji}\Tilde{L}_{\Tilde{b}}(x_i,x_j)|x_j,y_j,d_j]\\
&=\mathbb{E}\left\{\frac{d_j[y_j-\Tilde{m}_1(x_j;\gamma^*_1)]}{p^2(x_j)f(x_j)}[p(x_j)-d_i]\Tilde{L}_{\Tilde{b}}(x_i,x_j)\bigg|x_j,y_j,d_j\right\}\\
&=\frac{d_j[y_j-\Tilde{m}_1(x_j;\gamma^*_1)]}{p^2(x_j)f(x_j)}\mathbb{E}\left\{\mathbb{E}\left[(p(x_j)-d_i)\Tilde{L}_{\Tilde{b}}(x_i,x_j)|x_i,x_j,y_j,d_j\right]|x_j,y_j,d_j\right\}\\
&=\frac{d_j[y_j-\Tilde{m}_1(x_j;\gamma^*_1)]}{p^2(x_j)f(x_j)}\mathbb{E}\left\{[p(x_j)-p(x_i)]\Tilde{L}_{\Tilde{b}}(x_i,x_j)|x_j,y_j,d_j\right\}\\
&=\frac{d_j(y_j-\Tilde{m}_1(x_j;\gamma^*_1))}{p^2(x_j)f(x_j)}\int [p(x_j)-p(x)]\Tilde{L}_{\Tilde{b}}(x,x_j)f(x)dx\\
&=\frac{d_j(y_j-\Tilde{m}_1(x_j;\gamma^*_1))}{p^2(x_j)f(x_j)}\int [p(x_j)-p(x_j+t\Tilde{b})]\Tilde{L}(t)f(x_j+t\Tilde{b})dt\\
&=O(\Tilde{b}^s).
\end{align*}

\noindent The conditional expectation of $\left[\frac{H_{ij}+H_{ji}}{2}\right]\Tilde{L}_{\Tilde{b}}(x_i,x_j)$ is  $\frac{1}{2}\left(1-\frac{d_j}{p(x_j)}\right)[m_1(x_j)-\Tilde{m}_1(x_j;\gamma^*_1)]+O(\Tilde{b}^s)$. It follows that $\mathbb{E}\left[\frac{H_{ij}+H_{ji}}{2}\Tilde{L}_{\Tilde{b}}(x_i,x_j)\right]=O(\Tilde{b}^s)$.
Then we can calculate the projection of $U_n$. Based on Assumption \ref{nonparaassumption}, we have $E\big[||\frac{H_{ij}+H_{ji}}{2}\Tilde{L}_{\Tilde{b}}(x_i,x_j)||^2\big]=o(n)$. Applying Lemma 3.1 of \cite{Ustatistics} under Assumption \ref{nonparaassumption}, we obtain $\sqrt{n}B_n=\frac{1}{\sqrt{n}}\sum_{j=1}^n\left(1-\frac{d_j}{p(x_j)}\right)[m_1(x_j)-\Tilde{m}_1(x_j;\gamma^*_1)]+o_p(1)$. For $\sqrt{n}C_n$, similar to the derivation of $\sqrt{n}D_n$ in Section 6.2.1, we have $\sqrt{n}C_n=o_p(1)$. Consequently, combining the terms in (\ref{6.2.3}), we have
\begin{equation*}
\begin{split}
&\sqrt{n}(\hat{\theta}_1-\theta_1)\\
&=\frac{1}{\sqrt{n}}\sum_{i=1}^{n}\left\{\frac{d_iy_i}{p(x_i)}+\left(1-\frac{d_i}{p(x_i)}\right)\Tilde{m}_1(x_i;\gamma^*_1)-\theta_1\right\}\\
&\quad+\frac{1}{\sqrt{n}}\sum_{i=1}^n\left(1-\frac{d_i}{p(x_i)}\right)[m_1(x_i)-\Tilde{m}_1(x_i;\gamma^*_1)]\\
&\quad+o_p(1)\\
&=\frac{1}{\sqrt{n}}\sum_{i=1}^{n}\left\{\frac{d_iy_i}{p(x_i)}+\left(1-\frac{d_i}{p(x_i)}\right)m_1(x_i)-\theta_1\right\}+o_p(1).
\end{split}
\end{equation*}
Similarly, we can derive the form of $\sqrt{n}(\hat{\theta}_0-\theta_0)$. As a result,
\begin{equation*}
\begin{split}
&\sqrt{n}(\hat{\Delta}_3-\Delta)\\
&=\sqrt{n}[(\hat{\theta}_1-\hat{\theta}_0)-(\theta_1-\theta_0)]\\
&=\frac{1}{\sqrt{n}}\sum_{i=1}^{n}\left\{\frac{d_iy_i}{p(x_i)}+\left(1-\frac{d_i}{p(x_i)}\right)m_1(x_i)-\frac{(1-d_i)y_i}{1-p(x_i)}-\left(1-\frac{1-d_i}{1-p(x_i)}\right)m_0(x_i)-(\theta_1-\theta_0)\right\}+o_p(1)\\
&=\frac{1}{\sqrt{n}}\sum_{i=1}^{n}\Phi(x_i,y_i,d_i)+o_p(1).
\end{split}
\end{equation*}

\noindent Note that $\mathbb{E}\left\{\Phi(X,Y,D)\right\}=0$. We further assume $\mathbb{E}\left\{\Phi(X,Y,D)^2\right\} < \infty$. It follows from the Central Limit Theorem and Slutsky's Theorem that $\sqrt{n}(\hat{\Delta}_3-\Delta)$ converges in distribution to $N(0,\Sigma_1)$. In other words, the asymptotic variance of $\hat{\Delta}_3$ achieves the semiparametric efficiency bound no matter whether $OR$ model is globally/locally misspecified or not.

\subsubsection{Nonparametric $PS$ model and $OR$ model}

When $PS$ model and $OR $ model are both nonparametric, we have
\begin{eqnarray}\label{6.2.4}
\sqrt{n}(\hat{\theta}_1-\theta_1)
&=&\sqrt{n}\left\{n^{-1}\sum_{i=1}^{n}\left[\frac{d_iy_i}{\hat{p}(x_i)}+\left(1-\frac{d_i}{\hat{p}(x_i)}\right)\hat{m}_1(x_i)\right]-\theta_1\right\}\nonumber\\
&=&\frac{1}{\sqrt{n}}\sum_{i=1}^{n}\left\{\frac{d_iy_i}{p(x_i)}+\left(1-\frac{d_i}{p(x_i)}\right)m_1(x_i)-\theta_1\right\} \nonumber\\
&&+\frac{1}{\sqrt{n}}\sum_{i=1}^{n}d_i\left(\frac{1}{\hat{p}(x_i)}-\frac{1}{p(x_i)}\right)[y_i-m_1(x_i)] \nonumber\\
&&+\frac{1}{\sqrt{n}}\sum_{i=1}^{n}\left(1-\frac{d_i}{\hat{p}(x_i)}\right)[\hat{m}_1(x_i)-m_1(x_i)] \nonumber\\
&:=& \sqrt{n}A_n + \sqrt{n}B_n + \sqrt{n}C_n.
\end{eqnarray}

\newcontent Under Assumption \ref{nonparaassumption}, due to the consistency of $\hat{p}(x_i)$ to $p(x_i)$ mentioned in Section 6.2.2, we have $B_n=B^*_n+o_p(n^{-1/2})$, $C_n=C^*_n+o_p(n^{-1/2})$ with
\begin{equation*}
\begin{split}
&B^*_n=n^{-1}\sum_{i=1}^{n}d_i\frac{1}{p^2(x_i)}[p(x_i)-\hat{p}(x_i)][y_i-m_1(x_i),\\
&C^*_n=n^{-1}\sum_{i=1}^{n}\left(1-\frac{d_i}{p(x_i)}\right)[\hat{m}_1(x_i)-m_1(x_i)].
\end{split}
\end{equation*}

\noindent Consequently, following the derivations of $\sqrt{n}B_n$ in Section 6.2.3, we have $\sqrt{n}B_n=o_p(1)$. Similar to the derivations of $\sqrt{n}C_n$ in Section 6.2.2, we can get $\sqrt{n}C_n=o_p(1)$. Combing the terms in (\ref{6.2.4}), we have:
\begin{equation*}
\begin{split}
\sqrt{n}(\hat{\theta}_1-\theta_1)=\frac{1}{\sqrt{n}}\sum_{i=1}^{n}\left\{\frac{d_iy_i}{p(x_i)}+\left(1-\frac{d_i}{p(x_i)}\right)m_1(x_i)-\theta_1\right\}+o_p(1).
\end{split}
\end{equation*}
Similarly, we can derive the form of $\sqrt{n}(\hat{\theta}_0-\theta_0)$. As a result,
\begin{equation*}
\begin{split}
&\sqrt{n}(\hat{\Delta}_4-\Delta) \\
&=\sqrt{n}[(\hat{\theta}_1-\hat{\theta}_0)-(\theta_1-\theta_0)]\\
&=\frac{1}{\sqrt{n}}\sum_{i=1}^{n}\left\{\frac{d_iy_i}{p(x_i)}+\left(1-\frac{d_i}{p(x_i)}\right)m_1(x_i)-\frac{(1-d_i)y_i}{1-p(x_i)}-\left(1-\frac{1-d_i}{1-p(x_i)}\right)m_0(x_i)-(\theta_1-\theta_0)\right\}+o_p(1)\\
&=\frac{1}{\sqrt{n}}\sum_{i=1}^{n}\Phi(x_i,y_i,d_i)+o_p(1).
\end{split}
\end{equation*}
\noindent Note that $\mathbb{E}\left\{\Phi(X,Y,D)\right\}=0$. We further assume $\mathbb{E}\left\{\Phi(X,Y,D)^2\right\} < \infty$. It follows from the Central Limit Theorem and Slutsky's Theorem that $\sqrt{n}(\hat{\Delta}_4-\Delta)$ converges in distribution to $N(0,\Sigma_1)$. The asymptotic variance of $\hat{\Delta}_4$ achieves the semiparametric efficiency bound $\Sigma_1$.

\subsubsection{Semiparametric $PS$ model and parametric $OR$ model}

When $PS$ model is semiparametric and $OR$ model is parametric, we have
\begin{eqnarray} \label{6.2.5}
\sqrt{n}(\hat{\theta}_1-\theta_1)
&=&\sqrt{n}\left\{n^{-1}\sum_{i=1}^{n}\left[\frac{d_iy_i}{\hat{g}(\alpha^Tx_i)}+\left(1-\frac{d_i}{\hat{g}(\alpha^Tx_i)}\right)\Tilde{m}_1(x_i;\hat{\gamma}_1)\right]-\theta_1\right\}\nonumber\\
&=&\frac{1}{\sqrt{n}}\sum_{i=1}^{n}\left\{\frac{d_iy_i}{g(\alpha^Tx_i)}+\left(1-\frac{d_i}{g(\alpha^Tx_i)}\right)\Tilde{m}_1(x_i;\gamma^*_1)-\theta_1\right\} \nonumber\\
&&+\frac{1}{\sqrt{n}}\sum_{i=1}^{n}d_i\left(\frac{1}{\hat{g}(\alpha^Tx_i)}-\frac{1}{g(\alpha^Tx_i)}\right)[y_i-\Tilde{m}_1(x_i;\gamma^*_1)]\nonumber \\
&&+\frac{1}{\sqrt{n}}\sum_{i=1}^{n}\left(1-\frac{d_i}{g(\alpha^Tx_i)}\right)[\Tilde{m}_1(x_i;\hat{\gamma}_1)-\Tilde{m}_1(x_i;\gamma^*_1)]\nonumber \\
&&+\frac{1}{\sqrt{n}}\sum_{i=1}^{n}d_i\left(\frac{1}{g(\alpha^Tx_i)}-\frac{1}{\hat{g}(\alpha^Tx_i)}\right)[\Tilde{m}_1(x_i;\hat{\gamma}_1)-\Tilde{m}_1(x_i;\gamma^*_1)]\nonumber \\
&:=& \sqrt{n}A_n + \sqrt{n}B_n + \sqrt{n}C_n + \sqrt{n}D_n.
\end{eqnarray}

\noindent Under Assumption \ref{semiassumption}, we have $\sup_{x \in \mathcal{X}}|\hat{g}(\alpha^Tx)-g(\alpha^Tx)|=O_p(b^2+\sqrt{\frac{\log(n)}{nb}})$ (see \cite{nonparametric-reference1}, \cite{nonparametric-reference2} and \cite{nonparametric-reference3}). Due to the consistency of $\hat{g}(\alpha^Tx_i)$ to $g(\alpha^Tx_i)$, we have $B_n=B^*_n+o_p(n^{-1/2})$ and $D_n=D^*_n+o_p(n^{-1/2})$ with
\begin{equation*}
\begin{split}
&B^*_n=n^{-1}\sum_{i=1}^{n}d_i\frac{1}{g^2(\alpha^Tx_i)}[g(\alpha^Tx_i)-\hat{g}(\alpha^Tx_i)][y_i-\Tilde{m}_1(x_i;\gamma^*_1)],\\
&D^*_n=n^{-1}\sum_{i=1}^{n}d_i\frac{1}{g^2(\alpha^Tx_i)}[\hat{g}(\alpha^Tx_i)-g(\alpha^Tx_i)][\Tilde{m}_1(x_i;\hat{\gamma}_1)-\Tilde{m}_1(x_i;\gamma^*_1)].
\end{split}
\end{equation*}

\noindent Consequently,
\begin{equation*}
\begin{split}
\sqrt{n}B_n&=\sqrt{n}B^*_n+o_p(1) \\
&=\frac{1}{\sqrt{n}}\sum_{i=1}^{n}d_i\frac{1}{g^2(\alpha^Tx_i)}\left[g(\alpha^Tx_i)-\frac{\sum_{j=1}^{n}d_jL_b(\alpha^Tx_i,\alpha^Tx_j)}{\sum_{j=1}^{n}L_b(\alpha^Tx_i,\alpha^Tx_j)}\right][y_i-\Tilde{m}_1(x_i;\gamma^*_1)]+o_p(1)\\
&=\frac{1}{\sqrt{n}}\sum_{i=1}^{n}d_i\frac{1}{g^2(\alpha^Tx_i)}\frac{\sum_{j=1}^{n}[g(\alpha^Tx_i)-d_j]L_b(\alpha^Tx_i,\alpha^Tx_j)}{n\frac{1}{n}\sum_{j=1}^{n}L_b(\alpha^Tx_i,\alpha^Tx_j)}[y_i-\Tilde{m}_1(x_i;\gamma^*_1)]+o_p(1)\\
&=\frac{1}{\sqrt{n}}\sum_{i=1}^{n}d_i\frac{1}{g^2(\alpha^Tx_i)}\frac{\sum_{j=1}^{n}(g(\alpha^Tx_i)-d_j)L_b(\alpha^Tx_i,\alpha^Tx_j)}{n\hat{\Tilde{f}}(\alpha^Tx_i)}[y_i-\Tilde{m}_1(x_i;\gamma^*_1)]+o_p(1),
\end{split}
\end{equation*}
where $\hat{\Tilde{f}}(\alpha^Tx)$ is the kernel density estimator of $\Tilde{f}(\alpha^Tx)$. Under Assumption \ref{semiassumption}, we have $\sup_{x \in \mathcal{X}}|\hat{\Tilde{f}}(\alpha^Tx)-\Tilde{f}(\alpha^Tx)|=O_p(b^2+\sqrt{\frac{\log(n)}{nb}})$ (see \cite{nonparametric-reference1}, \cite{nonparametric-reference2} and \cite{nonparametric-reference3}). Due to the consistency of $\hat{\Tilde{f}}(\alpha^Tx)$ to $\Tilde{f}(\alpha^Tx)$, we have
\begin{equation*}
\begin{split}
\sqrt{n}B_n&=\frac{1}{\sqrt{n}}\sum_{i=1}^{n}d_i\frac{1}{g^2(\alpha^Tx_i)}\frac{\sum_{j=1}^{n}[g(\alpha^Tx_i)-d_j]L_b(\alpha^Tx_i,\alpha^Tx_j)}{n\Tilde{f}(\alpha^Tx_i)}[y_i-\Tilde{m}_1(x_i;\gamma^*_1)]+o_p(1)\\
&=\frac{1}{n\sqrt{n}}\sum_{i=1}^{n}\sum_{j=1}^{n}\frac{d_i[y_i-\Tilde{m}_1(x_i;\gamma^*_1)]}{g^2(\alpha^Tx_i)\Tilde{f}(\alpha^Tx_i)}[g(\alpha^Tx_i)-d_j]L_b(\alpha^Tx_i,\alpha^Tx_j)+o_p(1).
\end{split}
\end{equation*}

\noindent We further write $\sqrt{n}B_n$ in the form of U-statistics. Let
$H_{ij}=\frac{d_i(y_i-\Tilde{m}_1(x_i;\gamma^*_1))}{g^2(\alpha^Tx_i)\Tilde{f}(\alpha^Tx_i)}[g(\alpha^Tx_i)-d_j]$, we have
\begin{equation*}
\sqrt{n}B_n=\frac{n-1}{\sqrt{n}}\frac{1}{n(n-1)}\sum_{i=1}^{n}\sum_{j\neq i}^{n}\left[\frac{H_{ij}+H_{ji}}{2}\right]L_b(\alpha^Tx_i,\alpha^Tx_j)+o_p(1)=\frac{n-1}{\sqrt{n}}U_n +o_p(1),
\end{equation*}
where $U_n:=\frac{1}{n(n-1)}\sum_{i=1}^{n}\sum_{j\neq i}^{n}\left[\frac{H_{ij}+H_{ji}}{2}\right]L_b(\alpha^Tx_i,\alpha^Tx_j)$. Next, we compute the conditional expectation of $\left[\frac{H_{ij}+H_{ji}}{2}\right]L_b(\alpha^Tx_i,\alpha^Tx_j)$. We first compute
\begin{equation*}
\begin{split}
&\mathbb{E}[H_{ij}L_b(\alpha^Tx_i,\alpha^Tx_j)|x_j,y_j,d_j]\\
&=\mathbb{E}\left\{\frac{d_i(y_i-\Tilde{m}_1(x_i;\gamma^*_1))}{g^2(\alpha^Tx_i)\Tilde{f}(\alpha^Tx_i)}[g(\alpha^Tx_i)-d_j]L_b(\alpha^Tx_i,\alpha^Tx_j)\bigg|x_j,y_j,d_j\right\}\\
&=\mathbb{E}\left\{\mathbb{E}\left[\frac{d_i(y_i-\Tilde{m}_1(x_i;\gamma^*_1))}{g^2(\alpha^Tx_i)\Tilde{f}(\alpha^Tx_i)}[g(\alpha^Tx_i)-d_j]L_b(\alpha^Tx_i,\alpha^Tx_j)\bigg|\alpha^Tx_i,x_j,y_j,d_j\right]\bigg|x_j,y_j,d_j\right\}\\
&=\mathbb{E}\left\{\frac{\mathbb{E}[D(Y-\Tilde{m}_1(X;\gamma^*_1))|\alpha^Tx_i]}{g^2(\alpha^Tx_i)\Tilde{f}(\alpha^Tx_i)}[g(\alpha^Tx_i)-d_j]L_b(\alpha^Tx_i,\alpha^Tx_j)\bigg|x_j,y_j,d_j\right\}\\
&=\int w(\alpha^Tx)[g(\alpha^Tx)-d_j]L_b(\alpha^Tx,\alpha^Tx_j)\Tilde{f}(\alpha^Tx)d(\alpha^Tx).
\end{split}
\end{equation*}
Let $t=\frac{\alpha^Tx-\alpha^Tx_j}{b}$ and further apply Taylor series expansion. Based on Assumption \ref{semiassumption}, we have
\begin{equation*}
\begin{split}
&\mathbb{E}[H_{ij}L_b(\alpha^Tx_i,\alpha^Tx_j)|x_j,y_j,d_j]\\
&=\int w(\alpha^Tx_j+bt)[g(\alpha^Tx_j+bt)-d_j]L(t)\Tilde{f}(\alpha^Tx_j+bt)dt\\
&=\frac{\mathbb{E}[D(Y-\Tilde{m}_1(X;\gamma^*_1))|\alpha^Tx_j]}{g^2(\alpha^Tx_j)\Tilde{f}(\alpha^Tx_j)}[g(\alpha^Tx_j)-d_j]\Tilde{f}(\alpha^Tx_j)+O(b^2)\\
&=\frac{\mathbb{E}\left\{\mathbb{E}[D(Y-\Tilde{m}_1(X;\gamma^*_1))|X]|\alpha^Tx_j\right\}}{g^2(\alpha^Tx_j)}[g(\alpha^Tx_j)-d_j]+O(b^2)\\
&=\frac{\mathbb{E}\left\{p(X)[m_1(X)-\Tilde{m}_1(X;\gamma^*_1)]|\alpha^Tx_j\right\}}{g^2(\alpha^Tx_j)}[g(\alpha^Tx_j)-d_j]+O(b^2).
\end{split}
\end{equation*}

\noindent Similarly, we can obtain
\begin{equation*}
\begin{split}
&\mathbb{E}[H_{ji}L_b(\alpha^Tx_i,\alpha^Tx_j)|x_j,y_j,d_j]\\
&=\mathbb{E}\left\{\frac{d_j(y_j-\Tilde{m}_1(x_j;\gamma^*_1))}{g^2(\alpha^Tx_j)\Tilde{f}(\alpha^Tx_j)}[g(\alpha^Tx_j)-d_i]L_b(\alpha^Tx_i,\alpha^Tx_j)\bigg|x_j,y_j,d_j\right\}\\
&=\frac{d_j(y_j-\Tilde{m}_1(x_j;\gamma^*_1))}{g^2(\alpha^Tx_j)\Tilde{f}(\alpha^Tx_j)}\mathbb{E}\left\{\mathbb{E}\left[(g(\alpha^Tx_j)-d_i)L_b(\alpha^Tx_i,\alpha^Tx_j)|\alpha^Tx_i,x_j,y_j,d_j\right]|x_j,y_j,d_j\right\}\\
&=\frac{d_j(y_j-\Tilde{m}_1(x_j;\gamma^*_1))}{g^2(\alpha^Tx_j)\Tilde{f}(\alpha^Tx_j)}\mathbb{E}\left\{\left[g(\alpha^Tx_j)-g(\alpha^Tx_i)\right]L_b(\alpha^Tx_i,\alpha^Tx_j)|x_j,y_j,d_j\right\}\\
&=\frac{d_j(y_j-\Tilde{m}_1(x_j;\gamma^*_1))}{g^2(\alpha^Tx_j)\Tilde{f}(\alpha^Tx_j)}\int\left[g(\alpha^Tx_j)-g(\alpha^Tx)\right]L_b(\alpha^Tx,\alpha^Tx_j)\Tilde{f}(\alpha^Tx)d(\alpha^Tx)\\
&=\frac{d_j(y_j-\Tilde{m}_1(x_j;\gamma^*_1))}{g^2(\alpha^Tx_j)\Tilde{f}(\alpha^Tx_j)}\int [g(\alpha^Tx_j)-g(\alpha^Tx_j+bt)]L(t)\Tilde{f}(\alpha^Tx_j+bt)dt\\
&=O(b^2).
\end{split}
\end{equation*}

\noindent The conditional expectation of $\left[\frac{H_{ij}+H_{ji}}{2}\right]L_b(\alpha^Tx_i,\alpha^Tx_j)$ is $\frac{\mathbb{E}\{p(X)[m_1(X)-\Tilde{m}_1(X;\gamma^*_1)]|\alpha^Tx_j\}}{2g^2(\alpha^Tx_j)}[g(\alpha^Tx_j)-d_j]+O(b^2)$. It follows that $\mathbb{E}\left[\frac{H_{ij}+H_{ji}}{2}L_b(\alpha^Tx_i,\alpha^Tx_j)\right]=O(b^2)$. Then we can calculate the projection of $U_n$. Based on Assumption \ref{semiassumption}, we have $E\big[||\frac{H_{ij}+H_{ji}}{2}L_b(x_i,x_j)||^2\big]=o(n)$. Applying Lemma 3.1 of \cite{Ustatistics} under Assumption \ref{semiassumption}, we have $\sqrt{n}B_n=\frac{1}{\sqrt{n}}\sum_{j=1}^n\frac{\mathbb{E}\left\{p(X)[m_1(X)-\Tilde{m}_1(X;\gamma^*_1)]|\alpha^Tx_j\right\}}{g^2(\alpha^Tx_j)}[g(\alpha^Tx_j)-d_j]+o_p(1)$. Similar to the derivations of $\sqrt{n}D_n$ in Section 6.2.1, we can get
\begin{equation*}
\sqrt{n}C_n =\sqrt{n}(\hat{\gamma}_1-\gamma^*_1)^T \mathbb{E}\left[(1-\frac{p(X)}{g(\alpha^TX)})\frac{\partial \Tilde{m}_1(X;\gamma_1)}{\partial \gamma_1}\bigg|_{\gamma_1=\bar{\gamma}_1}\right]+o_p(1).
\end{equation*}

\noindent Recall that $\bar{\gamma}_1$ is defined in Section 6.2.1. For $\sqrt{n}D_n$, following the derivations of $\sqrt{n}B_n$, we have
\begin{equation*}
\begin{split}
\sqrt{n}D_n&=\sqrt{n}D^*_n+o_p(1) \\
&=\frac{1}{\sqrt{n}}\sum_{i=1}^{n}d_i\frac{1}{g^2(\alpha^Tx_i)}[\hat{g}(\alpha^Tx_i)-g(\alpha^Tx_i)][\Tilde{m}_1(x_i;\hat{\gamma}_1)-\Tilde{m}_1(x_i;\gamma^*_1)]+o_p(1)\\
&=\frac{1}{\sqrt{n}}\sum_{i=1}^{n}\frac{d_i(\hat{\gamma}_1-\gamma^*_1)^T\frac{\partial \Tilde{m}_1(x_i;\gamma_1)}{\partial\gamma_1}\Big|_{\gamma_1=\bar{\gamma}_1}}{g^2(\alpha^Tx_i)}[\hat{g}(\alpha^Tx_i)-g(\alpha^Tx_i)]+o_p(1)\\
&=\sqrt{n}(\hat{\gamma}_1-\gamma^*_1)^T\frac{1}{n}\sum_{i=1}^{n}\frac{d_i\frac{\partial\Tilde{m}_1(x_i;\gamma_1)}{\partial\gamma_1}\Big|_{\gamma_1=\bar{\gamma}_1}}{g^2(\alpha^Tx_i)}\left[\frac{\sum_{j=1}^{n}d_jL_b(\alpha^Tx_i,\alpha^Tx_j)}{\sum_{j=1}^{n}L_b(\alpha^Tx_i,\alpha^Tx_j)}-g(\alpha^Tx_i)\right]+o_p(1)\\
&=\sqrt{n}(\hat{\gamma}_1-\gamma^*_1)^T\frac{1}{n}\sum_{i=1}^{n}\frac{d_i\frac{\partial\Tilde{m}_1(x_i;\gamma_1)}{\partial\gamma_1}\Big|_{\gamma_1=\bar{\gamma}_1}}{g^2(\alpha^Tx_i)}\frac{\sum_{j=1}^{n}[d_j-g(\alpha^Tx_i)]L_b(\alpha^Tx_i,\alpha^Tx_j)}{n\frac{1}{n}\sum_{j=1}^{n}L_b(\alpha^Tx_i,\alpha^Tx_j)}+o_p(1)\\
&=\sqrt{n}(\hat{\gamma}_1-\gamma^*_1)^T\frac{1}{n^2}\sum_{i=1}^{n}\sum_{j=1}^{n}\frac{d_i\frac{\partial\Tilde{m}_1(x_i;\gamma_1)}{\partial\gamma_1}\Big|_{\gamma_1=\bar{\gamma}_1}}{g^2(A^Tx_i)\hat{\Tilde{f}}(\alpha^Tx_i)}(d_j-g(\alpha^Tx_i))L_b(\alpha^Tx_i,\alpha^Tx_j)+o_p(1)\\
&=\sqrt{n}(\hat{\gamma}_1-\gamma^*_1)^T\frac{1}{n^2}\sum_{i=1}^{n}\sum_{j=1}^{n}\frac{d_i\frac{\partial\Tilde{m}_1(x_i;\gamma_1)}{\partial\gamma_1}\Big|_{\gamma_1=\bar{\gamma}_1}}{g^2(\alpha^Tx_i)\Tilde{f}(\alpha^Tx_i)}(d_j-g(\alpha^Tx_i))L_b(\alpha^Tx_i,\alpha^Tx_j)+o_p(1)\\
&=\sqrt{n}(\hat{\gamma}_1-\gamma^*_1)^T\frac{n-1}{n\sqrt{n}}\left\{\frac{1}{\sqrt{n}}\sum_{j=1}^n\frac{\mathbb{E}\left[p(X)\frac{\partial\Tilde{m}_1(X;\gamma_1)}{\partial \gamma_1}\Big|_{\gamma_1=\bar{\gamma}_1}\bigg|\alpha^Tx_j\right]}{g^2(\alpha^Tx_j)}[d_j-g(\alpha^Tx_j)]+O(\sqrt{n}b^2)+o_p(1)\right\}+o_p(1).
\end{split}
\end{equation*}
\noindent Under Assumptions \ref{orassumption} and \ref{semiassumption}, it follows from the Central Limit Theorem and Slutsky's Theorem that $\sqrt{n}D_n=o_p(1)$. Now we can consider different cases as follows.

\newcontent (a) \textbf{Correctly specified $PS$ model and $OR$ model}

In this case, we have $p(x)=g(\alpha^Tx)$, $m_1(x)=\Tilde{m}_1(x;\gamma_{1,0})=\Tilde{m}_1(x;\gamma^*_1)$ and $m_0(x)=\Tilde{m}_0(x;\gamma_{0,0})=\Tilde{m}_0(x;\gamma^*_0)$. So $\sqrt{n}B_n=o_p(1)$ and $\sqrt{n}C_n=o_p(1)$. Combining the terms in (\ref{6.2.5}), we have
\begin{equation*}
\begin{split}
\sqrt{n}(\hat{\theta}_1-\theta_1)=\frac{1}{\sqrt{n}}\sum_{i=1}^{n}\left\{\frac{d_iy_i}{p(x_i)}+\left(1-\frac{d_i}{p(x_i)}\right)m_1(x_i)-\theta_1\right\}+o_p(1).
\end{split}
\end{equation*}
Similarly, we can derive the form of $\sqrt{n}(\hat{\theta}_0-\theta_0)$. As a result,
\begin{equation*}
\begin{split}
&\sqrt{n}(\hat{\Delta}_5-\Delta) \\
&=\sqrt{n}[(\hat{\theta}_1-\hat{\theta}_0)-(\theta_1-\theta_0)]\\
&=\frac{1}{\sqrt{n}}\sum_{i=1}^{n}\left\{\frac{d_iy_i}{p(x_i)}+\left(1-\frac{d_i}{p(x_i)}\right)m_1(x_i)-\frac{(1-d_i)y_i}{1-p(x_i)}-\left(1-\frac{1-d_i}{1-p(x_i)}\right)m_0(x_i)-(\theta_1-\theta_0)\right\}+o_p(1)\\
&=\frac{1}{\sqrt{n}}\sum_{i=1}^{n}\Phi(x_i,y_i,d_i)+o_p(1).
\end{split}
\end{equation*}

\noindent Note that $\mathbb{E}\left\{\Phi(X,Y,D)\right\}=0$. We further assume $\mathbb{E}\left\{\Phi(X,Y,D)^2\right\} < \infty$. It follows from the Central Limit Theorem and Slutsky's Theorem that $\sqrt{n}(\hat{\Delta}_5-\Delta)$ converges in distribution to $N(0,\Sigma_1)$. The asymptotic variance of $\hat{\Delta}_5$ achieves the semiparametric efficiency bound $\Sigma_1$ when both $PS$ and $OR$ models are correctly specified.

\newcontent (b) \textbf{Correctly specified $PS$ model and misspecified $OR$ model}

In this case, we have $p(x)=g(\alpha^Tx)$, $m_1(x)\neq \Tilde{m}_1(x;\gamma_{1,0})\neq\Tilde{m}_1(x;\gamma^*_1)$ and  $m_0(x)\neq\Tilde{m}_0(x;\gamma_{0,0})\neq\Tilde{m}_0(x;\gamma^*_0)$. So $\sqrt{n}C_n=o_p(1)$. Combining the terms in (\ref{6.2.5}), we have
\begin{equation*}
\begin{split}
&\sqrt{n}(\hat{\theta}_1-\theta_1)\\
&=\frac{1}{\sqrt{n}}\sum_{i=1}^{n}\left\{\frac{d_iy_i}{p(x_i)}+\left(1-\frac{d_i}{p(x_i)}\right)\Tilde{m}_1(x_i;\gamma^*_1)-\theta_1\right\} \\
&\quad+\frac{1}{\sqrt{n}}\sum_{i=1}^n\frac{\mathbb{E}\left\{p(X)[m_1(X)-\Tilde{m}_1(X;\gamma^*_1)]|\alpha^Tx_i\right\}}{p^2(x_i)}[p(x_i)-d_i]\\
&\quad+o_p(1).
\end{split}
\end{equation*}
Similarly, we can derive the form of $\sqrt{n}(\hat{\theta}_0-\theta_0)$. As a result,
\begin{equation*}
\begin{split}
&\sqrt{n}(\hat{\Delta}_5-\Delta)\\
&=\sqrt{n}[(\hat{\theta}_1-\hat{\theta}_0)-(\theta_1-\theta_0)]\\
&=\frac{1}{\sqrt{n}}\sum_{i=1}^{n}\left\{\frac{d_iy_i}{p(x_i)}+\left(1-\frac{d_i}{p(x_i)}\right)\Tilde{m}_1(x_i;\gamma^*_1)-\frac{(1-d_i)y_i}{1-p(x_i)}-\left(1-\frac{1-d_i}{1-p(x_i)}\right)\Tilde{m}_0(x_i;\gamma^*_0)-(\theta_1-\theta_0)\right\}\\
&\quad+\frac{1}{\sqrt{n}}\sum_{i=1}^n\frac{\mathbb{E}\left\{p(X)[m_1(X)-\Tilde{m}_1(X;\gamma^*_1)]|\alpha^Tx_i\right\}}{p^2(x_i)}[p(x_i)-d_i]\\
&\quad+\frac{1}{\sqrt{n}}\sum_{i=1}^n\frac{\mathbb{E}\left\{(1-p(X))[m_0(X)-\Tilde{m}_0(X;\gamma^*_0)]|\alpha^Tx_i\right\}}{(1-p(x_i))^2}[p(x_i)-d_i]\\
&\quad+o_p(1).
\end{split}
\end{equation*}

Recall the definition of locally misspecified $OR$ models in (\ref{localdef}), we have $m_1(x)=\Tilde{m}_1(x;\gamma_{1,0})+\delta_1\times s_1(x), m_0(x)=\Tilde{m}_0(x;\gamma_{0,0})+\delta_0 \times s_0(x)$. If the $OR$ models are locally misspecified, as shown in Section 6.2.1, the second and the third terms can be written as as $O(\delta_1)$ and $O(\delta_0)$ respectively by Taylor series expansion. These two terms converge to 0 as $\delta_1\to0$ and $\delta_0\to0$. Let $\Phi(x_i,y_i,d_i):=\frac{d_iy_i}{p(x_i)}+\left(1-\frac{d_i}{p(x_i)}\right)\Tilde{m}_1(x_i;\gamma^*_1)-\frac{(1-d_i)y_i}{1-p(x_i)}-\left(1-\frac{1-d_i}{1-p(x_i)}\right)\Tilde{m}_0(x_i;\gamma^*_0)-(\theta_1-\theta_0)$. Note that $\mathbb{E}\left\{\Phi(X,Y,D)\right\}=0$. We further assume $\mathbb{E}\left\{\Phi(X,Y,D)^2\right\} < \infty$. It follows from the Central Limit Theorem and Slutsky's Theorem that the first term converges in distribution to $N(0, \Sigma)$ with $$\Sigma=\mathrm{Var}\left[\frac{DY}{p(X)}+\left(1-\frac{D}{p(X)}\right)\Tilde{m}_1(X;\gamma^*_1)-\frac{(1-D)Y}{1-p(X)} -\left(1-\frac{1-D}{1-p(X)}\right)\Tilde{m}_0(X;\gamma^*_0)-(\theta_1-\theta_0)\right].$$ Note that $\Sigma$ converges to $\Sigma_1$ as $\delta_1\to0$ and $\delta_0\to0$. Consequently, the asymptotic variance of $\hat{\Delta}_5$ converges to $\Sigma_1$.

If $OR$ models are globally misspecified, we can further write $\sqrt{n}(\hat{\Delta}_5-\Delta)$ as
\begin{equation*}
\begin{split}
&\sqrt{n}(\hat{\Delta}_5-\Delta)\\
&=\frac{1}{\sqrt{n}}\sum_{i=1}^{n}\left\{\frac{d_iy_i}{p(x_i)}+\left(1-\frac{d_i}{p(x_i)}\right)\Tilde{m}_1(x_i;\gamma^*_1)-\frac{(1-d_i)y_i}{1-p(x_i)}-\left(1-\frac{1-d_i}{1-p(x_i)}\right)\Tilde{m}_0(x_i;\gamma^*_0)-(\theta_1-\theta_0)\right\}\\
&\qquad+w(x_i)[p(x_i)-d_i]\bigg\} +o_p(1)\\
&=\frac{1}{\sqrt{n}}\sum_{i=1}^{n}\Phi(x_i,y_i,d_i)+o_p(1),
\end{split}
\end{equation*}
where $w(x_i)=\frac{\mathbb{E}\{p(X)[m_1(X)-\Tilde{m}_1(X;\gamma^*_1)]|\alpha^Tx_i\}}{p^2(x_i)}+\frac{\mathbb{E}\{(1-p(X))[m_0(X)-\Tilde{m}_0(X;\gamma^*_0)]|\alpha^Tx_i\}}{(1-p(x_i))^2}$.

\noindent Note that $\mathbb{E}\left\{\Phi(X,Y,D)\right\}=0$. We further assume $\mathbb{E}\left\{\Phi(X,Y,D)^2\right\} < \infty$. It follows from the Central Limit Theorem and Slutsky's Theorem that $\sqrt{n}(\hat{\Delta}_5-\Delta)$ converges in distribution to $N(0,\Sigma)$ with
\begin{equation*}
\begin{split}
\Sigma &=\mathbb{E}\left\{\Phi(X,Y,D)^2\right\}\\
&=\Sigma_1+\mathbb{E}\bigg\{\Big[\sqrt{\frac{1}{p(X)}-1}[\Tilde{m}_1(X;\gamma^*_1)-m_1(X)]+\sqrt{\frac{1}{1-p(X)}-1}[\Tilde{m}_0(X;\gamma^*_0)-m_0(X)]\\
&\qquad+\sqrt{p(X)(1-p(X))}w(X)\Big]^2\bigg\} \geq \Sigma_1.
\end{split}
\end{equation*}

\noindent Therefore, the asymptotic variance of $\hat{\Delta}_5$ is enlarged compared to the semiparametric efficiency bound $\Sigma_1$ when $OR$ models are globally misspecified.

\newcontent (c) \textbf{Misspecified $PS$ model and correctly specified $OR$ model}

In this case, we have $p(x)\neq g(\alpha^Tx)$, $m_1(x)=\Tilde{m}_1(x;\gamma_{1,0})=\Tilde{m}_1(x;\gamma^*_1)$ and  $m_0(x)=\Tilde{m}_0(x;\gamma_{0,0})=\Tilde{m}_0(x;\gamma^*_0)$. So $\sqrt{n}B_n=o_p(1)$. As shown in Section 6.2.1,
\begin{equation*}
\begin{split}
&\sqrt{n}(\hat{\gamma}_1-\gamma_{1,0})=\frac{1}{\sqrt{n}}\sum^{n}_{i=1}\frac{1}{\sigma^2_{(1)}}I^{-1}(\gamma_{1,0})x_id_i[y_i-m_1(x_i)]+o_p(1),\\
&\sqrt{n}(\hat{\gamma}_0-\gamma_{0,0})=\frac{1}{\sqrt{n}}\sum^{n}_{i=1}\frac{1}{\sigma^2_{(0)}}I^{-1}(\gamma_{0,0})x_i(1-d_i)[y_i-m_0(x_i)]+o_p(1).
\end{split}
\end{equation*}
Combining the terms in (\ref{6.2.5}), we have
\begin{equation*}
\begin{split}
&\sqrt{n}(\hat{\theta}_1-\theta_1)\\
&=\frac{1}{\sqrt{n}}\sum_{i=1}^{n}\left\{\frac{d_iy_i}{g(\alpha^Tx_i)}+\left(1-\frac{d_i}{g(\alpha^Tx_i)}\right)m_1(x_i)-\theta_1\right\} \\
&\quad+\sqrt{n}(\hat{\gamma}_1-\gamma_{1,0})^T \mathbb{E}\left[(1-\frac{p(X)}{g(\alpha^TX)})\frac{\partial \Tilde{m}_1(X;\gamma_1)}{\partial \gamma_1}\bigg|_{\gamma_1=\bar{\gamma}_1}\right]\\
&\quad+o_p(1).
\end{split}
\end{equation*}
Similarly, we can derive $\sqrt{n}(\hat{\theta}_0-\theta_0)$. As a result,
\begin{equation*}
\begin{split}
&\sqrt{n}(\hat{\Delta}_5-\Delta)\\
&=\sqrt{n}[(\hat{\theta}_1-\hat{\theta}_0)-(\theta_1-\theta_0)]\\
&=\frac{1}{\sqrt{n}}\sum_{i=1}^{n}\bigg\{\frac{d_iy_i}{g(\alpha^Tx_i)}+\left(1-\frac{d_i}{g(\alpha^Tx_i)}\right)m_1(x_i)-\frac{(1-d_i)y_i}{1-g(\alpha^Tx_i)}-\left(1-\frac{1-d_i}{1-g(\alpha^Tx_i)}\right)m_0(x_i)-(\theta_1-\theta_0)\\
&\quad+w_1(x_i)d_i[y_i-m_1(x_i)]+w_0(x_i)(1-d_i)[y_1-m_0(x_i)]\bigg\}+o_p(1)\\
&=\frac{1}{\sqrt{n}}\sum_{i=1}^{n}\Phi(x_i,y_i,d_i)+o_p(1),
\end{split}
\end{equation*}
where
\begin{equation*}
\begin{split}
&w_1(x_i)=\frac{1}{\sigma^2_{(1)}}\left(I^{-1}(\gamma_{1,0})x_i\right)^T\mathbb{E}\left[(1-\frac{p(X)}{g(\alpha^TX)})\frac{\partial \Tilde{m}_1(X;\gamma_1)}{\partial \gamma_1}\big|_{\gamma_1=\bar{\gamma}_1}\right],  \\
&w_0(x_i)=\frac{1}{\sigma^2_{(0)}}\left(I^{-1}(\gamma_{0,0})x_i\right)^T\mathbb{E}\left[(\frac{1-p(X)}{1-g(\alpha^TX)}-1)\frac{\partial \Tilde{m}_0(X;\gamma_0)}{\partial \gamma_0}\big|_{\gamma_0=\bar{\gamma}_0}\right].
\end{split}
\end{equation*}

\noindent Note that $\mathbb{E}\left\{\Phi(X,Y,D)\right\}=0$. We further assume $\mathbb{E}\left\{\Phi(X,Y,D)^2\right\} < \infty$. It follows from the Central Limit Theorem and Slutsky's Theorem that $\sqrt{n}(\hat{\Delta}_5-\Delta) $ converges in distribution to $N(0,\Sigma)$, where
\begin{equation*}
\begin{split}
\Sigma&=\mathbb{E}\left\{\Phi(X,Y,D)^2\right\}\\
&=\Sigma_1+\mathbb{E}\left\{\frac{1}{p(X)}\mathrm{Var}[Y(1)|X]\left[\left(\frac{p(X)}{g(\alpha^TX)}+w_1(X)p(X)\right)^2-1\right]\right\}\\
&\quad+\mathbb{E}\left\{\frac{1}{1-p(X)}\mathrm{Var}[Y(0)|X]\left[\left(\frac{1-p(X)}{1-g(\alpha^TX)}+w_0(X)(1-p(X))\right)^2-1\right]\right\}.
\end{split}
\end{equation*}

\noindent Therefore, the asymptotic variance of $\hat{\Delta}_5$ is not necessarily enlarged compared to the semiparametric efficiency bound $\Sigma_1$ when $PS$ model is misspecified.

\subsubsection{Parametric $PS$ model and semiparametric $OR$ model}

When $PS$ model is parametric and $OR$ model is semiparametric, we have
\begin{eqnarray}\label{6.2.6}
\sqrt{n}(\hat{\theta}_1-\theta_1) &=& \sqrt{n}\left\{n^{-1}\sum_{i=1}^{n}\left[\frac{d_iy_i}{\Tilde{p}(x_i;\hat{\beta})}+\left(1-\frac{d_i}{\Tilde{p}(x_i;\hat{\beta})}\right)\hat{r}_1(\alpha^T_1x_i)\right]-\theta_1\right\} \nonumber\\
&=&\frac{1}{\sqrt{n}}\sum_{i=1}^{n}\left\{\frac{d_iy_i}{\Tilde{p}(x_i;\beta^*)}+\left(1-\frac{d_i}{\Tilde{p}(x_i;\beta^*)}\right)r_1(\alpha^T_1x_i)-\theta_1\right\} \nonumber\\
&&+\frac{1}{\sqrt{n}}\sum_{i=1}^{n}d_i\left(\frac{1}{\Tilde{p}(x_i;\hat{\beta})}-\frac{1}{\Tilde{p}(x_i;\beta^*)}\right)[y_i-r_1(\alpha^T_1x_i)] \nonumber\\
&&+\frac{1}{\sqrt{n}}\sum_{i=1}^{n}\left(1-\frac{d_i}{\Tilde{p}(x_i;\beta^*)}\right)[\hat{r}_1(\alpha^T_1x_i)-r_1(\alpha^T_1x_i)] \nonumber\\
&&+\frac{1}{\sqrt{n}}\sum_{i=1}^{n}d_i\left(\frac{1}{\Tilde{p}(x_i;\beta^*)}-\frac{1}{\Tilde{p}(x_i;\hat{\beta})}\right)[\hat{r}_1(\alpha^T_1x_i)-r_1(\alpha^T_1x_i)] \nonumber\\
&:=&\sqrt{n}A_n + \sqrt{n}B_n + \sqrt{n}C_n + \sqrt{n}D_n.
\end{eqnarray}

\noindent Recall that $\bar{\beta}$ is defined in Section 6.2.1. Similar to the derivations of $\sqrt{n}B_n$ in Section 6.2.1, we can get
\begin{equation*}
\sqrt{n}B_n=\sqrt{n}(\hat{\beta}-\beta^*)^T\mathbb{E}\left\{\frac{\partial\frac{1}{\Tilde{p}(X;\beta)}}{\partial \beta}\bigg|_{\beta=\bar{\beta}}p(X)[m_1(X)-r_1(\alpha^T_1X)]\right\}+o_p(1),
\end{equation*}
and
\begin{equation*}
\begin{split}
\sqrt{n}D_n&=\frac{1}{\sqrt{n}}\sum_{i=1}^{n}d_i\left(\frac{1}{\Tilde{p}(x_i;\beta^*)}-\frac{1}{\Tilde{p}(x_i;\hat{\beta})}\right)[\hat{r}_1(\alpha^T_1x_i)-r_1(\alpha^T_1x_i)]\\
&=\sqrt{n}(\beta^*-\hat{\beta})^T\frac{1}{n}\sum_{i=1}^{n}d_i\frac{\partial\frac{1}{\Tilde{p}(x_i;\beta)}}{\partial \beta}\bigg|_{\beta=\bar{\beta}}\left[\frac{\sum_{j=1}^{n}d_jy_jK_{h_{m_1}}(\alpha^T_1x_i,\alpha^T_1x_j)}{\sum_{j=1}^{n}d_jK_{h_{m_1}}(\alpha^T_1x_i,\alpha^T_1x_j)}-r_1(\alpha^T_1x_i)\right]\\
&=\sqrt{n}(\beta^*-\hat{\beta})^T\frac{1}{n}\sum_{i=1}^{n}d_i\frac{\partial\frac{1}{\Tilde{p}(x_i;\beta)}}{\partial \beta}\bigg|_{\beta=\bar{\beta}}\frac{\sum_{j=1}^{n}d_j[y_j-r_1(\alpha^T_1x_i)]K_{h_{m_1}}(\alpha^T_1x_i,\alpha^T_1x_j)}{n\frac{\sum_{j=1}^{n}d_jK_{h_{m_1}}(\alpha^T_1x_i,\alpha^T_1x_j)}{\sum_{j=1}^{n}K_{h_{m_1}}(\alpha^T_1x_i,\alpha^T_1x_j)}\frac{1}{n}\sum_{j=1}^{n}d_jK_{h_{m_1}}(\alpha^T_1x_i,\alpha^T_1x_j)}\\
&=\sqrt{n}(\beta^*-\hat{\beta})^T\frac{1}{n}\sum_{i=1}^{n}d_i\frac{\partial\frac{1}{\Tilde{p}(x_i;\beta)}}{\partial \beta}\bigg|_{\beta=\bar{\beta}}\frac{\sum_{j=1}^{n}d_j[y_j-r_1(\alpha^T_1x_i)]K_{h_{m_1}}(\alpha^T_1x_i,\alpha^T_1x_j)}{n\hat{\Tilde{f}}_1(\alpha^T_1x_i)\hat{q}_1(\alpha^T_1x_i)},
\end{split}
\end{equation*}
where $\hat{q}_1(\alpha^T_1x)$ is a semiparametric estimation of $PS$ model and $\hat{\Tilde{f}}_1(\alpha^T_1x)$ is the kernel density estimator of $\Tilde{f}_1(\alpha^T_1x)$. Under Assumption \ref{semiassumption}, we have $\sup_{x \in \mathcal{X}}|\hat{q}_1(\alpha^T_1x)-q_1(\alpha^T_1x)|=O_p(h^2_{m_1}+\sqrt{\frac{\log(n)}{nh_{m_1}}})$ and $\sup_{x \in \mathcal{X}}|\hat{\Tilde{f}}_1(\alpha^T_1x)-\Tilde{f}_1(\alpha^T_1x)|=O_p(h^2_{m_1}+\sqrt{\frac{\log(n)}{nh_{m_1}}})$ (see \cite{nonparametric-reference1}, \cite{nonparametric-reference2} and \cite{nonparametric-reference3}). Due to the consistency of $\hat{q}_1(\alpha^T_1x)$ to $q_1(\alpha^T_1x)$ and $\hat{\Tilde{f}}_1(\alpha^T_1x)$ to $\Tilde{f}_1(\alpha^T_1x)$ , we obtain
\begin{equation*}
\begin{split}
\sqrt{n}D_n&=\sqrt{n}(\beta^*-\hat{\beta})^T\frac{1}{n}\sum_{i=1}^{n}d_i\frac{\partial\frac{1}{\Tilde{p}(x_i;\beta)}}{\partial \beta}\bigg|_{\beta=\bar{\beta}}\frac{\sum_{j=1}^{n}d_j[y_j-r_1(\alpha^T_1x_i)]K_{h_{m_1}}(\alpha^T_1x_i,\alpha^T_1x_j)}{n\Tilde{f}_1(\alpha^T_1x_i)q_1(\alpha^T_1x_i)}+o_p(1)\\
&=\sqrt{n}(\beta^*-\hat{\beta})^T\frac{1}{n^2}\sum_{i=1}^{n}\sum_{j=1}^{n}\frac{d_i\frac{\partial\frac{1}{\Tilde{p}(x_i;\beta)}}{\partial \beta}\Big|_{\beta=\bar{\beta}}}{\Tilde{f}_1(\alpha^T_1x_i)q_1(\alpha^T_1x_i)}d_j[y_j-r_1(\alpha^T_1x_i)]K_{h_{m_1}}(\alpha^T_1x_i,\alpha^T_1x_j)+o_p(1).
\end{split}
\end{equation*}

\noindent We further write $\sqrt{n}D_n$ in the form of U-statistics. Let
$H_{ij}=\frac{d_i\frac{\partial\frac{1}{\Tilde{p}(x_i;\beta)}}{\partial \beta}\big|_{\beta=\bar{\beta}}}{\Tilde{f}_1(\alpha^T_1x_i)q_1(\alpha^T_1x_i)}d_j[y_j-r_1(\alpha^T_1x_i)]$, we have
\begin{equation*}
\begin{split}
\sqrt{n}D_n&=\sqrt{n}(\beta^*-\hat{\beta})^T\frac{n-1}{n}\frac{1}{n(n-1)}\sum_{i=1}^{n}\sum_{j\neq i}^{n} \frac{H_{ij}+H{ji}}{2}K_{h_{m_1}}(\alpha^T_1x_i,\alpha^T_1x_j)+o_p(1)\\
&=\sqrt{n}(\beta^*-\hat{\beta})^T\frac{n-1}{n}U_n+o_p(1),
\end{split}
\end{equation*}
where $U_n:=\frac{1}{n(n-1)}\sum_{i=1}^{n}\sum_{j\neq i}^{n} \frac{H_{ij}+H{ji}}{2}K_{h_{m_1}}(\alpha^T_1x_i,\alpha^T_1x_j)$. Next, we compute the conditional expectation of $\left[\frac{H_{ij}+H_{ji}}{2}\right]K_{h_{m_1}}(\alpha^T_1x_i,\alpha^T_1x_j)$. We first compute
\begin{equation*}
\begin{split}
&\mathbb{E}[H_{ij}K_{h_{m_1}}(\alpha^T_1x_i,\alpha^T_1x_j)|x_j,y_j,d_j]\\
&=\mathbb{E}\left\{\frac{d_i\frac{\partial\frac{1}{\Tilde{p}(x_i;\beta)}}{\partial \beta}\Big|_{\beta=\bar{\beta}}}{\Tilde{f}_1(\alpha^T_1x_i)q_1(\alpha^T_1x_i)}d_j[y_j-r_1(\alpha^T_1x_i)]K_{h_{m_1}}(\alpha^T_1x_i,\alpha^T_1x_j)\bigg|x_j,y_j,d_j\right\}\\
&=\mathbb{E}\left\{\mathbb{E}\left[\frac{d_i\frac{\partial\frac{1}{\Tilde{p}(x_i;\beta)}}{\partial \beta}\Big|_{\beta=\bar{\beta}}}{\Tilde{f}_1(\alpha^T_1x_i)q_1(\alpha^T_1x_i)}d_j[y_j-r_1(\alpha^T_1x_i)]K_{h_{m_1}}(\alpha^T_1x_i,\alpha^T_1x_j)\bigg|\alpha^T_1x_i,x_j,y_j,d_j\right]\bigg|x_j,y_j,d_j\right\}\\
&=E\left\{\frac{E\left[D\frac{\partial\frac{1}{\Tilde{p}(X;\beta)}}{\partial \beta}|_{\beta=\bar{\beta}}|\alpha^T_1x_i\right]}{\Tilde{f}_1(\alpha^T_1x_i)q_1(\alpha^T_1x_i)}d_j[y_j-r_1(\alpha^T_1x_i)]K_{h_{m_1}}(\alpha^T_1x_i,\alpha^T_1x_j)\bigg|x_j,y_j,d_j\right\}\\
&=\int w(\alpha^T_1x)d_j[y_j-r_1(\alpha^T_1x)]K_{h_{m_1}}(\alpha^T_1x,\alpha^T_1x_j)\Tilde{f}_1(\alpha^T_1x)d(\alpha^T_1x).
\end{split}
\end{equation*}
Let $t=\frac{\alpha^T_1x-\alpha^T_1x_j}{h_{m_1}}$ and further apply Taylor series expansion. Based on Assumption \ref{semiassumption}, we obtain
\begin{equation*}
\begin{split}
&\mathbb{E}[H_{ij}K_{h_{m_1}}(\alpha^T_1x_i,\alpha^T_1x_j)|x_j,y_j,d_j]\\
&=\int w(\alpha^T_1x_j+th_{m_1})d_j[y_j-r_1(\alpha^T_1x_j+th_{m_1})]K(t)\Tilde{f}_1(\alpha^T_1x_j+th_{m_1})dt\\
&=w(\alpha^T_1x_j)d_j[y_j-r_1(\alpha^T_1x_j)]\Tilde{f}_1(\alpha^T_1x_j)+O(h^2_{m_1})\\
&=\frac{\mathbb{E}\left[\mathbb{E}\left(D\frac{\partial\frac{1}{\Tilde{p}(X;\beta)}}{\partial \beta}\Big|_{\beta=\bar{\beta}}\Big|x_j\right)\Big|\alpha^T_1x_j\right]}{q_1(\alpha^T_1x_j)}d_j[y_j-r_1(\alpha^T_1x_j)]+O(h^2_{m_1})\\
&=\frac{\mathbb{E}\left[p(X)\frac{\partial\frac{1}{\Tilde{p}(X;\beta)}}{\partial \beta}\Big|_{\beta=\bar{\beta}}\Big|\alpha^T_1x_j\right]}{q_1(\alpha^T_1x_j)}d_j[y_j-r_1(\alpha^T_1x_j)]+O(h^2_{m_1}).
\end{split}
\end{equation*}

\noindent Similarly, we have
\begin{align*}
&\mathbb{E}[H_{ji}K_{h_{m_1}}(\alpha^T_1x_i,\alpha^T_1x_j)|x_j,y_j,d_j]\\
&=\mathbb{E}\left\{\frac{d_j\frac{\partial\frac{1}{\Tilde{p}(x_j;\beta)}}{\partial \beta}\Big|_{\beta=\bar{\beta}}}{\Tilde{f}_1(\alpha^T_1x_j)q_1(\alpha^T_1x_j)}d_i[y_i-r_1(\alpha^T_1x_j)]K_{h_{m_1}}(\alpha^T_1x_i,\alpha^T_1x_j)\bigg|x_j,y_j,d_j\right\}\\
&=\frac{d_j\frac{\partial\frac{1}{\Tilde{p}(x_j;\beta)}}{\partial \beta}\Big|_{\beta=\bar{\beta}}}{\Tilde{f}_1(\alpha^T_1x_j)q_1(\alpha^T_1x_j)}\mathbb{E}\left\{\mathbb{E}[d_i[y_i-r_1(\alpha^T_1x_j)]K_{h_{m_1}}(\alpha^T_1x_i,\alpha^T_1x_j)|\alpha^T_1x_i,x_j,y_j,d_j]|x_j,y_j,d_j\right\}\\
&=\frac{d_j\frac{\partial\frac{1}{\Tilde{p}(x_j;\beta)}}{\partial \beta}\Big|_{\beta=\bar{\beta}}}{\Tilde{f}_1(\alpha^T_1x_j)q_1(\alpha^T_1x_j)}\mathbb{E}\left\{\mathbb{E}(DY|\alpha^T_1x_i)K_{h_{m_1}}(\alpha^T_1x_i,\alpha^T_1x_j)-E(D|\alpha^T_1x_i)r_1(\alpha^T_1x_j)K_{h_{m_1}}(\alpha^T_1x_i,\alpha^T_1x_j)|x_j,y_j,d_j\right\}\\
&=\frac{d_j\frac{\partial\frac{1}{\Tilde{p}(x_j;\beta)}}{\partial \beta}\Big|_{\beta=\bar{\beta}}}{\Tilde{f}_1(\alpha^T_1x_j)q_1(\alpha^T_1x_j)}\bigg\{\int w_1(\alpha^T_1x)K_{h_{m_1}}(\alpha^T_1x,\alpha^T_1x_j)\Tilde{f}_1(\alpha^T_1x)d(\alpha^T_1x)\\
&\qquad-\int w_2(\alpha^T_1x)r_1(\alpha^T_1x_j)K_{h_{m_1}}(\alpha^T_1x,\alpha^T_1x_j)\Tilde{f}_1(\alpha^T_1x)d(\alpha^T_1x)\bigg\}\\
&=\frac{d_j\frac{\partial\frac{1}{\Tilde{p}(x_j;\beta)}}{\partial \beta}\Big|_{\beta=\bar{\beta}}}{\Tilde{f}_1(\alpha^T_1x_j)q_1(\alpha^T_1x_j)}\bigg\{\int w_1(\alpha^T_1x_j+th_{m_1})K(t)\Tilde{f}_1(\alpha^T_1x_j+th_{m_1})dt\\
&\qquad-\int w_2(\alpha^T_1x_j+th_{m_1})r_1(\alpha^T_1x_j)K(t)\Tilde{f}_1(\alpha^T_1x_j+th_{m_1})dt\bigg\}\\
&=\frac{d_j\frac{\partial\frac{1}{\Tilde{p}(x_j;\beta)}}{\partial \beta}\Big|_{\beta=\bar{\beta}}}{\Tilde{f}_1(\alpha^T_1x_j)q_1(\alpha^T_1x_j)}\left[w_1(\alpha^T_1x_j)\Tilde{f}_1(\alpha^T_1x_j)-w_2(\alpha^T_1x_j)r_1(\alpha^T_1x_j)\Tilde{f}_1(\alpha^T_1x_j)\right]+O(h^2_{m_1})\\
&=\frac{d_j\frac{\partial\frac{1}{\Tilde{p}(x_j;\beta)}}{\partial \beta}\Big|_{\beta=\bar{\beta}}}{q_1(\alpha^T_1x_j)}\mathbb{E}\left\{\mathbb{E}\left(D[Y-r_1(\alpha^T_1x_j)]|x_j\right)|\alpha^T_1x_j\right\}+O(h^2_{m_1})\\
&=\frac{d_j\frac{\partial\frac{1}{\Tilde{p}(x_j;\beta)}}{\partial \beta}|_{\beta=\bar{\beta}}}{q_1(\alpha^T_1x_j)}\left\{E\big[DY|\alpha^T_1x_j\big]-E\left[D|\alpha^T_1x_j\right]r_1(\alpha^T_1x_j)\right\}+O(h^2_{m_1})\\
&=\frac{d_j\frac{\partial\frac{1}{\Tilde{p}(x_j;\beta)}}{\partial \beta}|_{\beta=\bar{\beta}}}{q_1(\alpha^T_1x_j)}\left\{E\big[DY|\alpha^T_1x_j, D=1\big]P(D=1|\alpha^T_1x_j)-q_1(\alpha^T_1x_j)r_1(\alpha^T_1x_j)\right\}+O(h^2_{m_1})\\
&=\frac{d_j\frac{\partial\frac{1}{\Tilde{p}(x_j;\beta)}}{\partial \beta}|_{\beta=\bar{\beta}}}{q_1(\alpha^T_1x_j)}\left\{E\big[Y(1)|\alpha^T_1x_j, D=1\big]q_1(\alpha^T_1x_j)-q_1(\alpha^T_1x_j)r_1(\alpha^T_1x_j)\right\}+O(h^2_{m_1})\\
&=\frac{d_j\frac{\partial\frac{1}{\Tilde{p}(x_j;\beta)}}{\partial \beta}|_{\beta=\bar{\beta}}}{q_1(\alpha^T_1x_j)}\left\{r_1(\alpha^T_1x_j)q_1(\alpha^T_1x_j)-q_1(\alpha^T_1x_j)r_1(\alpha^T_1x_j)\right\}+O(h^2_{m_1})\\
&=O(h^2_{m_1})
\end{align*}
\noindent The conditional expectation of $\left[\frac{H_{ij}+H_{ji}}{2}\right]K_{h_{m_1}}(\alpha^T_1x_i,\alpha^T_1x_j)$ is $\frac{\mathbb{E}\left[p(X)\frac{\partial\frac{1}{\Tilde{p}(X;\beta)}}{\partial \beta}\Big|_{\beta=\bar{\beta}}\Big|\alpha^T_1x_j\right]}{2q_1(\alpha^T_1x_j)}d_j[y_j-r_1(\alpha^T_1x_j)]+O(h^2_{m_1})$. It follows that $\mathbb{E}\left[\frac{H_{ij}+H_{ji}}{2}K_{h_{m_1}}(\alpha^T_1x_i,\alpha^T_1x_j)\right]=O(h^2_{m_1})$. Then we can calculate the projection of $U_n$. Based on Assumption \ref{semiassumption}, we have $E\big[||\frac{H_{ij}+H_{ji}}{2}K_{h_{m_1}}(\alpha^T_1x_i,\alpha^T_1x_j)||^2\big]=o(n)$. Applying Lemma 3.1 of \cite{Ustatistics} under Assumption \ref{semiassumption}, we have $\sqrt{n}D_n=o_p(1)$. For $\sqrt{n}C_n$, following the derivations of $\sqrt{n}D_n$, we have
\begin{align*}
\sqrt{n}C_n &= \frac{1}{\sqrt{n}}\sum_{i=1}^{n}\left(1-\frac{d_i}{\Tilde{p}(x_i;\beta^*)}\right)[\hat{r}_1(\alpha^T_1x_i)-r_1(\alpha^T_1x_i)]\\
&=\frac{1}{\sqrt{n}}\sum_{i=1}^{n}\left(1-\frac{d_i}{\Tilde{p}(x_i;\beta^*)}\right)\left[\frac{\sum_{j=1}^{n}d_jy_jK_{h_{m_1}}(\alpha^T_1x_i,\alpha^T_1x_j)}{\sum_{j=1}^{n}d_jK_{h_{m_1}}(\alpha^T_1x_i,\alpha^T_1x_j)}-r_1(\alpha^T_1x_i)\right]\\
&=\frac{1}{\sqrt{n}}\sum_{i=1}^{n}\left(1-\frac{d_i}{\Tilde{p}(x_i;\beta^*)}\right)\frac{\sum_{i=1}^{n}d_j(y_j-r_1(\alpha^T_1x_i))K_{h_{m_1}}(\alpha^T_1x_i,\alpha^T_1x_j)}{\sum_{j=1}^{n}d_jK_{h_{m_1}}(\alpha^T_1x_i,\alpha^T_1x_j)}\\
&=\frac{1}{\sqrt{n}}\sum_{i=1}^{n}\left(1-\frac{d_i}{\Tilde{p}(x_i;\beta^*)}\right)\frac{\sum_{i=1}^{n}d_j(y_j-r_1(\alpha^T_1x_i))K_{h_{m_1}}(\alpha^T_1x_i,\alpha^T_1x_j)}{n\frac{\sum_{j=1}^{n}d_jK_{h_{m_1}}(\alpha^T_1x_i,\alpha^T_1x_j)}{\sum_{j=1}^{n}K_{h_{m_1}}(\alpha^T_1x_i,\alpha^T_1x_j)} \frac{1}{n}\sum_{j=1}^{n}K_{h_{m_1}}(\alpha^T_1x_i,\alpha^T_1x_j)}\\
&=\frac{1}{\sqrt{n}}\sum_{i=1}^{n}\left(1-\frac{d_i}{\Tilde{p}(x_i;\beta^*)}\right)\frac{\sum_{i=1}^{n}d_j(y_j-r_1(\alpha^T_1x_i))K_{h_{m_1}}(\alpha^T_1x_i,\alpha^T_1x_j)}{n\hat{\Tilde{f}}(\alpha^T_1x_i)\hat{q}_1(\alpha^T_1x_i)}\\
&=\frac{1}{\sqrt{n}}\sum_{i=1}^{n}\left(1-\frac{d_i}{\Tilde{p}(x_i;\beta^*)}\right)\frac{\sum_{i=1}^{n}d_j(y_j-r_1(\alpha^T_1x_i))K_{h_{m_1}}(\alpha^T_1x_i,\alpha^T_1x_j)}{nf(\alpha^T_1x_i)q_1(\alpha^T_1x_i)}+o_p(1)\\
&=\frac{1}{n\sqrt{n}}\sum_{i=1}^{n}\sum_{j=1}^{n}\frac{1-\frac{d_i}{\Tilde{p}(x_i;\beta^*)}}{f(\alpha^T_1x_i)q_1(\alpha^T_1x_i)}d_j[y_j-r_1(\alpha^T_1x_i)]K_{h_{m_1}}(\alpha^T_1x_i,\alpha^T_1x_j)+o_p(1)\\
&=\frac{1}{\sqrt{n}}\sum_{j=1}^n\frac{E\big[1-\frac{p(X)}{\Tilde{p}(X;\beta^*)}|\alpha^T_1x_j\big]}{q_1(\alpha^T_1x_j)}d_j[y_j-r_1(\alpha^T_1x_j)]+O(\sqrt{n}h^2_{m_1})+o_p(1)
\end{align*}
\noindent Under Assumption \ref{semiassumption}, we have $\sqrt{n}C_n=\frac{1}{\sqrt{n}}\sum_{j=1}^n\frac{E\big[1-\frac{p(X)}{\Tilde{p}(X;\beta^*)}|\alpha^T_1x_j\big]}{q_1(\alpha^T_1x_j)}d_j[y_j-r_1(\alpha^T_1x_j)]+o_p(1)$. Now we can consider different cases as follows.

\newcontent (a) \textbf{Correctly specified $PS$ model and $OR$ model}

In this case, we have $p(x)=\Tilde{p}(x;\beta_0)=\Tilde{p}(x;\beta^*)$, $m_1(x)=r_1(\alpha^T_1x)$ and $m_0(x)=r_0(\alpha^T_0x)$. So $\sqrt{n}B_n=o_p(1)$ and $\sqrt{n}C_n=o_p(1)$. Combining the terms in (\ref{6.2.6}), we have
\begin{equation*}
\begin{split}
\sqrt{n}(\hat{\theta}_1-\theta_1)=\frac{1}{\sqrt{n}}\sum_{i=1}^{n}\left\{\frac{d_iy_i}{p(x_i)}+\left(1-\frac{d_i}{p(x_i)}\right)m_1(x_i)-\theta_1\right\}+o_p(1).
\end{split}
\end{equation*}
Similarly, we can derive the form of $\sqrt{n}(\hat{\theta}_0-\theta_0)$. As a result,
\begin{equation*}
\begin{split}
&\sqrt{n}(\hat{\Delta}_6-\Delta) \\
&=\sqrt{n}[(\hat{\theta}_1-\hat{\theta}_0)-(\theta_1-\theta_0)]\\
&=\frac{1}{\sqrt{n}}\sum_{i=1}^{n}\left\{\frac{d_iy_i}{p(x_i)}+\left(1-\frac{d_i}{p(x_i)}\right)m_1(x_i)-\frac{(1-d_i)y_i}{1-p(x_i)}-\left(1-\frac{(1-d_i)}{1-p(x_i)}\right)m_0(x_i)-(\theta_1-\theta_0)\right\}+o_p(1)\\
&=\frac{1}{\sqrt{n}}\sum_{i=1}^{n}\Phi(x_i,y_i,d_i)+o_p(1).
\end{split}
\end{equation*}
\noindent Note that $\mathbb{E}\left\{\Phi(X,Y,D)\right\}=0$. We further assume $\mathbb{E}\left\{\Phi(X,Y,D)^2\right\} < \infty$. It follows from the Central Limit Theorem and Slutsky's Theorem that $\sqrt{n}(\hat{\Delta}_6-\Delta)$ converges in distribution to $N(0,\Sigma_1)$. The asymptotic variance of $\hat{\Delta}_6$ achieves the semiparametric efficiency bound $\Sigma_1$ when both $PS$ and $OR$ models are correctly specified.

\newcontent (b) \textbf{Correctly specified $PS$ model and misspecified $OR$ model}

In this case, we have $p(x)=\Tilde{p}(x;\beta_0)=\Tilde{p}(x;\beta^*)$, $m_1(x)\neq r_1(\alpha^T_1x)$ and $m_0(x)\neq r_0(\alpha^T_0x)$. So $\sqrt{n}C_n=o_p(1)$. Combining the terms in (\ref{6.2.6}), we have
\begin{equation*}
\begin{split}
&\sqrt{n}(\hat{\theta}_1-\theta_1)\\
&=\frac{1}{\sqrt{n}}\sum_{i=1}^{n}\left\{\frac{d_iy_i}{p(x_i)}+\left(1-\frac{d_i}{p(x_i)}\right)r_1(\alpha^T_1x_i)-\theta_1\right\}\\
&\quad+\sqrt{n}(\hat{\beta}-\beta^*)^T\mathbb{E}\left\{\frac{\partial\frac{1}{\Tilde{p}(X;\beta)}}{\partial \beta}\bigg|_{\beta=\bar{\beta}}p(X)[m_1(X)-r_1(\alpha^T_1X)]\right\}\\
&\quad+o_p(1).
\end{split}
\end{equation*}
Similarly, the form of $\sqrt{n}(\hat{\theta}_0-\theta_0)$ can be derived. As shown in Section 6.2.1, we have
$$\sqrt{n}(\hat{\beta}-\beta_0)=\frac{1}{\sqrt{n}}\sum^n_{i=1}I^{-1}(\beta_0) x_i\left[d_i-p(x_i)\right]+o_p(1).$$
As a result,
\begin{equation*}
\begin{split}
&\sqrt{n}(\hat{\Delta}_6-\Delta) \\
&=\sqrt{n}[(\hat{\theta}_1-\hat{\theta}_0)-(\theta_1-\theta_0)]\\
&=\frac{1}{\sqrt{n}}\sum_{i=1}^{n}\bigg\{\frac{d_iy_i}{p(x_i)}+\left(1-\frac{d_i}{p(x_i)}\right)r_1(\alpha^T_1x_i)-\frac{(1-d_i)y_i}{1-p(x_i)}-\left(1-\frac{(1-d_i)}{1-p(x_i)}\right)r_0(\alpha^T_0x_i)-(\theta_1-\theta_0)\\
&\qquad+w(x_i)[p(x_i)-d_i]\bigg\}+o_p(1)\\
&=\frac{1}{\sqrt{n}}\sum_{i=1}^{n}\Phi(x_i,y_i,d_i)+o_p(1),
\end{split}
\end{equation*}
where
\begin{equation*}
\begin{split}
w(x_i)&=\left(I^{-1}(\beta_0)x_i\right)^T\Bigg\{\mathbb{E}\Big\{\frac{\partial\frac{1}{1-\Tilde{p}(X;\beta)}}{\partial \beta}\Big|_{\beta=\bar{\beta}}(1-p(X))[m_0(X)-r_0(\alpha^T_0X)]\Big\}\\
&\qquad-\mathbb{E}\Big\{\frac{\partial\frac{1}{\Tilde{p}(X;\beta)}}{\partial \beta}\Big|_{\beta=\bar{\beta}}p(X)[m_1(X)-r_1(\alpha^T_1X)]\Big\}\bigg\}.
\end{split}
\end{equation*}

\noindent Note that $\mathbb{E}\left\{\Phi(X,Y,D)\right\}=0$. We further assume $\mathbb{E}\left\{\Phi(X,Y,D)^2\right\} < \infty$. It follows from the Central Limit Theorem and Slutsky's Theorem that $\sqrt{n}(\hat{\Delta}_6-\Delta)$ converges in distribution to $N(0,\Sigma)$ with
\begin{equation*}
\begin{split}
\Sigma&=\mathbb{E}\left\{\Phi(X,Y,D)^2\right\}\\
&=\Sigma_1+\mathbb{E}\bigg\{\Big[\sqrt{\frac{1}{p(X)}-1}[r_1(\alpha^T_1X)-m_1(X)]+\sqrt{\frac{1}{1-p(X)}-1}[r_0(\alpha^T_0X)-m_0(X)]\\
&\qquad+\sqrt{p(X)(1-p(X))}w(X)\Big]^2\bigg\} \geq \Sigma_1.
\end{split}
\end{equation*}
\noindent Therefore, the asymptotic variance of $\hat{\Delta}_6$ is enlarged compared to the semiparametric efficiency bound $\Sigma_1$ when $OR$ models are misspecified.

\newcontent (c) \textbf{Misspecified $PS$ model and correctly specified $OR$ model}

In this case, we have $p(x)\neq\Tilde{p}(x;\beta_0)\neq\Tilde{p}(x;\beta^*)$, $m_1(x)=r_1(\alpha^T_1x)$ and $m_0(x)=r_0(\alpha^T_0x)$. So $\sqrt{n}B_n=o_p(1)$. Combining the terms in (\ref{6.2.6}),
\begin{equation*}
\begin{split}
&\sqrt{n}(\hat{\theta}_1-\theta_1)\\
&=\frac{1}{\sqrt{n}}\sum_{i=1}^{n}\left\{\frac{d_iy_i}{\Tilde{p}(x_i;\beta^*)}+\left(1-\frac{d_i}{\Tilde{p}(x_i;\beta^*)}\right)m_1(x_i)-\theta_1\right\}\\
&\quad+\frac{1}{\sqrt{n}}\sum_{i=1}^n\frac{\mathbb{E}\left[1-\frac{p(X)}{\Tilde{p}(X;\beta^*)}|\alpha^T_1x_i\right]}{q_1(\alpha^T_1x_i)}d_i[y_i-m_1(x_i)]\\
&\quad+o_p(1).
\end{split}
\end{equation*}
Similarly, the form of $\sqrt{n}(\hat{\theta}_0-\theta_0)$ can be derived. As a result,
\begin{equation*}
\begin{split}
&\sqrt{n}(\hat{\Delta}_6-\Delta)\\
&=\sqrt{n}[(\hat{\theta}_1-\hat{\theta}_0)-(\theta_1-\theta_0)]\\
&=\frac{1}{\sqrt{n}}\sum_{i=1}^{n}\left\{\frac{d_iy_i}{\Tilde{p}(x_i;\beta^*)}+\left(1-\frac{d_i}{\Tilde{p}(x_i;\beta^*)}\right)m_1(x_i)-\frac{(1-d_i)y_i}{1-\Tilde{p}(x_i;\beta^*)}-\left(1-\frac{1-d_i}{1-\Tilde{p}(x_i;\beta^*)}\right)m_0(x_i)-(\theta_1-\theta_0)\right\}\\
&\quad+\frac{1}{\sqrt{n}}\sum_{i=1}^n\left\{\frac{\mathbb{E}\left[1-\frac{p(X)}{\Tilde{p}(X;\beta^*)}\Big|\alpha^T_1x_i\right]}{q_1(\alpha^T_1x_i)}d_i[y_i-m_1(x_i)]-\frac{\mathbb{E}\left[1-\frac{1-p(X)}{1-\Tilde{p}(X;\beta^*)}\Big|\alpha^T_0x_i\right]}{q_0(\alpha^T_0x_i)}(1-d_i)[y_i-m_0(x_i)]\right\}\\
&\quad+o_p(1).
\end{split}
\end{equation*}

Recall the definition of locally misspecified $PS$ model in (\ref{localdef}), we have $p(x)=\Tilde{p}(x;\beta_0)(1+\delta \times s(x))$. If the $PS$ model is locally misspecified, as shown in Section 6.2.1, we can observe that the second term can be written as $O(\delta)$. It converges to 0 as $\delta\to0$. Let $\Phi(x_i,y_i,d_i):=\frac{d_iy_i}{\Tilde{p}(x_i;\beta^*)}+\left(1-\frac{d_i}{\Tilde{p}(x_i;\beta^*)}\right)m_1(x_i)-\frac{(1-d_i)y_i}{1-\Tilde{p}(x_i;\beta^*)}-\left(1-\frac{1-d_i}{1-\Tilde{p}(x_i;\beta^*)}\right)m_0(x_i)-(\theta_1-\theta_0)$. Note that $\mathbb{E}\left\{\Phi(X,Y,D)\right\}=0$. We further assume $\mathbb{E}\left\{\Phi(X,Y,D)^2\right\} < \infty$. It follows from the Central Limit Theorem and Slutsky's Theorem that the first term converges in distribution to $N(0, \Sigma)$ with $$\Sigma=\mathrm{Var}\left[\frac{DY}{\Tilde{p}(X;\beta^*)}+\left(1-\frac{D}{\Tilde{p}(X;\beta^*)}\right)m_1(X)-\frac{(1-D)Y}{1-\Tilde{p}(X;\beta^*)} -\left(1-\frac{1-D}{1-\Tilde{p}(X;\beta^*)}\right)m_0(X)-(\theta_1-\theta_0)\right].$$ Note that $\Sigma$ converges to $\Sigma_1$ as $\delta\to0$. Consequently, the asymptotic variance of $\hat{\Delta}_6$ converges to $\Sigma_1$.

If $PS$ model is globally misspecified, we can further write $\sqrt{n}(\hat{\Delta}_6-\Delta)$ as
\begin{equation*}
\begin{split}
&\sqrt{n}(\hat{\Delta}_6-\Delta)\\
&=\sqrt{n}[(\hat{\theta}_1-\hat{\theta}_0)-(\theta_1-\theta_0)]\\
&=\frac{1}{\sqrt{n}}\sum_{i=1}^{n}\bigg\{\frac{d_iy_i}{\Tilde{p}(x_i;\beta^*)}+\left(1-\frac{d_i}{\Tilde{p}(x_i;\beta^*)}\right)m_1(x_i)-\frac{(1-d_i)y_i}{1-\Tilde{p}(x_i;\beta^*)}-\left(1-\frac{1-d_i}{1-\Tilde{p}(x_i;\beta^*)}\right)m_0(x_i)-(\theta_1-\theta_0)\\
&\qquad+w_1(x_i)d_i[y_i-m_1(x_i)]+w_0(x_i)(1-d_i)[y_i-m_0(x_i)]\bigg\}+o_p(1)\\
&=\frac{1}{\sqrt{n}}\sum_{i=1}^{n}\Phi(x_i,y_i,d_i)+o_p(1),
\end{split}
\end{equation*}
\noindent where $w_1(x_i)=\frac{\mathbb{E}\left[1-\frac{p(X)}{\Tilde{p}(X;\beta^*)}\big|\alpha^T_1x_i\right]}{q_1(\alpha^T_1x_i)}$ and $w_0(x_i)=\frac{\mathbb{E}\left[\frac{1-p(X)}{1-\Tilde{p}(X;\beta^*)}-1\big|\alpha^T_0x_i\right]}{q_0(\alpha^T_0x_i)}$.

\noindent Note that $\mathbb{E}\left\{\Phi(X,Y,D)\right\}=0$. We further assume $\mathbb{E}\left\{\Phi(X,Y,D)^2\right\} < \infty$. It follows from the Central Limit Theorem and Slutsky's Theorem that $\sqrt{n}(\hat{\Delta}_6-\Delta)$ converges in distribution to $N(0, \Sigma)$ with
\begin{equation*}
\begin{split}
\Sigma &=\mathbb{E}\left\{\Phi(X,Y,D)^2\right\}\\
&=\Sigma_1+\mathbb{E}\left\{\frac{1}{p(X)}\mathrm{Var}[Y(1)|X]\left[\left(\frac{p(X)}{\Tilde{p}(X;\beta^*)}+w_1(X)p(X)\right)^2-1\right]\right\}\\
&\quad+\mathbb{E}\left\{\frac{1}{1-p(X)}\mathrm{Var}[Y(0)|X]\left[\left(\frac{1-p(X)}{1-\Tilde{p}(X;\beta^*)}+w_0(X)(1-p(X))\right)^2-1\right]\right\}.
\end{split}
\end{equation*}
\noindent Therefore, the asymptotic variance of $\hat{\Delta}_6$ is not necessarily enlarged compared to the semiparametric efficiency bound $\Sigma_1$ when $PS$ is globally misspecified.

\subsubsection{Semiparametric $PS$ model and nonparametric $OR$ model}

When $PS$ model is semiparametric and $OR$ model is nonparametric, we have
\begin{eqnarray}\label{6.2.7}
\sqrt{n}(\hat{\theta}_1-\theta_1) &=& \sqrt{n}\left\{n^{-1}\sum_{i=1}^{n}\left[\frac{d_iy_i}{\hat{g}(\alpha^Tx_i)}+\left(1-\frac{d_i}{\hat{g}(\alpha^Tx_i)}\right)\hat{m}_1(x_i)\right]-\theta_1\right\} \nonumber\\
&=&\frac{1}{\sqrt{n}}\sum_{i=1}^{n}\left\{\frac{d_iy_i}{g(\alpha^Tx_i)}+\left(1-\frac{d_i}{g(\alpha^Tx_i)}\right)m_1(x_i)-\theta_1\right\} \nonumber\\
&&+\frac{1}{\sqrt{n}}\sum_{i=1}^{n}d_i\left(\frac{1}{\hat{g}(\alpha^Tx_i)}-\frac{1}{g(\alpha^Tx_i)}\right)[y_i-m_1(x_i)]\nonumber \\
&&+\frac{1}{\sqrt{n}}\sum_{i=1}^{n}\left(1-\frac{d_i}{\hat{g}(\alpha^Tx_i)}\right)[\hat{m}_1(x_i)-m_1(x_i)] \nonumber\\
&:=& \sqrt{n}A_n + \sqrt{n}B_n + \sqrt{n}C_n.
\end{eqnarray}

\noindent Under Assumption \ref{semiassumption}, due to the consistency of $\hat{g}(\alpha^Tx_i)$ to $g(\alpha^Tx_i)$ mentioned in Section 6.2.5, we have $B_n=B^*_n+o_p(n^{-1/2})$, $C_n=C^*_n+o_p(n^{-1/2})$ with
\begin{equation*}
\begin{split}
&B^*_n=n^{-1}\sum_{i=1}^{n}d_i\frac{1}{g^2(\alpha^Tx_i)}[g(\alpha^Tx_i)-\hat{g}(\alpha^Tx_i)][y_i-m_1(x_i)],\\
&C^*_n=n^{-1}\sum_{i=1}^{n}\left(1-\frac{d_i}{g(\alpha^Tx_i)}\right)[\hat{m}_1(x_i)-m_1(x_i)].
\end{split}
\end{equation*}

\noindent Similar to the derivations of $\sqrt{n}B_n$ in Section 6.2.5, we have $\sqrt{n}B_n=o_p(1)$. Similar to the derivations of $\sqrt{n}C_n$ in Section 6.2.2, we have $\sqrt{n}C_n=\frac{1}{\sqrt{n}}\sum_{j=1}^n\Big\{\big(\frac{1}{p(x_j)}-\frac{1}{g(\alpha^Tx_j)}\big)d_j[y_j-m_1(x_j)]\Big\}+o_p(1)$. Consequently, combining the terms in (\ref{6.2.7}), we have
\begin{equation*}
\begin{split}
&\sqrt{n}(\hat{\theta}_1-\theta_1)\\
&=\frac{1}{\sqrt{n}}\sum_{i=1}^{n}\left\{\frac{d_iy_i}{g(\alpha^Tx_i)}+\left(1-\frac{d_i}{g(\alpha^Tx_i)}\right)m_1(x_i)-\theta_1\right\}\\
&\quad+\frac{1}{\sqrt{n}}\sum_{i=1}^n\left\{\left(\frac{1}{p(x_i)}-\frac{1}{g(\alpha^Tx_i)}\right)d_i[y_i-m_1(x_i)]\right\}\\
&\quad+o_p(1)\\
&=\frac{1}{\sqrt{n}}\sum_{i=1}^{n}\left\{\frac{d_iy_i}{p(x_i)}+\left(1-\frac{d_i}{p(x_i)}\right)m_1(x_i)-\theta_1\right\}+o_p(1).
\end{split}
\end{equation*}
Similarly, we can derive the form of $\sqrt{n}(\hat{\theta}_0-\theta_0)$. As a result,
\begin{equation*}
\begin{split}
&\sqrt{n}(\hat{\Delta}_7-\Delta)\\
&=\sqrt{n}[(\hat{\theta}_1-\hat{\theta}_0)-(\theta_1-\theta_0)]\\
&=\frac{1}{\sqrt{n}}\sum_{i=1}^{n}\left\{\frac{d_iy_i}{p(x_i)}+\left(1-\frac{d_i}{p(x_i)}\right)m_1(x_i)-\frac{(1-d_i)y_i}{1-p(x_i)}-\left(1-\frac{1-d_i}{1-p(x_i)}\right)m_0(x_i)-(\theta_1-\theta_0)\right\}+o_p(1)\\
&=\frac{1}{\sqrt{n}}\sum_{i=1}^{n}\Phi(x_i,y_i,d_i)+o_p(1).
\end{split}
\end{equation*}

\noindent Note that $\mathbb{E}\left\{\Phi(X,Y,D)\right\}=0$. We further assume $\mathbb{E}\left\{\Phi(X,Y,D)^2\right\} < \infty$. It follows from the Central Limit Theorem and Slutsky's Theorem that $\sqrt{n}(\hat{\Delta}_7-\Delta)$ converges in distribution to $N(0,\Sigma_1)$. In other words, the asymptotic variance of $\hat{\Delta}_7$ achieves the semiparametric efficiency bound $\Sigma_1$ no matter whether $PS$ model is correctly specified or not.

\subsubsection{Nonparametric $PS$ model and semiparametric $OR$ model}

When $PS$ model is nonparametric and $OR$ model is semiparametric, we have
\begin{eqnarray}\label{6.2.8}
\sqrt{n}(\hat{\theta}_1-\theta_1)&=&\sqrt{n}\left\{n^{-1}\sum_{i=1}^{n}\left[\frac{d_iy_i}{\hat{p}(x_i)}+\left(1-\frac{d_i}{\hat{p}(x_i)}\right)\hat{r}_1(\alpha^T_1x_i)\right]-\theta_1\right\} \nonumber\\
&=&\frac{1}{\sqrt{n}}\sum_{i=1}^{n}\left\{\frac{d_iy_i}{p(x_i)}+\left(1-\frac{d_i}{p(x_i)}\right)r_1(\alpha^T_1x_i)-\theta_1\right\} \nonumber\\
&&+\frac{1}{\sqrt{n}}\sum_{i=1}^{n}d_i\left(\frac{1}{\hat{p}(x_i)}-\frac{1}{p(x_i)}\right)[y_i-r_1(\alpha^T_1x_i)]\nonumber \\
&&+\frac{1}{\sqrt{n}}\sum_{i=1}^{n}\left(1-\frac{d_i}{\hat{p}(x_i)}\right)[\hat{r}_1(\alpha^T_1x_i)-r_1(\alpha^T_1x_i)]\nonumber\\
&:=& \sqrt{n}A_n + \sqrt{n}B_n + \sqrt{n}C_n.
\end{eqnarray}

\noindent Under Assumption \ref{nonparaassumption}, due to the consistency of $\hat{p}(x)$ to $p(x)$ mentioned in Section 6.2.2, we have $B_n=B^*_n+o_p(n^{-1/2})$, $C_n=C^*_n+o_p(n^{-1/2})$ with
\begin{equation*}
\begin{split}
&B^*_n=n^{-1}\sum_{i=1}^{n}d_i\frac{1}{p^2(x_i)}[p(x_i)-\hat{p}(x_i)][y_i-r_1(\alpha^T_1x_i)],\\
&C^*_n=n^{-1}\sum_{i=1}^{n}\left(1-\frac{d_i}{p(x_i)}\right)[\hat{r}_1(\alpha^T_1x_i)-r_1(\alpha^T_1x_i)].
\end{split}
\end{equation*}

\noindent Similar to the derivations of $\sqrt{n}B_n$ in Section 6.2.3, we have $\sqrt{n}B_n=\frac{1}{\sqrt{n}}\sum_{j=1}^n[m_1(x_j)-r_1(\alpha^T_1x_j)]\big(1-\frac{d_j}{p(x_j)}\big)+o_p(1)$. Similar to the derivations of $\sqrt{n}C_n$ in Section 6.2.6, we have $\sqrt{n}C_n=o_p(1)$. Combining the terms in (\ref{6.2.8}), we have
\begin{equation*}
\begin{split}
\sqrt{n}(\hat{\theta}_1-\theta_1)&=\frac{1}{\sqrt{n}}\sum_{i=1}^{n}\left\{\frac{d_iy_i}{p(x_i)}+\left(1-\frac{d_i}{p(x_i)}\right)r_1(\alpha^T_1x_i)-\theta_1\right\}\\
&\quad+\frac{1}{\sqrt{n}}\sum_{i=1}^n[m_1(x_i)-r_1(\alpha^T_1x_i)]\big(1-\frac{d_i}{p(x_i)}\big)\\
&\quad+o_p(1)\\
&=\frac{1}{\sqrt{n}}\sum_{i=1}^{n}\left\{\frac{d_iy_i}{p(x_i)}+\left(1-\frac{d_i}{p(x_i)}\right)m_1(x_j)-\theta_1\right\}+o_p(1).
\end{split}
\end{equation*}
\noindent Similarly, we can derive the form of $\sqrt{n}(\hat{\theta}_0-\theta_0)$. As a result,
\begin{equation*}
\begin{split}
&\sqrt{n}(\hat{\Delta}_8-\Delta) \\
&=\sqrt{n}[(\hat{\theta}_1-\hat{\theta}_0)-(\theta_1-\theta_0)]\\
&=\frac{1}{\sqrt{n}}\sum_{i=1}^{n}\left\{\frac{d_iy_i}{p(x_i)}+\left(1-\frac{d_i}{p(x_i)}\right)m_1(x_j)-\frac{(1-d_i)y_i}{1-p(x_i)}-\left(1-\frac{1-d_i}{1-p(x_i)}\right)m_0(x_j)-(\theta_1-\theta_0)\right\}+o_p(1)\\
&=\frac{1}{\sqrt{n}}\sum_{i=1}^{n}\Phi(x_i,y_i,d_i)+o_p(1).
\end{split}
\end{equation*}
\noindent Note that $\mathbb{E}\left\{\Phi(X,Y,D)\right\}=0$. We further assume $\mathbb{E}\left\{\Phi(X,Y,D)^2\right\} < \infty$. It follows from the Central Limit Theorem and Slutsky's Theorem that $\sqrt{n}(\hat{\Delta}_8-\Delta)$ converges in distribution to $N(0,\Sigma_1)$. In other words, the asymptotic variance of $\hat{\Delta}_8$ achieves the semiparametric efficiency bound $\Sigma_1$ no matter whether $OR$ model is correctly specified or not.

\subsubsection{Semiparametric $PS$ model and $OR$ model}

When $PS$ model and $OR$ model are both semiparametric, we have
\begin{eqnarray}\label{6.2.9}
\sqrt{n}(\hat{\theta}_1-\theta_1)
&=&\sqrt{n}\left\{n^{-1}\sum_{i=1}^{n}\left[\frac{d_iy_i}{\hat{g}(\alpha^Tx_i)}+\left(1-\frac{d_i}{\hat{g}(\alpha^Tx_i)}\right)\hat{r}_1(\alpha^T_1x_i)\right]-\theta_1\right\}\nonumber\\
&=&\frac{1}{\sqrt{n}}\sum_{i=1}^{n}\left\{\frac{d_iy_i}{g(\alpha^Tx_i)}+\left(1-\frac{d_i}{g(\alpha^Tx_i)}\right)r_1(\alpha^T_1x_i)-\theta_1\right\} \nonumber\\
&&+\frac{1}{\sqrt{n}}\sum_{i=1}^{n}d_i\left(\frac{1}{\hat{g}(\alpha^Tx_i)}-\frac{1}{g(\alpha^Tx_i)}\right)[y_i-r_1(\alpha^T_1x_i)] \nonumber\\
&&+\frac{1}{\sqrt{n}}\sum_{i=1}^{n}\left(1-\frac{d_i}{\hat{g}(\alpha^Tx_i)}\right)[\hat{r}_1(\alpha^T_1x_i)-r_1(\alpha^T_1x_i)] \nonumber\\
&:=& \sqrt{n}A_n + \sqrt{n}B_n + \sqrt{n}C_n.
\end{eqnarray}

\newcontent Under Assumption \ref{semiassumption}, due to the consistency of $\hat{g}(\alpha^Tx)$ to $g(\alpha^Tx)$ mentioned in Section 6.2.5, we have $B_n=B^*_n+o_p(n^{-1/2})$, $C_n=C^*_n+o_p(n^{-1/2})$ with
\begin{equation*}
\begin{split}
&B^*_n=n^{-1}\sum_{i=1}^{n}d_i\frac{1}{g^2(\alpha^Tx_i)}[g(\alpha^Tx_i)-\hat{g}(\alpha^Tx_i)][y_i-r_1(\alpha^T_1x_i)],\\
&C^*_n=n^{-1}\sum_{i=1}^{n}\left(1-\frac{d_i}{g(\alpha^Tx_i)}\right)[\hat{r}_1(\alpha^T_1x_i)-r_1(\alpha^T_1x_i)].
\end{split}
\end{equation*}

\noindent Similar to the derivations of $\sqrt{n}B_n$ in Section 6.2.5, we have
$$\sqrt{n}B_n=\frac{1}{\sqrt{n}}\sum_{j=1}^n\frac{E\big\{p(X)[m_1(X)-r_1(\alpha^T_1X)]|\alpha^Tx_j\big\}}{g^2(\alpha^Tx_j)}[g(\alpha^Tx_j)-d_j]+o_p(1).$$
\noindent Similar to the derivations of $\sqrt{n}C_n$ in Section 6.2.6, we have
$$\sqrt{n}C_n=\frac{1}{\sqrt{n}}\sum_{j=1}^n\frac{E\big[1-\frac{p(X)}{g(\alpha^TX)}|\alpha^T_1x_j\big]}{q_1(\alpha^T_1x_j)}d_j[y_j-r_1(\alpha^T_1x_j)]+o_p(1).$$
\noindent Now we can consider different cases as follows.

\newcontent (a) \textbf{Correctly specified $PS$ and $OR$ models}
\newcontent

In this case, $p(x)=g(\alpha^Tx)$, $m_1(x)=r_1(\alpha^T_1x)$ and $m_0(x)=r_0(\alpha^T_0x)$. So we have $\sqrt{n}B_n=o_p(1)$ and $\sqrt{n}C_n=o_p(1)$. Combining the terms in (\ref{6.2.9}), we have
\begin{equation*}
\begin{split}
\sqrt{n}(\hat{\theta}_1-\theta_1)=\frac{1}{\sqrt{n}}\sum_{i=1}^{n}\left\{\frac{d_iy_i}{p(x_i)}+\left(1-\frac{d_i}{p(x_i)}\right)m_1(x_i)-\theta_1\right\}+o_p(1).
\end{split}
\end{equation*}
Similarly, we can derive the form of $\sqrt{n}(\hat{\theta}_0-\theta_0)$. As a result,
\begin{equation*}
\begin{split}
&\sqrt{n}(\hat{\Delta}_9-\Delta)\\
&=\sqrt{n}[(\hat{\theta}_1-\hat{\theta}_0)-(\theta_1-\theta_0)]\\
&=\frac{1}{\sqrt{n}}\sum_{i=1}^{n}\left\{\frac{d_iy_i}{p(x_i)}+\left(1-\frac{d_i}{p(x_i)}\right)m_1(x_i)-\frac{(1-d_i)y_i}{1-p(x_i)}-\left(1-\frac{1-d_i}{1-p(x_i)}\right)m_0(x_i)-(\theta_1-\theta_0)\right\}+o_p(1)\\
&=\frac{1}{\sqrt{n}}\sum_{i=1}^{n}\Phi(x_i,y_i,d_i)+o_p(1).
\end{split}
\end{equation*}

\noindent Note that $\mathbb{E}\left\{\Phi(X,Y,D)\right\}=0$. We further assume $\mathbb{E}\left\{\Phi(X,Y,D)^2\right\} < \infty$. It follows from the Central Limit Theorem and Slutsky's Theorem that $\sqrt{n}(\hat{\Delta}_9-\Delta)$ converges in distribution to $N(0,\Sigma_1)$. The asymptotic variance of $\hat{\Delta}_9$ achieves the semiparametric efficiency bound $\Sigma_1$ when both $PS$ and $OR$ models are correctly specified.

\newcontent (b) \textbf{Correctly specified $PS$ model and misspecified $OR$ model}

In this case, $p(x)=g(\alpha^Tx)$, $m_1(x)\neq r_1(\alpha^T_1x)$ and $m_0(x)\neq r_0(\alpha^T_0x)$. So we have $\sqrt{n}C_n=o_p(1)$. Combining the terms in (\ref{6.2.9}), we have
\begin{equation*}
\begin{split}
&\sqrt{n}(\hat{\theta}_1-\theta_1)\\
&=\frac{1}{\sqrt{n}}\sum_{i=1}^{n}\left\{\frac{d_iy_i}{p(x_i)}+\left(1-\frac{d_i}{p(x_i)}\right)r_1(\alpha^T_1x_i)-\theta_1\right\}\\
&\quad+\frac{1}{\sqrt{n}}\sum_{i=1}^n\frac{\mathbb{E}\left\{p(X)[m_1(X)-r_1(\alpha^T_1X)]|\alpha^Tx_i\right\}}{p^2(x_i)}[p(x_i)-d_i]\\
&\quad+o_p(1).
\end{split}
\end{equation*}
Similarly, we can derive the form of $\sqrt{n}(\hat{\theta}_0-\theta_0)$. As a result,
\begin{equation*}
\begin{split}
&\sqrt{n}(\hat{\Delta}_9-\Delta)\\
&=\sqrt{n}[(\hat{\theta}_1-\hat{\theta}_0)-(\theta_1-\theta_0)]\\
&=\frac{1}{\sqrt{n}}\sum_{i=1}^{n}\bigg\{\frac{d_iy_i}{p(x_i)}+\left(1-\frac{d_i}{p(x_i)}\right)r_1(\alpha^T_1x_i)-\frac{(1-d_i)y_i}{1-p(x_i)}-\left(1-\frac{1-d_i}{1-p(x_i)}\right)r_0(\alpha^T_0x_i)-(\theta_1-\theta_0)\\
&\qquad+w(x_i)[p(x_i)-d_i]\bigg\}+o_p(1)\\
&=\frac{1}{\sqrt{n}}\sum_{i=1}^{n}\Phi(x_i,y_i,d_i)+o_p(1),
\end{split}
\end{equation*}
where $w(x_i)=\frac{\mathbb{E}\left\{p(X)[m_1(X)-r_1(\alpha^T_1X)]|\alpha^Tx_i\right\}}{p^2(x_i)}+\frac{\mathbb{E}\left\{(1-p(X))[m_0(X)-r_0(\alpha^T_0X)]|\alpha^Tx_i\right\}}{(1-p(x_i))^2}.$

\noindent Note that $\mathbb{E}\left\{\Phi(X,Y,D)\right\}=0$. We further assume $\mathbb{E}\left\{\Phi(X,Y,D)^2\right\} < \infty$. It follows from the Central Limit Theorem and Slutsky's Theorem that $\sqrt{n}(\hat{\Delta}_9-\Delta)$ converges in distribution to $N(0,\Sigma)$ with
\begin{equation*}
\begin{split}
\Sigma &=\mathbb{E}\left\{\Phi(X,Y,D)^2\right\}\\
&=\Sigma_1+\mathbb{E}\bigg\{\Big[\sqrt{\frac{1}{p(X)}-1}[r_1(\alpha^T_1X)-m_1(X)]+\sqrt{\frac{1}{1-p(X)}-1}[r_0(\alpha^T_0X)-m_0(X)]\\
&\quad+\sqrt{p(X)(1-p(X))}w(X)\Big]^2\bigg\} \geq \Sigma_1.
\end{split}
\end{equation*}
\noindent Therefore, the asymptotic variance of $\hat{\Delta}_9$ is enlarged compared to the semiparametric efficiency bound $\Sigma_1$ when $OR$ models are misspecified.

\newcontent (c) \textbf{Misspecified $PS$ model and correctly specified $OR$ model}

In this case, $p(x)\neq g(\alpha^Tx)$, $m_1(x)=r_1(\alpha^T_1x)$ and $m_0(x)=r_1(\alpha^T_0x)$. So $\sqrt{n}B_n=o_p(1)$.  Combining the terms in (\ref{6.2.9}), we have
\begin{equation*}
\begin{split}
&\sqrt{n}(\hat{\theta}_1-\theta_1)\\
&=\frac{1}{\sqrt{n}}\sum_{i=1}^{n}\left\{\frac{d_iy_i}{g(\alpha^Tx_i)}+\left(1-\frac{d_i}{g(\alpha^Tx_i)}\right)m_1(x_i)-\theta_1\right\}\\
&\quad+\frac{1}{\sqrt{n}}\sum_{i=1}^n\frac{\mathbb{E}\left[1-\frac{p(X)}{g(\alpha^TX)}\Big|\alpha^T_1x_i\right]}{q_1(\alpha^T_1x_i)}d_i[y_i-m_1(x_i)]\\
&\quad+o_p(1).
\end{split}
\end{equation*}
Similarly, we can derive the form of $\sqrt{n}(\hat{\theta}_0-\theta_0)$. As a result,
\begin{equation*}
\begin{split}
&\sqrt{n}(\hat{\Delta}_9-\Delta)\\
&=\frac{1}{\sqrt{n}}\sum_{i=1}^{n}\bigg\{\frac{d_iy_i}{g(\alpha^Tx_i)}+\left(1-\frac{d_i}{g(\alpha^Tx_i)}\right)m_1(x_i)-\frac{(1-d_i)y_i}{1-g(\alpha^Tx_i)}-\left(1-\frac{1-d_i}{1-g(\alpha^Tx_i)}\right)m_0(x_i)-(\theta_1-\theta_0)\\
&\qquad+w_1(x_i)d_i[y_i-m_1(x_i)]+w_0(x_i)(1-d_i)[y_i-m_0(x_i)]\bigg\}+o_p(1)\\
&=\frac{1}{\sqrt{n}}\sum_{i=1}^{n}\Phi(x_i,y_i,d_i)+o_p(1),
\end{split}
\end{equation*}
where $w_1(x_i)=\frac{\mathbb{E}\left[1-\frac{p(X)}{g(\alpha^TX)}\big|\alpha^T_1x_i\right]}{q_1(\alpha^T_1x_i)}$ and $w_0(x_i)=\frac{\mathbb{E}\left[\frac{1-p(X)}{1-g(\alpha^TX)}-1\big|\alpha^T_0x_i\right]}{q_0(\alpha^T_0x_i)}$.

\noindent Note that $\mathbb{E}\left\{\Phi(X,Y,D)\right\}=0$. We further assume $\mathbb{E}\left\{\Phi(X,Y,D)^2\right\} < \infty$. It follows from the Central Limit Theorem and Slutsky's Theorem that $\sqrt{n}(\hat{\Delta}_9-\Delta)$ converges in distribution to $N(0, \Sigma)$ with
\begin{equation*}
\begin{split}
\Sigma &=\mathbb{E}\left\{\Phi(X,Y,D)^2\right\}\\
&=\Sigma_1+\mathbb{E}\left\{\frac{1}{p(X)}Var[Y(1)|X]\left[\left(\frac{p(X)}{g(\alpha^TX)}+w_1(X)p(X)\right)^2-1\right]\right\}\\
&\quad+\mathbb{E}\left\{\frac{1}{1-p(X)}Var[Y(0)|X]\left[\left(\frac{1-p(X)}{1-g(\alpha^TX)}+w_0(X)(1-p(X))\right)^2-1\right]\right\}.
\end{split}
\end{equation*}

\noindent Therefore, the asymptotic variance of $\hat{\Delta}_9$ is not necessarily enlarged compared to the semiparametric efficiency bound $\Sigma_1$ when $PS$ model is misspecified.


\begin{thebibliography}{9}
\bibitem{nonparametric-reference1}
Abrevaya, J., Hsu, Y.-C. and Lieli, R. P. (2015). Estimating Conditional Average Treatment Effects.
\textit{Journal of Business and Economic Statistics}, 33(4), 485–505.

\bibitem{mle1}
Davison, A. C. (2009). {\it Statistical models}. Cambridge Univ. Press.

\bibitem{simulation1}
Ding, X. and Wang, Q. (2011). Fusion-Refinement Procedure for Dimension Reduction With Missing Response at Random.
\textit{Journal of the American Statistical Association}, 106(495), 1193–1207.

\bibitem{drproof}
Funk, M. J., Westreich, D., Wiesen, C., Stürmer, T., Brookhart, M. A. and Davidian, M. (2011). Doubly Robust Estimation of Causal Effects.
\textit{American Journal of Epidemiology}, 173(7), 761–767.

\bibitem{semidr}
Guo, X., Fang, Y., Zhu, X., Xu, W. and Zhu, L. (2018). Semiparametric double robust and efficient estimation for mean functionals with response missing at random.
\textit{Computational Statistics and Data Analysis}, 128, 325-339.

\bibitem{misspecification-reference1}
Han, F. (2018). Doubly robust estimation of the causal effects in the causal inference with missing outcome data.
\textit{Journal of Ambient Intelligence and Humanized Computing}.

\bibitem{efficiencybound}
Hahn, J. (1998). On the Role of the Propensity Score in Efficient Semiparametric Estimation of Average Treatment Effects.
\textit{Econometrica}, 66(2), 315-331.

\bibitem{psnonpara}
Hirano, K., G. Imbens and G. Ridder. (2003). Effcient Estimation of Average
Treatment Effects Using the Estimated Propensity Score.
\textit{Econometrica}, 71(4): 1161-1189.

\bibitem{ornonpara}
Imbens, G. W., Newey, W. K. and Ridder, G. (2005). Mean-square-error Calculations for Average Treatment Effects. IEPR Working Paper No. 05.34

\bibitem{nonparametric-reference2}
Ichimura, H. and Linton, O. (2005). Asymptotic Expansions for Some Semiparametric Program Evaluation Estimators.
\textit{Identification and Inference for Econometric Models}, 149–170.

\bibitem{simulation2}
Kang, J. D. Y. and Schafer, J. L. (2007). Rejoinder: Demystifying Double Robustness: A Comparison of Alternative Strategies for Estimating a Population Mean from Incomplete Data.
\textit{Statistical Science}, 22(4), 574–580.

\bibitem{drpara+pspara}
Lunceford, J. K. and Davidian, M. (2004). Stratification and weighting via the propensity score in estimation of causal treatment effects: a comparative study.
\textit{Statistics in Medicine}, 23(19), 2937–2960.

\bibitem{misspecification-reference2}
Lefebvre, G. and Gustafson, P. (2010). Impact of Outcome Model Misspecification on Regression and Doubly-Robust Inverse Probability Weighting to Estimate Causal Effect.
\textit{The International Journal of Biostatistics}, 6(2), Article 15.

\bibitem{nonparametric-reference3}
Masry, E. (1996), Multivariate Local Polynomial Regression for Time Series: Uniform Strong Consistency and Rates.
\textit{Journal of Time Series Analysis}, 17, 571-599.

\bibitem{framework1}
Splawa-Neyman, J., Dabrowska, D. M. and Speed, T. P. (1990). On the Application of Probability Theory to Agricultural Experiments. Essay on Principles. Section 9.
\textit{Statistical Science}, 5(4), 465–472.

\bibitem{Ustatistics}
Powell, J. L., Stock, J. H. and Stoker, T. M. (1989). Semiparametric Estimation of Index Coefficients. \textit{Econometrica}, 57(6), 1403-1430.

\bibitem{framework2}
Rubin, D. B. (1974). Estimating causal effects of treatments in randomized and nonrandomized studies.
\textit{Journal of Educational Psychology}, 66(5), 688–701.

\bibitem{pspara}
Rubin, D. (1977). Assignment to Treatment Group on the Basis of a Covariate.
\textit{Journal of Educational Statistics}, 2(1), 1-26.

\bibitem{dr}
Robins, J., Rotnitzky, A. and Zhao, L. (1994). Estimation of Regression Coefficients When Some Regressors Are Not Always Observed.
\textit{Journal of the American Statistical Association}, 89(427), 846-866.

\bibitem{denonpara}
Rothe, C. and Firpo, S. (2013). Semiparametric Estimation and Inference Using Doubly Robust Moment Conditions. IDEAS Working Paper Series from RePEc, 2013.

\bibitem{comment}
Tan, Z. (2007). Comment: Understanding or, PS and DR.
\textit{Statistical Science}, 22(4), 560-568.

\bibitem{mle2}
White, H. (1982). Maximum Likelihood Estimation of Misspecified Models.
\textit{Econometrica}, 50(1), 1-25.

\bibitem{mvt}
Wade, W. R. (2018). {\it An introduction to analysis}. New York, NY: Pearson.

\end{thebibliography}
\end{document}